\documentclass[final]{siamltex}

\def \bfw {\mbox{\boldmath$  \mathcal{W}$}}
\def \bfell {\mbox{\boldmath$\ell$}}
\def \bfr {\mbox{\boldmath$r$}}

\def \bftf {\mbox{\boldmath $\tilde {\mathcal F} $}}
\usepackage{subeqn}
\usepackage {subfig}

\usepackage{graphicx}
\usepackage{epstopdf}

\title{Variable relaxed schemes for multidimensional hyperbolic conservation laws}

\author{Shalini Krishnamurthy\thanks{Institute for Computational \&
Mathematical Engineering, Stanford University, Stanford, CA 94305 ({\tt
shalinik@alumni.stanford.edu}).} \and Margot Gerritsen\thanks{Department of Energy Resources
Engineering, Stanford University, Stanford, CA 94305 ({\tt
margot.gerritsen@stanford.edu}).}}

\begin{document}

\DeclareGraphicsExtensions{.pdf,.gif,.jpg}

\maketitle


\begin{abstract}
We present a new class of component-wise numerical schemes that are in the family of relaxation formulations, originally
introduced by [S. Jin and Z. P. Xin, Comm. Pure Appl. Math., 48(1995), pp. 235-277]. The relaxation framework enables the
construction of schemes that are free of nonlinear Riemann solvers and are independent of the underlying eigenstructure of
the problem. The constant relaxation schemes proposed by Jin \& Xin  can however introduce strong numerical diffusion,
especially when the maximum characteristic speeds are high compared to the average speeds in the domain. We propose a
general class of  variable relaxation formulations for multidimensional systems of conservation laws which utilizes
estimates of local maximum and  minimum speeds to arrive at more accurate relaxation schemes, irrespective of the contrast
in maximum and average characteristic speeds. First and second order variable relaxation methods are presented for general
nonlinear systems  in one and two spatial dimensions, along with monotonicity and TVD (Total Variation Diminishing)
properties for the 1D schemes. The effectiveness of the schemes is demonstrated on a test suite that includes Burgers'
equation, the weakly hyperbolic Engquist-Runborg problem, as well as the weakly hyperbolic gas injection displacements
that are governed by strong nonlinear coupling thus making them highly sensitive to numerical diffusion.
In the latter examples the second order Jin-Xin scheme fails to capture the fronts reasonably, when both the first
and second order variable relaxed schemes produce the displacement profiles sharply.
\end{abstract}



\pagestyle{myheadings}
\thispagestyle{plain}
\markboth{S. KRISHNAMURTHY AND M. GERRITSEN}{Variable relaxed schemes for multidimensional conservation laws}


\section{Introduction}
Nonlinear hyperbolic conservation laws are often solved with schemes that are based on upwind differencing motivated by
Courant-Isaacson-Rees \cite{cir} and used within Godunov's reconstruction-evolution-average (REA) framework \cite{god}.
Upwind schemes require exact or approximate solution of a nonlinear Riemann problem at the cell interfaces, which in
turn necessitates decomposition of the Riemann fan, to determine the contribution from each characteristic variable. This
traditional application of upwind schemes is not possible when Riemann solutions are not available or when  the system does
not have a full set of eigenvectors.  One such example is the two-phase gas injection displacement in subsurface
formations \cite{orr,gd} which motivated this work and are discussed in detail in section \ref{gasinjection}.
The governing systems of equations are characterized \cite{kg,gd} by weak
hyperbolicity at isolated points in space, strong nonlinear coupling, and flux evaluations that require computationally
expensive thermal equilibrium calculations. Furthermore, Riemann solutions
for these  problems are available only for simplified phase behavior. Component-wise upwinding can be used, but only when gravitational effects can be
ignored so that all wave speeds have the same sign. For general problems, which are characterized by both  positive as well as negative
characteristic speeds, component-wise limiting gives rise to spurious oscillations \cite{qs, heoc}.
This motivated us to investigate other families of schemes for these interesting and challenging problems that are independent of the
eigenstructure of the system and do not require Riemann solvers.

Central schemes \cite{nt,kt,knp, jllot} are one class of numerical methods  which are independent of the
eigenstructure of the system and do not require Riemann solvers. The fully-discrete higher order Nessyahu-Tadmor central
scheme (NT scheme) \cite{nt}  which is a natural extension to Lax-Freidrichs (LxF) scheme, is constructed by viewing the
LxF scheme as a Godunov approach where Riemann solvers and characteristic decompositions are avoided by evolving
staggered cell averages. This approach requires alternating between two staggered grids which becomes particularly
combersome near domain boundaries. A procedure to avoid staggering, while still retaining the resolution of the
NT scheme was introduced in  \cite{jllot}.
However, because these schemes, like the LxF scheme, do not take advantage of  characteristic information,
they generally introduce much higher numerical diffusion than upwind schemes of the same order.
The  higher order central schemes of Kurganov and Tadmor (KT schemes) \cite{kt}, which can viewed as the higher order
extensions of  local-Lax-Friedrichs (LLF) scheme \cite{rus},
reduce numerical diffusion significantly by using some characteristic information, i.e., the  maximum absolute,
local characteristic speeds. The central-upwind schemes \cite{knp, kl} further reduce diffusion
by including the information of local maximum and minimum characteristic speeds in each direction.
However, as explained in detail in section 6,  in
gas injection processes, the strong nonlinear coupling demands more phase equilibrium calculations in a higher order
central/central-upwind framework.
Since the phase equilibrium calculations form the computational bottleneck for these
processes, we took a different approach that leads to schemes similar to central/central-upwind schemes.

Instead of a staggered cell approach, the LxF scheme for nonlinear conservation laws can also be viewed as an update
obtained by upwinding on the characteristic variables of a nearby linear hyperbolic system with nonlinear source terms.
This is the basis for the relaxation schemes introduced by Jin \& Xin \cite{jx}.
Like the central schemes, the relaxation schemes do not require a full eigensystem or Riemann solutions of the original nonlinear
problem. They also require lesser flux evaluations than the central framework (see section 6) and hence computationally more viable to
gas injection processes. The relaxation model for nonlinear systems was first studied by  Liu \cite{liu}.
Using relaxation framework  for development  of higher order schemes was first introduced by Jin and Xin \cite{jx}.
Since then it has been extensively studied \cite{an,chal,lp,tt}. Relaxation schemes have been
tried on problems like equations of gas dynamics \cite{ban}, shallow water systems \cite{dk, ak}, and
weakly-hyperbolic, conservation-formulation of Hamilton-Jacobi equations \cite{jx98}.

Though amenable for weakly hyperbolic systems, the Jin \& Xin's schemes suffer from excessive numerical diffusion that
increases as the difference between the global maximum speed and the average speed of the problem grows. In multidimensional gas
injection problems, this may result in second order solution profiles that are more diffusive than that of a first order LLF scheme.
In \cite{kt}, Kurganov and Tadmor mention about the relationship between central and relaxation schemes, and that a special
choice of relaxation matrix ${A = \rho\left(  \frac{\partial {f}({u)}}{\partial u}\right)I }$ can lead to high-resolution schemes
that will be similar to their KT central scheme. In \cite{lp}, Leveque and Pelanti
mention that relaxation schemes can be viewed as a means for defining an approximate Riemann solver.
Motivated by the aforementioned works, in \cite{kg}, we proposed a  variable relaxation formulation for gas injection processes,
based on a nearby linear hyperbolic system with locally variable eigenvalues. In his seminal work \cite{liu} in 1987, Liu
presented the fundamentals of local subcharacteristics in a relaxation model. This forms the basis for our choice of local
relaxation parameters.

Other local relaxation approaches have been presented in \cite{hjs} and more recently in \cite{ak}.
The approach presented in \cite{hjs} yields a first order scheme that does not
obey discrete conservation (since the numerical diffusion terms do not telescope upon summation, see section 7 of \cite{hjs}),
except in the trivial case where the system has constant velocities.
This lead us to base our work on the fundamentals presented by Liu \cite{liu}.
In  \cite{ak}, which was presented around the same time as our initial work \cite{kg}, a relaxation formulation with local subcharacteristics, specifically for the
1D conditionally hyperbolic two layer shallow water system has been proposed. Here the authors use Roe linearization to solve
the  relaxation formulation. Their approach is somewhat similar to our work in \cite{kg}, but we go into more detailed analysis and understanding
of the local subcharacteristics of the relaxation schemes.

In this work we present a general class of
1D variable relaxation schemes for general nonlinear systems of conservation laws (section 3) and provide a framework for
multidimensional extensions (section 4). We also verify the stability properties - monotonicity for
the first order schemes and the TVD property for the second order schemes - in section 3.  Further,
we attempt to understand the connection between central schemes and local relaxation schemes.
In section 3, it can be seen that
particular choices of local subcharacteristic speeds in the variable relaxation formulations give rise to
schemes that are similar to KT and KNP schemes. In fact, the first order variable relaxed schemes are the same as those of
KT and KNP schemes. For the second order schemes, as noticed in our 1D experiments,
 the accuracy of the central schemes and the relaxation schemes are quite close:
Jin-Xin scheme is similar to NT scheme, the variable relaxed schemes with symmetric and optimal choice of speeds (section 3)
are similar to KT and KNP schemes respectively (see \cite{kg}). We present this work with hope that it will take us a step closer
to understanding the interesting nature of relationship between the various approaches to arrive at component-wise schemes.

In section 5, we present the results of our higher order schemes on Burgers equation (order of accuracy test),
the 2D weakly hyperbolic Engquist-Runborg problem \cite{er} and  the weakly hyperbolic gas injection displacements.
The  improved resolution of our relaxation schemes is especially noticeable in the 1D and 2D gas injection problems.
For completeness we first recapitulate, from \cite{kg}, the analysis of constant relaxation schemes in the following section.

\section{Constant relaxation in 1D}
In the relaxation approach, the  component-wise updates for the nonlinear conservation systems are obtained by
applying upwind schemes to a  linear, strongly hyperbolic system
that is close to the original nonlinear conservation system.
The nonlinearity is moved to a stiff source term.
In constant relaxation formulations, the nearby linear system also has constant eigenvalues and eigenvectors.
The Jin-Xin (JX) relaxation scheme \cite{jx} is a special case of constant relaxation, which we discuss first.
\subsection{Jin-Xin relaxation schemes}

Given a system of conservation laws,
\begin{equation}
\label{eq:eqn2.1}
\frac{\partial {\bf C}}{\partial t}
+ \frac{\partial {\bf F}({\bf C)}}{\partial x}
= 0,
\hspace{.2in}{\bf C}, {\bf F}\in R^{N},
\end{equation}
with  the initial conditions
\[ {\bf C}(x,0) = g(x) \]
the  relaxation system of Jin \& Xin replaces the original system (\ref{eq:eqn2.1}) by
\begin{eqnarray}
&\frac{\partial{\bf C}}{\partial t}  &
+\frac{\partial{\bf V}}{\partial x}
= 0 ,
\hspace{.3in}{\bf C}, {\bf F}, {\bf V} \in R^{N}
\nonumber\\[-1.5ex]
\label{eq:eqn2.2} \\[-1.5ex]
&\frac{\partial {\bf V}}{\partial t} &
+ {\bf A}^{2}\frac{\partial {\bf C}}{\partial x}
= \frac{1}{\varepsilon}({\bf F(C)}-{\bf V}),
\nonumber
\end{eqnarray}
with  the additional initial conditions ${\bf V}(x,0)={\bf F(C}(x,0))$.

Here ${\bf A}=diag\left(a_{1} ,a_{2} \cdots a_{N} \right)$ is a positive diagonal matrix. Its
diagonal elements $\left\{a_{p} \right\}$  are called the subcharacteristic speeds.
The parameter $\varepsilon $, which is called the relaxation rate, is a small positive number $0<\varepsilon <<1$.
We can use the Chapman-Enskog expansion \cite{cc} to represent ${\bf V}$ as
\begin{equation}
\label{eq:eqn2.3} {\bf V}(x,t) =
{\bf F(C}(x,t)) + \varepsilon {\bf V}_{1}(x,t) + \varepsilon ^{2} {\bf
V}_{2} (x,t) + \varepsilon ^{3} {\bf V}_{3}(x,t) + \cdots  ,
\end{equation}
since in the limit  $\varepsilon \to 0$ we have that ${\bf V}\to {\bf F\left(C\right)}$ and the relaxation system  (\ref{eq:eqn2.2})
approaches the original conservation law (\ref{eq:eqn2.1}).
Substituting this expansion for ${\bf V}$ in (\ref{eq:eqn2.2}), the relaxation system can be seen as an
approximation to the original conservation law with a small dissipative correction
\begin{equation}
\label{eq:eqn2.4}
\frac{\partial {\bf C}}{\partial t}
+ \frac{\partial {\bf F(C)}}{\partial x}
= \varepsilon \frac{\partial }{\partial x}
\left( \left({\bf A}^{2}-{\bf F'(C)}^{2}  \right) \frac{\partial {\bf C}}{\partial x} \right)
+ {\rm O} {\rm (}\varepsilon ^{{\rm 2}} {\rm )},
\end{equation}
where ${\bf F'(C)}$ is the Jacobian of the flux function ${\bf F(C)}$.

In  (\ref{eq:eqn2.4}) $\varepsilon \frac{\partial }{\partial x}
\left( \left( {\bf A}^{2}-{\bf F'(C)}^{2}  \right) \frac{\partial {\bf C}}{\partial x} \right)$
is a ${\rm O} \left(\varepsilon \right)$ diffusive term, with $\varepsilon \left({\bf A}^{2} - {\bf F'(C)}^{2} \right)$
 the diffusion coefficient matrix. For (\ref{eq:eqn2.4}) to be well-posed $\left({\bf A}^{2} - {\bf F'(C)}^{2} \right)$
must be positive semi-definite for all ${\bf C}$ everywhere in the domain. This requirement on the diffusion coefficient matrix
$\left({\bf A}^{2} - {\bf F'(C)}^{2} \right)$, which is called the subcharacteristic condition, controls the
magnitude of the subcharacteristic speeds. In 1D, it is equivalent to
\begin{equation}
\label{eq:eqn2.5} \lambda ^{2}\le a^{2}, \;\;\;\;  \lambda
=\mathop{\max }\limits_{1\le p\le N} \left|\lambda _{p} \right|,  \;\;\;\;
a=\mathop{\min }\limits_{1\le p\le N} a_{p},
\end{equation}
where $\lambda _{p}$ are the eigenvalues of  the Jacobian ${\bf F'(C)}$.

For ${\bf C}$ in a bounded domain, the subcharacteristic condition can always be satisfied by choosing
subcharacteristic speeds that are sufficiently larger than the global maximum speed.
But the subcharacteristic speeds are also the characteristic speeds of the relaxation system and so
larger $\left\{a_{p} \right\}$ will necessitate  time steps smaller than that would have been needed
if the original system (\ref{eq:eqn2.1}) was solved directly.
As will be seen later in this section, large subcharacteristic speeds also increase numerical diffusion. Therefore,
$\left\{a_{p} \right\}$ is usually set to the smallest value that meets the stability criteria (\ref{eq:eqn2.4}),
typically the spectral radius of the Jacobian ${\bf F'(C)}$.

The stiff source term of the relaxation system can be effectively handled by operator splitting
\cite{jx, kg}, where the system is split into two sub-problems, a homogenous hyperbolic system
\begin{eqnarray}
\label{eq:eqn2.51}
\left[\begin{array}{cc}{{\bf C}}\\{{\bf V}}\end{array}\right]_{t}
+ {\bf B}
\left[\begin{array}{l} {{\bf C}}\\{{\bf V}}\end{array}\right]_{x}
= \left[\begin{array}{l} {{\rm {\bf 0}}} \\ {{\rm {\bf 0}}}\end{array}\right],
\;\;
{\bf B} = \left[\begin{array}{cc}{{\bf 0}}&{{\bf I}}\\{{\bf A}^{2}}&{{\bf 0}}\end{array}\right] ,
\end{eqnarray}
and a stiff ODE system
\begin{eqnarray}
\label{eq:eqn2.52}
\left[\begin{array}{l} {{\bf C}} \\ {{\bf V}}
\end{array}\right]_{t}
= \left[\begin{array}{l}
{\bf 0} \\
{\frac{1}{\varepsilon}
({\bf F(C)-V)}}
\end{array}\right] \nonumber ,
\end{eqnarray}
that are solved sequentially every time step. The ODE part can be solved exactly as
\begin{eqnarray}
{\bf C} &=& {\bf C}^{0}  ,
\nonumber\\[-1.5ex]
\label{eq:eqn2.6} \\[-1.5ex]
{\bf V} &=& {\bf F(C)}
\left(1-e^{-{\Delta t \mathord{\left/{\vphantom{\Delta t \varepsilon}}\right.\kern-\nulldelimiterspace}\varepsilon}} \right)
+ {\bf V}^{0} e^{-{\Delta t \mathord{\left/{\vphantom{\Delta t \varepsilon }}\right.\kern-\nulldelimiterspace}\varepsilon}}
.
\nonumber
\end{eqnarray}
A splitting technique introduced by Jin \cite{j95}, which retains the order of accuracy of the
underlying time-stepping \cite{kg}, is used throughout this work.
Since the matrix ${\bf B}$ is constant, the homogeneous hyperbolic system can be solved by first
obtaining the the characteristic variables through diagonalization of ${\bf B}$ and then upwinding on these variables
to obtain the numerical flux for the relaxation system. Of course, there are other approaches to handle the
stiff source term, like the IMEX approach proposed by \cite{paru}. However, as explained later in the paper,
the significant source of numerical diffusion is in the set up of homogenous system. Hence, in this work, we focus on improving
the homogenous subsystem, and for simplicity, we use Jin's splitting approach to handle the stiff source term.

For very small $\varepsilon $, ${\bf V}\approx {\bf F}\left({\bf C}\right)$.
So setting ${\bf V}={\bf F}\left({\bf C}\right)$ in the numerical flux, we get a first order, component-wise, semi-discrete
approximation, which Jin \& Xin refer to as the relaxed scheme. It is given by
\begin {subequations}
\label{eqn2.8}
\begin{equation}
\label{eq:eqn2.8a}
\frac{\partial {C}_{p,j} }{\partial t}
+ \frac{1}{\Delta x}
\left({\mathcal F}_{p,j+{\frac{1}{2}}}
- {\mathcal F}_{p,j-{\frac{1}{2}}}\right)= 0,
\end{equation}
\vspace {-0.2in}
\begin{equation}
\label{eq:eqn2.8b}
{\mathcal F}_{p,j+1/2}
= \frac{1}{2} \left({F}_{p,j} + {F}_{p,j+1}\right)
- \frac{1}{2} a_{p}\left({C}_{p,j+1} -{C}_{p,j} \right), \;\; p=1,2,\cdots N.
\end{equation}
\end{subequations}
The fully discrete version of this can be viewed as a generalized version of Lax-Friedrichs scheme;
for a specific choice of subcharacteristic speeds
$\left\{a_{p} = a = {\Delta x \mathord{\left/{\vphantom{\Delta x \Delta t}}\right.\kern-\nulldelimiterspace} \Delta t}
\right\}$, the Lax-Friedrichs scheme results.

To obtain a second order, semi-discrete, component-wise scheme, Jin \& Xin use van Leer's MUSCL \cite{vL}
reconstruction on the characteristic variables of (\ref{eq:eqn2.51}) and then set ${\bf V}={\bf F}\left({\bf C}\right)$. This leads to

\begin {subequations}
\label{eqn2.9}
\begin{equation}
\label{eq:eqn2.9a}
\frac{\partial {C}_{p,j} }{\partial t}
+ \frac{1}{\Delta x}
\left({\mathcal F}_{p,j+{\frac{1}{2}}}
- {\mathcal F}_{p,j-{\frac{1}{2}}}\right)
+ \frac{1}{\Delta x}
\left({\tilde {\mathcal F}}_{p,j+{\frac{1}{2}}}
- \tilde {\mathcal F}_{p,j-{\frac{1}{2}}}\right)= 0,
\end{equation}
where ${\mathcal F}_{p,j+1/2} $ is given in (\ref {eq:eqn2.8b}) and
${\tilde {\mathcal F}_{p,j+1/2}}$ is the second order correction
\begin{eqnarray}
\label{eq:eqn2.9c}
\tilde {\mathcal F}_{p,j+1/2}
&=& \frac{\Delta x}{4} \left( {\sigma} _{p,j}^{+} - {\sigma} _{p,j+1}^{-} \right),
\end{eqnarray}
with
\begin{eqnarray}
\sigma _{p,j}^{\pm } =  \frac{1}{\Delta x}
\left(
\left(F_{p} \pm a_{p} C_{p} \right)_{j+1}
- \left(F_{p} \pm a_{p} C_{p} \right)_{j}
\right)
\phi \left(\theta _{p,j}^{\pm } \right),
\label{eq:eqn2.9d}
\end{eqnarray}
\begin{equation}
\label{eq:eqn2.9e}
\theta _{p,j}^{\pm}
= \frac
{ \left(F_{p} \pm a_{p} C_{p} \right)_{j} - \left(F_{p} \pm a_{p} C_{p} \right)_{j-1} }
{\left(F_{p} \pm a_{p} C_{p} \right)_{j+1} -\left(F_{p} \pm a_{p} C_{p} \right)_{j} }.
\end{equation}
\end{subequations}
Throughout this paper we use the van Leer limiter for ${\phi (\theta )}$, but other limiters can be used as well.

While the JX relaxation methodology promises an efficient way of arriving at component-wise schemes, the JX scheme  introduces excessive numerical diffusion which gets
exagerrated in multidimensional problems, as shown in section 5.  A modified equation analysis \cite{kg}
on the scalar conservation law
shows that the  numerical diffusion coefficient
${\frac{\Delta x}{2} a \left(1-\frac{\Delta t}{\Delta x} \frac{F'}{a} F'\right)}$
of the first order JX scheme  is always greater than that of the corresponding upwind scheme,
${\frac{\Delta x}{2} F'\left(1-\frac{\Delta t}{\Delta x} F'\right)}$,
except in the trivial case of linear advection.
This is due to the restriction imposed on the minimum value of the subcharacteristic speed
 by (\ref{eq:eqn2.5}). The numerical diffusion increases not only with the magnitude of subcharacteristic speed,
but also with the increasing contrast between the subcharacteristic speed and local speeds of the original system.


\subsection{General constant relaxation systems}

In \cite{lp}, LeVeque and Pelanti presented a theory where relaxation systems are viewed as a means for defining approximate
Riemann solvers and present generalizations of the relaxation system, which they anticipate to lead to improved relaxation
schemes.
Instead of a relaxation system with symmetric characteristic speeds (i.e. negative eigenvalues = -positive eigenvalues),
a general relaxation system can be formulated as
\begin{eqnarray}
\label{eq:eqn2.10}
\left[  \begin{array}{c} {\bf C} \\ {\bf V} \end{array}  \right]_{t}
+ \left[  \begin{array}{cc} {\bf 0}&{\bf I}\\{\bf A}_{prod}&{\bf A}_{sum}  \end{array}  \right]
\left[  \begin{array}{c} {\bf C} \\ {\bf V} \end{array}  \right]_{x}
=  \left[ \begin{array}{c}
{\bf 0} \\ {\frac{1}{\varepsilon } \left( {\bf F}\left({\bf C}\right) - {\bf V} \right)}
 \end{array} \right],
\end{eqnarray}
where ${{\bf A}_{prod} = -{\bf A}_R{\bf A}_L, \;\;\;\;  {\bf A}_{sum} = {\bf A}_R + {\bf A}_L}$, with
${\bf A}_R=diag\left(a_{R,1} ,a_{R,2} \cdots a_{R,N} \right)$  and
${\bf A}_L=diag\left(a_{L,1} ,a_{L,2} \cdots a_{L,N} \right)$. As we will see below, the JX relaxation system is a
special case of the general relaxation (\ref{eq:eqn2.10}). This system has characteristic variables
$\left\{{\rm V}_{p} - a_{R,p} {\rm C}_{p} \right\}$ and $\left\{{\rm V}_{p} - a_{L,p} {\rm C}_{p} \right\}$, traveling
with speeds ${a_{R,p}}$ and ${a_{L,p}}$, respectively.  Here, the subscript R is used to denote waves moving to the right (positive speeds) and L to denote waves moving to the left (negative speeds). The subcharacteristic condition for this system is
\begin{eqnarray}
\left({\bf A}_{R} - {\bf F'(C)} \right)\left({\bf F'(C) - {\bf A}_{L}} \right) \ge 0,
\nonumber\\[-1.5ex]
\label{eq:eqn2.11} \\[-1.5ex]
\textnormal{i.e.}, \qquad
\mathop{\max }\limits_{1\le p\le N} \lambda _{p}  \le a_{R,p}, \qquad
\mathop{\min }\limits_{1\le p\le N} \lambda _{p} \ge a_{L,p},
\nonumber
\end{eqnarray}
where $\lambda _{p}$ are the eigenvalues of  the Jacobian ${\bf F'(C)}$.

Depending on the eigenvalues of the original conservation system and the choice of subcharacteristic speeds,
the general relaxation scheme can be either
\begin{enumerate}
\item
a one-sided system having either nonnegative or nonpositive  speeds
\begin{eqnarray}
a_{L,p}  =0,             \; a_{R,p} \ge \mathop{\max }\limits_{1\le p\le N} \left|\lambda _{p} \right|
\;\;\;  \textnormal{or} \;\;\;
 a_{L,p} \le - \mathop{\max }\limits_{1\le p\le N} \left|\lambda _{p} \right| ,      \; a_{R,p} =0
\nonumber ,
\end{eqnarray}
\item
a symmetric system, which is in fact the JX system with
\begin{eqnarray}
-a_{L,p} = a_{R,p} \ge \mathop{\max }\limits_{1\le p\le N} \left|\lambda _{p} \right|
\nonumber ,
\end{eqnarray}
\item
an optimal two-sided system, where the subcharacteristic speeds are chosen optimally based on the eigenvalues of the
original system,
\begin{eqnarray}
a_{L,p} \le \min \left( \min \limits_{1\le p\le N} \lambda _{p}  , 0 \right)
\;\;\; \textnormal{and} \;\;\;
a_{R,p} \ge \max \left(\max \limits_{1\le p\le N} \lambda _{p} , 0 \right)
\nonumber .
\end{eqnarray}
\end{enumerate}
The first choice, which can be used only when the original system has one-sided speeds,
has the least diffusion of the three. In the presence of mixed speeds,
this system fails the subcharacteristic condition  and the  solution blows up.
The optimal two-sided relaxation system reduces to the optimal one-sided system,
if the eigenvalues of the original system are either nonnegative or nonpositive. It reduces to the JX system when
the original system has both positive and negative eigenvalues, and the magnitude of the minimum negative eigenvalue is equal to
the magnitude of the maximum positive eigenvalue. While this choice is the most
promising of the three, a component-wise scheme which is developed from this optimal two-sided system, can still
exhibit large numerical diffusion when global maximum and minimum speeds are far from the average speeds.
The only way to develop a relaxation scheme that dynamically adapts numerical diffusion,
like the upwind scheme, is by using locally optimal subcharacteristic speeds. This motivated us to develop the variable relaxation scheme discussed below.

\section{Variable relaxation in 1D}
One way to construct a locally varying relaxation system is the nonconservative formulation \cite{kg}
\begin{eqnarray}
\label{eq:eqn3.1}
\hspace{-0.2in}
&{\bf q}_{t} \; + \; {\bf B}(x,t){\bf q}_{x} = {\bf s}, &
\\ \nonumber \\[-1.5ex]
\hspace{-0.7in}\textnormal{where}\; \;\;
{\bf q}=&\left[\begin{array}{c} {{\bf C}} \\ {{\bf V}} \end{array}\right] , \;
{\bf B}=\left[\begin{array}{cc} 0&I \\ A\left(x,t\right)^{2}& 0 \end{array}\right], \;
\textnormal{and}\; \;
{\bf s}=\left[ \begin{array}{c}
{\bf 0} \\ {\frac{1}{\varepsilon } \left( {\bf F}\left({\bf C}\right) - {\bf V} \right)}
 \end{array} \right]. &
\nonumber
\end{eqnarray}
Note that while ${{\bf C}}$ is a conserved variable, the relaxation variable ${{\bf V}}$ is not. As before,
${\bf A}\left(x,t\right)=diag\left(a_{1} \left(x,t\right),a_{2} \left(x,t\right)\cdots a_{\small{N}} \left(x,t\right)\right)$
is a positive diagonal matrix.
But the diagonal elements $\left\{a_{p} \left(x,t\right)\right\}$ are now the \textit {local} subcharacteristic speeds.
Again, using the Chapman-Enskog expansion to represent ${{\bf V}}$, the variable relaxation system can be seen as an
approximation to the original conservation law, plus a diffusive term
\begin{equation}
\label{eq:eqn3.2}
\frac{\partial {\bf C}}{\partial t}
+ \frac{\partial {\bf F}({\bf C})}{\partial x}
= \varepsilon \frac{\partial }{\partial x}
\left(\left( {\bf A}(x,t)^{2} - {\bf F'}({\bf C})^{2} \right )
\frac{\partial {\bf C}}{\partial x}  \right)
+ {\rm O} {\rm (}\varepsilon ^{{\rm 2}} {\rm )}.
\end{equation}
This leads to the requirement that $\left( {\bf A}(x,t)^{2} - {\bf F}'({\bf C})^{2} \right)$ must
be positive semi-definite for all ${\bf C}$, that is
\begin{equation}
\label{eq:eqn3.3}
\lambda(x,t) ^{2}\le a(x,t)^{2}, \;\;\;\;
\lambda(x,t) = \mathop{\max }\limits_{1\le p\le N} \left|\lambda _{p}(x,t) \right|,  \;\;\;\;
a(x,t)=\mathop{\min }\limits_{1\le p\le N} a_{p}(x,t),
\end{equation}
where  $\lambda _{p} (x,t)$ are the local speeds (for rarefactions or shocks) of the original conservation law. The subcharacteristic speeds must be chosen so that the positive
semi-definiteness of the diffusion coefficient matrix is guaranteed in either case \cite{liu}. In JX relaxation
this was done by setting the subcharacteristic speed to the global maximum speed,
\begin{eqnarray}
a=\mathop{\max }\limits_{1\le p\le N -1} \left|\lambda _{p} (x,t)\right|, \qquad \forall \left(x,t\right).
\nonumber
\end{eqnarray}
Here, the  local subcharacteristic speeds can be chosen in different ways as outlined in section 3.2.
\subsection{Variable relaxed schemes}
The approach to numerically solving the variable relaxation system is similar to Jin \& Xin's approach.
The only difference is in the way the homogenous part of the relaxation system is
solved. In the JX system, the relaxation matrix ${\bf B}$ is constant and hence diagonalizable into decoupled
characteristic variables. Here, the matrix ${\bf B}$ and hence the eigenvalues and eigenvectors are variable in ${x}$
and ${t}$. An attempt to diagonalize leads to a coupled system with  complicated source terms made up of
derivatives of the eigenvector matrix (see \cite{kg}). Another operator splitting on this system might not be helpful
since the
splitting error cannot be quantified and controlled easily. Instead, the homogenous part of the variable relaxation
system can be posed and solved easily as a set of Riemann problems at the cell interfaces
$\left\{x_{j+\frac{1}{2} } \right\}$  \cite[ch. 9]{lvq02}.

Given the vectors of unknowns of the relaxation system
\begin{eqnarray}
{\bf q}_{j}^{n} = \left[ {\bf C}_{j}^{n},\; {\bf V}_{j}^{n} \right]^{{\rm T}}
=\left[C_{1,j}^{n},\; C_{2,j}^{n},\; \ldots \; C_{N,j}^{n},\; V_{1,j}^{n},\; V_{2,j}^{n},\; \ldots \; V_{N,j}^{n}
\right] ^{{\rm T}} ,
\nonumber
\end{eqnarray}
the rate of change of the component average of the cell ${\mathcal C}_{{\it j}}$,
using the wave propagation form of the REA (Reconstruct-Evolve-Average) algorithm \cite{lvq02},
is given by the sum of the right-going  fluctuations at the left edge and  the left-going fluctuations at the right-edge,
that is
\begin{eqnarray}
\label{eq:eqn3.4}
\frac{\partial {\bf q}_{j}^{}}{\partial t}
= -\frac{1}{\Delta x} \left(
{\bf R}_{j-{\frac{1}{2}} } {\bf \Lambda}_{j-{\frac{1}{2}} }^{+} {\bf L}_{j-{\frac{1}{2}}}
\Delta {\bf q}_{j-{\frac{1}{2}}}
+ {\bf R}_{j+{\frac{1}{2}}} {\bf  \Lambda} _{j+{\frac{1}{2}}}^{-}{\bf L}_{j+{\frac{1}{2}}}
\Delta {\bf q}_{j+{\frac{1}{2}}} \right),
\end{eqnarray}
where $\Delta {\bf q}_{j-{\frac{1}{2}} } ={\bf q}_{j} -{\bf q}_{j-1} $. Here,
${{\Lambda} _{j-{\frac{1}{2}}}={\Lambda} _{j-{\frac{1}{2}}}^{+}+{\Lambda} _{j-{\frac{1}{2}}}^{-}}$
consists of the  eigenvalues of the relaxation matrix ${\bf B}_{j-{\frac{1}{2}}}$ given by
\begin{eqnarray}
{\bf \Lambda}_{j-{\frac{1}{2}}}^{+}
&=&diag\left(a_{1,j-{\frac{1}{2}}}^{+},\; 0, \; a_{2,j-{\frac{1}{2}}}^{+}, \; 0\;
\cdots\;  a_{N,j-{\frac{1}{2}}}^{+},\; 0\right),
\nonumber \\ [-1.5ex]
\label{eq:eqn3.5a} \\[-1.5ex]
{\bf \Lambda}_{j-{\frac{1}{2}}}^{-}
&=&diag\left(0,a_{1,j-{\frac{1}{2}}}^{-}, \; 0, \; a_{2,j-{\frac{1}{2}}}^{-}, \; 0 \;
\cdots \;  0, \; a_{N,j-{\frac{1}{2}}}^{-} \right),
\nonumber
\end{eqnarray}
where ${a_{p,j-{\frac{1}{2}} }^{\pm}}$ are the local subcharacteristic speeds. These speeds are assumed to be piecewise constant and
are reset in every time step as outlined in section 3.2.

\noindent ${\bf R}_{j-{\frac{1}{2}} } = \left[\bfr_{1,j-{\frac{1}{2}}}^{+}, \; \bfr_{1,j-{\frac{1}{2}} }^{-}, \;
\ldots\; \bfr_{N,j-{\frac{1}{2}}}^{+}, \;  \bfr_{N,j-{\frac{1}{2}}}^{-} \right]$ is the eigenvector matrix of
${\bf B}_{j-{\frac{1}{2}}}$ with
${{\bfr_{p,j-{\frac{1}{2}} }^{+}}}$ the columns corresponding to the left  going waves
\begin{subequations}
\label {eqn3.6}
\vspace {-0.1in}
\begin{eqnarray}
\begin{array}{l} {\bfr_{p,j-{\frac{1}{2}} }^{+}
=\left[{\rm 0\; }\ldots {\rm \; \; 1\; \; 0\; \; }\cdots
{\rm \; 0\; \; }a_{p,j-{\frac{1}{2}} }^{+} {\rm \; 0\; \; }\cdots {\rm \; 0}\right]^{{\rm T}} },
\\ {\hspace{1.03in}\uparrow
\hspace{0.7in}\uparrow }
\\  \hspace{0.8in} {  p^{{\rm th}}
{\rm element}\hspace{0.2in}\left(N +p\right)^{{\rm th}} {\rm element}}
\end{array}
\label{eq:eqn3.6a}
\end{eqnarray}
and  ${\bfr_{p,j-{\frac{1}{2}}}^{-}}$  the columns corresponding to  the right going waves
\vspace {-0.1in}
\begin{eqnarray}
\label{eq:eqn3.6b}
\begin{array}{l} {\bfr_{p,j-{\frac{1}{2}}}^{-}
=\left[{\rm 0\; }\ldots {\rm \; \; 1\; \; 0\; \; }
\cdots {\rm \; 0\; \; }a_{p,j-{\frac{1}{2}}}^{-} {\rm \; 0\; \; }\cdots {\rm \; 0}\right]^{{\rm T}} }.
\\ {\hspace{1.03in}\uparrow
\hspace{0.7in}\uparrow }
\\\hspace{0.8in} {p^{{\rm th}}
{\rm element} \hspace{0.2in} \left(N+p\right)^{{\rm th}} {\rm element}}
\end{array}
\end{eqnarray}
\end{subequations}
${\bf L}_{j-{\frac{1}{2}} }^{n}=\left[\bfell _{1,j-{\frac{1}{2}}}^{+}, \; \bfell _{1,j-{\frac{1}{2}} }^{-}, \;
\ldots \; \bfell _{N,j-{\frac{1}{2}}}^{+}, \; \bfell _{N,j-{\frac{1}{2}}}^{-} \right]^{T} $
is the inverse eigenmatrix
with ${\bfell _{p,j-{\frac{1}{2}}}^{+}}$ the rows corresponding to the left  going waves
\begin{subequations}
\label {eqn3.7}
\vspace {-0.1in}
\begin{eqnarray}
\label{eq:eqn3.7a}
\begin{array}{l}
{\bfell _{p,j-{\frac{1}{2}}}^{+}
=\left[{\rm 0\; }\ldots {\rm \; \; }\frac{-a_{p,j-{\frac{1}{2}}}^{-} }{a_{p,j-{\frac{1}{2}}}^{+}-a_{p,j-{\frac{1}{2}} }^{-}}
{\rm \; \; 0\; }\cdots {\rm \; 0\; \; }\frac{1}{a_{p,j-{\frac{1}{2}}}^{+} -a_{p,j-{\frac{1}{2}}}^{-}}
{\rm \; 0\; }\cdots {\rm \; 0}\right]{\rm \; }},
\\ {\hspace{1.4in}\uparrow
\hspace{1.4in}\uparrow }
\\ {\hspace{1.2in}p^{{\rm th}}
{\rm element}    \hspace{0.5in}
\left(N+p\right)^{{\rm th}} {\rm element}}
\end{array}
\end{eqnarray}
and ${\bfell _{p,j-{\frac{1}{2}}}^{-}}$  the rows corresponding to the right going waves
\vspace {-0.1in}
\begin{eqnarray}
\label{eq:eqn3.7b}
\begin{array}{l} \hspace{1in}{\bfell _{p,j-{\frac{1}{2}}}^{-}
=\left[{\rm 0\; }\ldots {\rm \; \; }\frac{a_{p,j-{\frac{1}{2}} }^{+}}{a_{p,j-{\frac{1}{2}}}^{+}-a_{p,j-{\frac{1}{2}}}^{-}}
{\rm \; \; 0\; }\cdots {\rm \; 0\; }\frac{-1}{a_{p,j-{\frac{1}{2}}}^{+} -a_{p,j-{\frac{1}{2}}}^{-}}
{\rm \; 0\; }\cdots {\rm \; 0}\right]}.
\\ {\hspace{2.3in}\uparrow
\hspace{1.4in}\uparrow }
\\ \hspace{2in} {p^{{\rm th}}
{\rm element} \hspace{0.5in}
\left(N+p\right)^{{\rm th}} {\rm element}}
 \end{array}
\end{eqnarray}
\end{subequations}
Using (\ref{eqn3.6}) and (\ref{eqn3.7}) in (\ref{eq:eqn3.4}), a first order, component-wise, relaxed,
semi-discrete update for the ${p^{th}}$ component follows as
\begin {subequations}
\label {eqn3.8}
\begin{equation}
\label{eq:eqn3.8a}
\frac{\partial C_{p,j} }{\partial t}
+ \frac{1}{\Delta x} \left( \mathcal {F}_{p,j+{\frac{1}{2}}} - \mathcal{F}_{p,j-{\frac{1}{2}}} \right) = 0,
\end{equation}
with
\begin{eqnarray}
\label{eq:eqn3.8b} \;\;\; \qquad
\mathcal{F}_{p,j-{\frac{1}{2}}}
=  F_{p,j-1}
+ \frac{-a_{p,j-{\frac{1}{2}}}^{-}} {a_{p,j-{\frac{1}{2}}}^{+} -a_{p,j-{\frac{1}{2}}}^{-} }
\left[ F_{p,j} -F_{p,j-1} -a_{p,j-{\frac{1}{2}}}^{+} \left(C_{p,j} -C_{p,j-1} \right)\right].
\end{eqnarray}
\end {subequations}

The second order update is obtained by adding high resolution correction vectors to the matrix form of the
first order update (\ref{eq:eqn3.4}).
In the JX scheme, the high resolution terms  were obtained directly by limiting on the characteristic variables.
Here, since the characteristic variables are not available, the high resolution corrections are obtained
by limiting on the change in  ${2N}$ characteristic waves (see \cite[p.182]{lvq02}).
Across the interface $\left\{x_{j+\frac{1}{2} } \right\}$ the change in the
left going and right going characteristic waves are constructed as, respectively,
\begin{eqnarray}
\bfw_{p,j-{\frac{1}{2}}}^{+}
=\left( \bfell _{p,j-{\frac{1}{2}}}^{+} \bullet \Delta {\bf q}_{j-{\frac{1}{2}}} \right)
{\bfr}_{p,j-{\frac{1}{2}}}^{+}
\quad \textnormal{and} \quad
\bfw_{p,j-{\frac{1}{2}}}^{-}
= \left( \bfell _{p,j-{\frac{1}{2}}}^{-} \bullet \Delta {\bf q}_{j-{\frac{1}{2}}} \right)
 {\bfr}_{p,j-{\frac{1}{2}}}^{-} .
\nonumber
\end{eqnarray}
The limiting parameters for each component are then given by
\begin{eqnarray}
\theta _{p,j-{\frac{1}{2}} }^{+}
=\frac
{\bfw_{p,j-1-{\frac{1}{2}}}^{+} \bullet \bfw_{p,j-{\frac{1}{2}}}^{+}} 
{\bfw_{p,j-{\frac{1}{2}}}^{+} \bullet \bfw_{p,j-{\frac{1}{2}}}^{+}} 
, \qquad
\theta _{p,j-{\frac{1}{2}} }^{-}
=\frac
{\bfw_{p,j+1-{\frac{1}{2}}}^{-} \bullet \bfw_{p,j-{\frac{1}{2}}}^{- }} 
{\bfw_{p,j-{\frac{1}{2}}}^{-} \bullet \bfw_{p,j-{\frac{1}{2}}}^{-}}   
\nonumber ,
\end{eqnarray}
leading to the high resolution correction vector
\begin{equation}
\label{eq:eqn3.9} \quad
\bftf_{j-{\frac{1}{2}}}
=\frac{1}{2} \sum _{p=1,2...N}
a_{p,j-{\frac{1}{2}}}^{+} \phi \left(\theta _{p,j-{\frac{1}{2}}}^{+} \right) \bfw_{p,j-{\frac{1}{2}}}^{+}
- a_{p,j-{\frac{1}{2}}}^{-} \phi \left(\theta _{p,j-{\frac{1}{2}}}^{-} \right)\bfw_{p,j-{\frac{1}{2}}}^{-}.
\end{equation}
The above corrections are then added to the matrix form of the first order update (\ref{eq:eqn3.4})
to obtain a second order, semi-discrete  update
\begin{eqnarray}
\frac{\partial {\bf q}_{j} }{\partial t}
=& -&\frac{1}{\Delta x} \left(
{\bf R}_{j-{\frac{1}{2}} } {\bf \Lambda}_{j-{\frac{1}{2}} }^{+} {\bf L}_{j-{\frac{1}{2}}}
\Delta {\bf q}_{j-{\frac{1}{2}}}
+ {\bf R}_{j+{\frac{1}{2}}} {\bf  \Lambda} _{j+{\frac{1}{2}}}^{-}{\bf L}_{j+{\frac{1}{2}}}
\Delta {\bf q}_{j+{\frac{1}{2}}}
\right)
\nonumber\\[-1.5ex]
\label{eq:eqn3.10} \\[-1.5ex]
&-& \frac{1}{\Delta x} \left(\bftf_{j+{\frac{1}{2}}}
-\bftf_{j-{\frac{1}{2}}} \right).
\nonumber
\end{eqnarray}
The individual component-wise  updates are
\begin {subequations}
\label {eqn3.11}
\begin{equation}
\label{eq:eqn3.11a}
\frac{\partial C_{p,j} }{\partial t}
= - \frac{1}{\Delta x}
\left( \mathcal{F}_{p,j+{\frac{1}{2}}} - \mathcal{F}_{p,j-{\frac{1}{2}}} \right)
- \frac{1}{\Delta x}
\left( \tilde{\mathcal {F}}_{p,j+{\frac{1}{2}}} -\tilde{\mathcal{F}}_{p,j-{\frac{1}{2}}} \right),
\end{equation}
where ${\mathcal{F}_{p,j-{\frac{1}{2}}}}$ is given by (\ref{eq:eqn3.8b}) and, the component-wise correction terms are
\begin{eqnarray}
\tilde{\mathcal{F}}_{p,j-{\frac{1}{2}}}
&=& \frac
{ a_{p,j-{\frac{1}{2}}}^{+} }
{a_{p,j-{\frac{1}{2}}}^{+} -a_{p,j-{\frac{1}{2}}}^{-}}
\frac {\phi \left(\theta _{p,j-{\frac{1}{2}}}^{+} \right)} {2}
\left[ - a_{p,j-{\frac{1}{2}}}^{-}\left(C_{p,j} -C_{p,j-1} \right)
+ \left(F_{p,j} -F_{p,j-1} \right) \right]
\nonumber\\
\qquad & - &\frac
{a_{p,j-{\frac{1}{2}}}^{-}}
{a_{p,j-{\frac{1}{2}}}^{+} - a_{p,j-{\frac{1}{2}}}^{-}}
\frac{ \phi \left(\theta _{p,j-{\frac{1}{2}}}^{-} \right)} {2}
\left[ a_{p,j-{\frac{1}{2}}}^{+}\left(C_{p,j} -C_{p,j-1}  \right)
-\left(F_{p,j} -F_{p,j-1} \right) \right],
\label{eq:eqn3.11b}
\end{eqnarray}
with the component-wise limiter-parameters
\begin{eqnarray}
\qquad \qquad \theta _{p,j-{\frac{1}{2}}}^{+}
&=& \frac{
\frac  {1 + a_{p,j-{\frac{1}{2}}}^{+} a_{p,j-{\frac{3}{2}}}^{+} }
{a_{p,j-{\frac{3}{2}}}^{+} - a_{p,j-{\frac{3}{2}}}^{-} }
\left[ -a_{p,j-{\frac{3}{2}}}^{-}\left(C_{p,j-1} -C_{p,j-2} \right)
+ \left( F_{p,j-1} -F_{p,j-2} \right)  \right]}
{
\frac{1 + a_{p,j-{\frac{1}{2}}}^{+} a_{p,j-{\frac{1}{2}}}^{+} }
{a_{p,j-{\frac{1}{2}}}^{+} -a_{p,j-{\frac{1}{2}}}^{-} }
\left[ -a_{p,j-{\frac{1}{2}} }^{-} \left(C_{p,j} - C_{p,j-1} \right)
+ \left( F_{p,j} -F_{p,j-1} \right) \right]
} ,
\nonumber\\[-1.5ex]
\label{eq:eqn3.11c} \\[-1.5ex]
\theta _{p,j-{\frac{1}{2}}}^{-}
&=& \frac{
\frac  {1+a_{p,j+{\frac{1}{2}}}^{-} a_{p,j-{\frac{1}{2}}}^{-}}
{a_{p,j+{\frac{1}{2}}}^{+} -a_{p,j+{\frac{1}{2}}}^{-}}
\left[
a_{p,j+{\frac{1}{2}}}^{+} \left(C_{p,j+1} -C_{p,j} \right)
-  \left( F_{p,j+1} -F_{p,j}   \right)
\right]
}
{
\frac{1+a_{p,j-{\frac{1}{2}}}^{-} a_{p,j-{\frac{1}{2}}}^{-}}
{a_{p,j-{\frac{1}{2}}}^{+} - a_{p,j-{\frac{1}{2}}}^{-}}
\left[ a_{p,j-{\frac{1}{2}}}^{+}\left(C_{p,j} -C_{p,j-1} \right)
 -  \left( F_{p,j} - F_{p,j-1}  \right)   \right]
} .
\nonumber
\end{eqnarray}
\end {subequations}
\subsection {Choice of subcharacteristic speeds}

For a choice of constant subcharacteristics
\begin{eqnarray}
- a_{p,j-{\frac{1}{2}}}^{-} = a_{p,j-{\frac{1}{2}}}^{+} = a_{p},
\nonumber
\end{eqnarray}
the equations (\ref{eqn3.8}) and (\ref{eqn3.11}) of the variable relaxed scheme
 reduce to those of  the JX scheme (\ref{eqn2.8}) and (\ref{eqn2.9}).
There is a subtle difference between the limiting parameters of the left-going wave in (\ref{eq:eqn2.9e}) and
(\ref{eq:eqn3.11c}). The limiter-parameter  ${\theta _{p,j-{\frac{1}{2}}}^{-}}$ as given in
(\ref{eq:eqn3.11c}), with a choice of constant subcharacteristics,
 is the reciprocal of the JX limiter parameter  ${\theta _{p,j}^{-}}$ (\ref{eq:eqn2.9e}).
However the  second order correction terms of (\ref{eq:eqn3.11b}) will still be equal to the
JX second order correction terms (\ref{eq:eqn2.9c}-\ref{eq:eqn2.9d})
because the van Leer limiter (and also other limiters like minmod, superbee, MC) obeys the symmetry condition
$ \frac{\phi \left(\theta \right)}{ \theta } = \phi \left( \frac{1}{\theta}  \right).  $

Rather than constant subcharacteristics, we allow them to vary locally. Such locally varying subcharacteristics can be chosen in two ways: \\
1. {\bf Symmetric speeds}: We can set
\begin{eqnarray}
\label{eq:eqn3.12}
a_{p,j-{\frac{1}{2}}}^{+} = - a_{p,j-{\frac{1}{2}}}^{-} = a_{p,j-{\frac{1}{2}}}
 =  \mathop{\max }\limits_{1\le p\le N} {\left|\lambda _{p}\left( {\bf C} \right) \right|},
\end{eqnarray}
for all ${\bf C}$ between ${{\bf C}_{j}}$ and ${{\bf C}_{j-1}}$ with $\lambda _{p}$ the eigenvalues of the
Jacobian ${\bf F'(C)}$. Then, the first order numerical flux from (\ref{eq:eqn3.8b}) reduces to  the LLF flux
\begin{equation}
\label{eq:eqn3.13}
\qquad \mathcal{F}_{p,j-{\frac{1}{2}}}^{n}
=  F_{p,j-1}^{n}
+ \frac{1} {2}
\left[\left( F_{p,j}^{n} -F_{p,j-1}^{n} \right)-a_{p,j-{\frac{1}{2}}} \left(C_{p,j}^{n} -C_{p,j-1}^{n} \right)\right].
\end{equation}
The second order terms of (\ref{eq:eqn3.11b})  and (\ref{eq:eqn3.11c}) will simplify as
\begin {subequations}
\label {eqn3.14}
\begin{eqnarray}
\tilde{\mathcal{F}}_{p,j-{\frac{1}{2}}}^{n}
&=& \frac {\phi \left(\theta _{p,j-{\frac{1}{2}}}^{+} \right)} {4}
\left[ a_{p,j-{\frac{1}{2}}} \left(C_{p,j}^{n} -C_{p,j-1}^{n} \right)
+ \left(F_{p,j}^{n} -F_{p,j-1}^{n} \right) \right]
\nonumber\\[-1.5ex]
\label{eq:eqn3.14a} \\[-1.5ex]
& + &
\frac{ \phi \left(\theta _{p,j-{\frac{1}{2}}}^{-} \right)} {4}
\left[ a_{p,j-{\frac{1}{2}}} \left(C_{p,j}^{n} -C_{p,j-1}^{n}  \right)
-\left(F_{p,j}^{n} -F_{p,j-1}^{n} \right) \right],
\nonumber
\end{eqnarray}
\begin{eqnarray}
\qquad \theta _{p,j-{\frac{1}{2}}}^{+}
&=& \frac{
\frac  {1 + a_{p,j-{\frac{1}{2}}} a_{p,j-{\frac{3}{2}}} }
{2a_{p,j-{\frac{3}{2}}}}
\left[a_{p,j-{\frac{3}{2}}}\left(C_{p,j-1}^{n} -C_{p,j-2}^{n} \right)
+ \left( F_{p,j-1}^{n} -F_{p,j-2}^{n} \right)  \right]}
{
\frac{1 + a_{p,j-{\frac{1}{2}}} a_{p,j-{\frac{1}{2}}}}
{2a_{p,j-{\frac{1}{2}}}}
\left[ a_{p,j-{\frac{1}{2}}} \left(C_{p,j}^{n} - C_{p,j-1}^{n} \right)
+ \left( F_{k,j}^{n} -F_{k,j-1}^{n} \right) \right]
},
\nonumber\\[-1.5ex]
\label{eq:eqn3.14b} \\[-1.5ex]
\theta _{p,j-{\frac{1}{2}}}^{-}
&=& \frac{
\frac  {1+a_{p,j+{\frac{1}{2}}} a_{p,j-{\frac{1}{2}}}}
{2a_{p,j+{\frac{1}{2}}}}
\left[
a_{p,j+{\frac{1}{2}}} \left(C_{p,j+1}^{n} -C_{p,j}^{n} \right)
-  \left( F_{p,j+1}^{n} -F_{p,j}^{n}   \right)
\right]
}
{
\frac{1+a_{p,j-{\frac{1}{2}}} a_{p,j-{\frac{1}{2}}}}
{2a_{p,j-{\frac{1}{2}}}}
\left[ a_{p,j-{\frac{1}{2}}}\left(C_{p,j}^{n} -C_{p,j-1}^{n} \right)
 -  \left( F_{p,j}^{n} - F_{p,j-1}^{n}  \right)   \right]
}.
\nonumber
\end{eqnarray}
\end {subequations}
2. {\bf Optimal speeds}: We can also prescribe asymmetric speeds that lead to a  optimal scheme that adapts itself to
become upwind in the presence of one-sided fluxes, as
\begin{eqnarray}
\label{eq:eqn3.15}
&{a_{p,j-{\frac{1}{2}}}^{+} = a_{j-{\frac{1}{2}}}^{+}
 = \max \left(\max \limits_{1\le p\le N} \lambda _{p}\left( {\bf C} \right) , 0 \right)}&
\;\;\; \textnormal{and} \;\;\;  \\
&{a_{p,j-{\frac{1}{2}}}^{-} = a_{j-{\frac{1}{2}}}^{-}
= \min \left( \min \limits_{1\le p\le N} \lambda _{p}\left( {\bf C} \right) , 0 \right)}&
\nonumber
\end{eqnarray}
over all {\bf C} between ${{\bf C_{j}^{n}}}$   and ${{\bf C_{j-1}^{n}}}$.  This is the same as the HLL solver developed by Harten, Lax and van Leer \cite{hll}. The
first order KNP scheme \cite{knp} also coincides with the first order variable relaxed scheme with the subcharacteristic speeds
chosen as above.
The same optimal scheme can also be derived by formulating an asymmetric general variable relaxation system
\begin{eqnarray}
\left[ \begin{array}{c} {\bf C} \\ {\bf V} \end{array}  \right]_{t}
+ \left[\begin{array}{cc} {\bf 0}&{\bf I}\\{\bf A}_{prod}\left(x,t\right) &{\bf A}_{sum}\left(x,t\right) \end{array} \right]
\left[  \begin{array}{c} {\bf C} \\ {\bf V} \end{array}  \right]_{x}
=  \left[ \begin{array}{c}
{\bf 0} \\ {\frac{1}{\varepsilon } \left( {\bf F}\left({\bf C}\right) - {\bf V} \right)}
 \end{array} \right],
\nonumber
\end{eqnarray}
with the subcharacteristic condition
\begin{eqnarray}
\left({\bf A}^{+}\left(x,t\right) - {\bf F'(C)} \right)
\left({\bf F'(C) - {\bf A}^{-}} \left(x,t)\right) \right) \ge 0,
\nonumber
\end{eqnarray}
where $\;\;{{\bf A}_{prod}(x,t) = -{\bf A}^{+}(x,t){\bf A}^{-}(x,t), \;\;\;\;
{\bf A}_{sum} = {\bf A}^{+}(x,t) + {\bf A}^{-}(x,t)}$, with
\begin{eqnarray}
&{\bf A}^{+}(x,t)&=diag\left(a_{1}^{+}(x,t)  ,a_{2}^{+}(x,t) \cdots a_{N}^{+}(x,t) \right),
\nonumber \\
&{\bf A}^{-}(x,t)&=diag\left(a_{1}^{-}(x,t) , a_{2}^{-}(x,t) \cdots a_{N}^{-}(x,t) \right),
\nonumber \\
\nonumber
\end{eqnarray}
and ${\left\{ a_{p}^{\pm}(x,t)\right\}}$ chosen as in (\ref{eq:eqn3.15}).

\subsection{Stability properties}

We can show that the first order, fully discrete, variable relaxed schemes proposed above are monotonic under appropriate conditions on the subcharacteristics:

\begin{theorem}
\label{th:th3.1}
The first order, fully discrete,  variable relaxed scheme
\begin{eqnarray}
{C}_{j}^{n+1} = {C}_{j}^{n}
 - \frac{\Delta t}{\Delta x } \left({\mathcal F}_{j+{\frac{1}{2}}}^{n} - {\mathcal F}_{j-{\frac{1}{2}}}^{n} \right),
\label {eq:eqn3.16}
\end{eqnarray}
with symmetric speeds, where ${\mathcal F}_{j\pm{\frac{1}{2}}}^{n}$ is given by (\ref {eq:eqn3.13}),
is monotonic under the local subcharacteristic condition
\begin{equation}
a_{j-{\frac{1}{2}} } \ge \left|F'\left(C\right)_{j-1}^{n}\right| \;\; \textnormal{and} \;\;
a_{j-{\frac{1}{2}}} \ge \left|F'\left(C\right)_{j}^{n}\right|,
\label{eq:eqn3.17}
\end{equation}
and the time step restriction $\frac{\Delta t}{\Delta x} a_{\max } \le 1$, where $a_{\max }$ is the
maximum subcharacteristic speed.
The first order scheme (\ref {eq:eqn3.16})  with optimal speeds,
where ${\mathcal F}_{j\pm{\frac{1}{2}}}^{n}$ is given by (\ref {eq:eqn3.8b}),  is monotonic under
the local subcharacteristic condition
\begin{equation}
a_{j-{\frac{1}{2}} }^{-} \le    \min \left( F'(C)_{j-1}^{n}, 0\right), \;\; \textnormal{and} \;\;
a_{j-{\frac{1}{2}} }^{+} \ge    \max \left( F'(C)_{j}^{n}, 0\right),
\label{eq:eqn3.18}
\end{equation}
and the time step restriction $\frac{\Delta t}{\Delta x} a_{\max } \le \frac {1}{2}$.
\end {theorem}
\newline
Proof given in the appendix.
\vspace{0.1in}
\begin{theorem}
\label{th:th3.2}
The variable relaxed scheme with second order spatial discretization and forward Euler time stepping
\begin{equation}
\label{eq:eqn3.19}
C_{j}^{n+1} =C_{j}^{n} -  \frac{1}{\Delta x}
\left( \mathcal{F}_{j+{\frac{1}{2}}}^{n} - \mathcal{F}_{j-{\frac{1}{2}}}^{n} \right)
- \frac{1}{\Delta x}
\left( \tilde{\mathcal {F}}_{j+{\frac{1}{2}}}^{n} -\tilde{\mathcal{F}}_{j-{\frac{1}{2}}}^{n} \right),
\end{equation}
with $\mathcal{F}_{j \pm {\frac{1}{2}}}^{n}$  given by (\ref{eq:eqn3.8b}) and
$\tilde{\mathcal {F}}_{j \pm {\frac{1}{2}}}^{n}$  given by (\ref{eq:eqn3.11b}),
is TVD  under the CFL condition $\frac{\Delta t}{\Delta x} a_{\max } \le \frac{1}{2} $, and the local subcharacteristic
condition
\begin{equation}
a_{j-{\frac{1}{2}} }
\ge \left| \frac{ F\left(C\right)_{j}^{n} -  F\left(C\right)_{j-1}^{n} }{C_{j}^{n} -C_{j-1}^{n} } \right|,
\label{eq:eqn3.20}
\end{equation}
for the symmetric case and,
\begin{eqnarray}
\label{eq:eqn3.21}
\hspace {0.4in}
a_{j-{\frac{1}{2}} }^{-} \le
\min \left( \frac{ F\left(C\right)_{j}^{n} -  F\left(C\right)_{j-1}^{n} }{C_{j}^{n} -C_{j-1}^{n} }, 0 \right),
a_{j-{\frac{1}{2}} }^{+} \ge \max
 \left( \frac{ F\left(C\right)_{j}^{n} -  F\left(C\right)_{j-1}^{n} }{C_{j}^{n} -C_{j-1}^{n} }, 0 \right),
\end{eqnarray}
for the optimal case.
\end{theorem}
\newline
Proof given in the appendix.
\newline
Using the lemma due to Shu and Osher \cite{so} and the above theorem,
a second order variable relaxed scheme with RK-2 time stepping can also be shown to be TVD in a straightforward manner.

\section{Multidimensional Relaxation}
In multidimensions the relaxed schemes can be obtained either by considering relaxation dimension-by-dimension or, equivalently, by
simply applying the 1D flux (\ref {eq:eqn2.8b}, \ref {eq:eqn2.9c}, \ref {eq:eqn3.8b}, \ref {eq:eqn3.11b}) dimension-wise.
The restriction due to the subcharacteristic condition however,
becomes more severe with increase in the number of dimensions. In this section we present a framework for choosing  local
subcharacteristic speeds and deriving variable  schemes for general multidimensional conservation laws and specifically
show the semi-discrete updates for the 2D case. The multidimensional JX relaxation is revisited and analyzed here for completeness.
\subsection {Jin-Xin relaxation in multidimensions}
Consider the \textit{m}-dimensional conservation system
\begin{equation}
\label{eq:eqn4.1}
\frac{\partial {\bf C}}{\partial t}
+ \sum _{k=1}^{m}\frac{\partial {\bf F}_{k} }{\partial x_{k} }  = 0.
\hspace{.2in} {\bf C}, {\bf F}  \in R^{N},
\end{equation}
Jin \& Xin formulate the corresponding relaxation system as
\begin{eqnarray}
\frac{\partial {\bf C}}{\partial t}
&+& \sum _{k=1}^{m}\frac{\partial {\bf V}_{k} }{\partial x_{k} }  = 0,
\hspace{.3in}  {\bf C}{\bf ,V}_{k} \in R^{N},
\nonumber\\[-1.5ex]
\label{eq:eqn4.2}  \\[-1.5ex]
\frac{\partial {\bf V}_{k} }{\partial t}
&+&  {\bf A}_{k}^{2} \frac{\partial {\bf C}}{\partial x_{k} }
= \frac{1}{\varepsilon } \left( {\bf F}_{k} \left ({\bf C} \right) - {\bf V}_{k} \right),
\hspace{.3in} k=1,2,\cdots m ,
\nonumber
\end{eqnarray}
where ${\bf A}_{k}=diag\left(a_{1}^k ,a_{2}^k \cdots a_{N}^k \right)$
is a positive diagonal matrix.
As before, the relaxation system can be expressed as an approximation
to the original system of conservation laws with a small dissipative correction,  that is
\begin{equation}
\label{eq:eqn4.3}
\frac{\partial {\bf C}}{\partial t}
+ \sum _{k=1}^{m}\frac{\partial {\bf F}_{k} }{\partial x_{k} }
=   \varepsilon \sum _{k,\tilde{k}=1}^{m} \frac{\partial }{\partial x_{k}}
\left[
\left( \partial _{k\tilde{k}} {\bf A}_{k}^{2} -{\bf F}_{k} {\bf F}_{\tilde{k}} \right)
\frac{\partial {\bf C}}{\partial x_{\tilde{k}} }
\right]
+ {\rm O} \left( {\varepsilon ^2} \right),
\end{equation}
where $\partial _{k\tilde{k}} $ is the Kronecker delta.  Since $\varepsilon $ is positive, the system
(\ref{eq:eqn4.2}, \ref{eq:eqn4.3}) will be stable if
$\left(\partial _{k\tilde{k}} {\bf A}_{k}^{2} -{\bf F}_{k} {\bf F}_{\tilde{k}} \right) \ge 0$.
This subcharacteristic condition can be expressed in terms of the eigenvalues $\left\{ \lambda _{p}^k \right\}$
of the Jacobian ${\bf F}_{k} \left ({\bf C} \right)$ and the subcharacteristic speeds $\left\{ a _{p}^k \right\}$
for each dimension $k$ as
\begin{equation}
\label{eq:eqn4.4}
\frac{\lambda _{^{1,\max } }^{2} }{a_{^{1,min} }^{2} }
+\frac{\lambda _{^{2,\max } }^{2} }{a_{^{2,min} }^{2} }
+\cdots
+\frac{\lambda _{^{m,\max } }^{2} }{a_{^{m,min}}^{2} } \le 1 ,
\end{equation}
where $\lambda _{k, \max } = \mathop{\max }\limits_{1\le p\le N} \left|\lambda _{p}^k \right|$ and
$a_{k, \min}  = \mathop{\min }\limits_{1\le p\le N} \left|a _{p}^k \right|$. The subcharacteristic variables are
chosen to minimize numerical diffusion, while obeying the condition (\ref{eq:eqn4.4}).
Note that this is a more severe restriction than the 1D subcharacteristic condition, and hence can result in very
diffusive solutions.

To derive the component-wise schemes, operator splitting can be used like in the 1D approach. In 1D, the homogenous
hyperbolic part  in the JX relaxation is diagonalized to obtain the characteristic variables. In the
multidimensional case  diagonalization is possible only if the relaxation matrices commute. This is not the case for Jin \& Xin's relaxation matrices. Therefore, the JX scheme is derived by diagonalizing the relaxation system dimension-by-dimension and upwinding on the
resulting dimension-wise characteristic variables $\left\{{\rm V}_{k,p} \pm \; a_{p}^k {\rm C}_{p} \right\}$.

We remark that a commutative formulation is  possible
only for a relaxation formulation that leads to one-sided eigenvalues, such as in the first case discussed in Section 2.2. Such as system can only be formulated
for a conservation system with one-sided speeds.  Even then, the commutative formulation
necessitates the same subcharacteristic speed in every dimension. In problems where there is preferential
flow in some directions, this generally results in much more diffusive solutions than those obtained with JX relaxation.

We consider the 2D system
\begin{equation}
\label{eq:eqn4.5}
\frac{\partial {\bf C}}{\partial t}
+ \frac{\partial {\bf F}}{\partial x} + \frac{\partial {\bf G}}{\partial y} = 0.
\end{equation}
To derive the component-wise schemes, the relaxation system is sequentially diagonalized in one dimension first
\begin {subequations}
\label {eqn4.6}
\begin{eqnarray}
\frac{\partial {\bf C}}{\partial t}
&+& \frac{\partial {\bf V}}{\partial x} = 0 ,
\nonumber\\[-1.5ex]
\label{eq:eqn4.6a}  \\[-1.5ex]
\frac{\partial {\bf V}}{\partial t}
&+& {\bf A}_{X}^{2} \frac{\partial {\bf C}}{\partial x}
= \frac{1}{\varepsilon} \left({\bf F}\left({\bf C}\right) - {\bf V}\right) ,
\nonumber
\end{eqnarray}
and then in the other dimension
\begin{eqnarray}
\frac{\partial {\bf C}}{\partial t}
&+& \frac{\partial {\bf W}}{\partial y} = 0 ,
\nonumber\\[-1.5ex]
\label{eq:eqn4.6b}  \\[-1.5ex]
\frac{\partial {\bf W}}{\partial t}
&+& {\bf A}_{Y}^{2} \frac{\partial {\bf C}}{\partial y}
= \frac{1}{\varepsilon} \left({\bf G}\left({\bf C}\right) - {\bf W}\right),
\nonumber
\end{eqnarray}
\end {subequations}
in each time step, to obtain the characteristic variables $\left\{{\rm V}_{p} \pm a_{p}^x {\rm C}_{p} \right\}$
and $\left\{{\rm W}_{p} \pm a_{p}^y {\rm C}_{p} \right\}$. The relaxed updates are then obtained by
upwinding on these characteristic variables dimension-wise and setting
${\bf V}={\bf F}\left({\bf C}\right)$ and ${\bf W}={\bf G}\left({\bf C}\right)$.
So, a first order, component-wise, semi-discrete update is given as
\begin {subequations}
\label {eqn4.7}
\begin{equation}
\label{eq:eqn4.7a}
\frac{\partial C_{p,i,j} }{\partial t}
= - \frac{1}{\Delta x}
\left( \mathcal{F}_{p,i+{\frac{1}{2}},j} - \mathcal{F}_{p,i-{\frac{1}{2}},j} \right)
- \frac{1}{\Delta y}
\left( \mathcal {G}_{p,i,j+{\frac{1}{2}}} - \mathcal{G}_{p,i,j-{\frac{1}{2}}} \right),
\end{equation}
with
\begin{eqnarray}
\qquad \mathcal{F}_{p,i-{\frac{1}{2}},j}
= \frac{1} {2} \left[
\left( F_{p,i,j} + F_{p,i-1,j} \right)-a_{p}^{x} \left(C_{p,i,j} -C_{p,i-1,j} \right)
\right] ,
\nonumber\\[-1.5ex]
\label{eq:eqn4.7b} \\[-1.5ex]
\qquad \mathcal{G}_{p,i,j-{\frac{1}{2}}}
=  \frac{1} {2} \left[
\left( G_{p,i,j} + G_{p,i,j-1} \right)-a_{p}^{y} \left(C_{p,i,j} - C_{p,i,j-1} \right)
\right].
\nonumber
\end{eqnarray}
\end {subequations}
A second order, component-wise, semi-discrete update is given by
\begin {subequations}
\label {eqn4.8}
\begin{eqnarray}
\frac{\partial C_{p,i,j} }{\partial t}
&=& - \frac{1}{\Delta x}
\left( \mathcal{F}_{p,i+{\frac{1}{2}},j} - \mathcal{F}_{p,i-{\frac{1}{2}},j} \right)
- \frac{1}{\Delta x}
\left( \tilde{\mathcal{F}}_{p,i+{\frac{1}{2}},j} - \tilde{\mathcal{F}}_{p,i-{\frac{1}{2}},j} \right)
\nonumber\\[-1.5ex]
\label{eq:eqn4.8a} \\[-1.5ex]
&-& \frac{1}{\Delta y}
\left( \mathcal {G}_{p,i,j+{\frac{1}{2}}} - \mathcal{G}_{p,i,j-{\frac{1}{2}}} \right)
- \frac{1}{\Delta y}
\left( \tilde{\mathcal {G}}_{p,i,j+{\frac{1}{2}}} - \tilde{\mathcal{G}}_{p,i,j-{\frac{1}{2}}} \right),
\nonumber
\end{eqnarray}
where the fluxes ${\mathcal{F}_{p,i-{\frac{1}{2}},j}}$  and ${\mathcal{G}_{p,i,j-{\frac{1}{2}}}}$  are as
given in (\ref{eq:eqn4.7b}) and
\begin {eqnarray}
\label{eq:eqn4.8b}
\tilde {\mathcal F}_{p,i+1/2,j}
= \frac{\Delta x}{4} \left( {\sigma} _{p,i,j}^{x,+} - {\sigma} _{p,i+1,j}^{x,-} \right), \;\;
\tilde {\mathcal G}_{p,i,j+1/2}
= \frac{\Delta y}{4} \left( {\sigma} _{p,i,j}^{y,+} - {\sigma} _{p,i,j+1}^{y,-} \right),
\end{eqnarray}
with the slopes in ${x}$ and ${y}$ dimensions defined as
\begin{eqnarray}
\sigma _{p,i,j}^{x,\pm }
&=& \frac{1}{\Delta x} \left(
\left(F_{p} \pm a_{p}^x C_{p} \right)_{i+1,j} -\left(F_{p} \pm a_{p}^{x} C_{p} \right)_{i,j}
\right) \phi \left(\theta _{p,i,j}^{x,\pm } \right),
\nonumber\\[-1.5ex]
\label{eq:eqn4.8c}\\ [-1.5ex]
\sigma _{p,i,j}^{y,\pm }
&=& \frac{1}{\Delta y} \left(
\left(G_{p} \pm a_{p}^{y} C_{p} \right)_{i,j+1} -\left(G_{p} \pm a_{p}^{y} C_{p} \right)_{i,j}
\right)  \phi \left(\theta _{p,i,j}^{y,\pm } \right),
\nonumber
\end{eqnarray}
and with the limiting parameters
\begin{eqnarray}
\theta _{p,i,j}^{x,\pm }
&=&\frac
{ \left(F_{p} \pm a_{p}^x C_{p} \right)_{i,j}  -\left(F_{p} \pm a_{p}^x C_{p} \right)_{i-1,j} } 
{ \left(F_{p} \pm a_{p}^x C_{p} \right)_{i+1,j} -\left(F_{p} \pm a_{p}^x C_{p} \right)_{i,j} }, 
\nonumber\\[-1.5ex]
\label{eq:eqn4.8d} \\[-1.5ex]
\theta _{p,i,j}^{y,\pm }
&=&\frac
{ \left(G_{p} \pm a_{p}^y C_{p} \right)_{i,j} -\left(G_{p} \pm a_{p}^y C_{p} \right)_{i,j-1} } 
{ \left(G_{p} \pm a_{p}^y C_{p} \right)_{i,j+1} -\left(G_{p} \pm a_{p}^y C_{p} \right)_{i,j} }. 
\nonumber
\end{eqnarray}
\end{subequations}
We note that the above updates (\ref{eqn4.7}) and (\ref{eqn4.8}) correspond to
applying the 1D numerical fluxes (\ref{eq:eqn2.8b}, \ref{eq:eqn2.9c}) in a
dimension-by-dimension fashion.

\subsection {Variable relaxation in multidimensions}
An \textit{m}-dimensional variable relaxation system can be formulated as
\begin{eqnarray}
\frac{\partial {\bf C}}{\partial t}
&+& \sum _{k=1}^{m}\frac{\partial {\bf V}_{k} }{\partial x_{k} }  = 0,
\hspace{.3in}  {\bf C}{\bf ,V}_{k} \in R^{N}
\nonumber\\[-1.5ex]
\label{eq:eqn4.9}  \\[-1.5ex]
\frac{\partial {\bf V}_{k} }{\partial t}
&+&  {\bf A}_{k}({\bf x},t)^{2} \frac{\partial {\bf C}}{\partial x_{k} }
= \frac{1}{\varepsilon } \left( {\bf F}_{k} \left ({\bf C} \right) - {\bf V}_{k} \right),
\hspace{.3in} k=1,2,\cdots m
\nonumber .
\end{eqnarray}
This system, when expressed as a dissipative approximation to the original system (\ref{eq:eqn4.1}),
has the stability requirement
$\left(\partial _{k\tilde{k}} {\bf A}_{k}({\bf x},t)^{2} -{\bf F}_{k} {\bf F}_{\tilde{k}} \right) \ge 0$,
$\forall {\rm \; } {\bf C} \left({\bf x}, t\right)$. For symmetric subcharacteristic speeds
(i.e. $a_{p}^{k,+}=-a_{p}^{k,-}=a_{p}^k $) the subcharacteristic condition can be represented
in terms of the local speeds ${\left\{ \lambda _{p}^k({\bf x},t) \right\}}$ of the
Jacobian ${\bf F}_{k} \left ({\bf C} \left( {\bf x}, t\right)\right)$
and the local subcharacteristic speeds $ \left\{a _{p}^k({\bf x},t) \right\}$  as
\begin{equation}
\label{eq:eqn4.10}
\frac{\lambda _{^{1,\max } }({\bf x},t)^{2} }{a_{^{1,min} }({\bf x},t)^{2} }
+\frac{\lambda _{^{2,\max } }({\bf x},t)^{2} }{a_{^{2,min} }({\bf x},t)^{2} }
+\cdots
+\frac{\lambda _{^{m,\max } }({\bf x},t)^{2} }{a_{^{m,min}}({\bf x},t)^{2} } \le 1,
\end{equation}
where $\lambda _{k, \max }({\bf x},t) = \mathop{\max }\limits_{1\le p\le N} \left|\lambda _{p}^k({\bf x},t) \right|$ and
$a_{k, \min}({\bf x},t)  = \mathop{\min }\limits_{1\le p\le N} \left|a _{p}^k({\bf x},t) \right|$. The multidimensional
relaxed schemes can be obtained either by solving the relaxation system dimension-by-dimension or by
dimension-wise extension of the 1D fluxes (\ref{eq:eqn3.8b}, \ref{eq:eqn3.11b}). A
 first order, semi-discrete update with symmetric subcharacteristics can be written as
\begin {subequations}
\label {eqn4.11}
\begin{equation}
\label{eq:eqn4.11a}
\frac{\partial C_{p,i,j}}{\partial t}
= - \frac{1}{\Delta x}
\left( \mathcal{F}_{p,i+{\frac{1}{2}},j} - \mathcal{F}_{p,i-{\frac{1}{2}},j} \right)
- \frac{1}{\Delta y}
\left( \mathcal {G}_{p,i,j+{\frac{1}{2}}} - \mathcal{G}_{p,i,j-{\frac{1}{2}}} \right),
\end{equation}
with
\begin{eqnarray}
\qquad \mathcal{F}_{p,i-{\frac{1}{2}},j}
=  F_{p,i-1,j}
+ \frac{1} {2} \left[
\left( F_{p,i,j} - F_{p,i-1,j} \right)-a_{p,i-{\frac{1}{2}},j}^{x} \left(C_{p,i,j} - C_{p,i-1,j} \right)
\right],
\nonumber\\[-1.5ex]
\label{eq:eqn4.11b} \\[-1.5ex]
\qquad \mathcal{G}_{p,i,j-{\frac{1}{2}}}
=  G_{p,i,j-1}
+ \frac{1} {2} \left[
\left( G_{p,i,j} - G_{p,i,j-1} \right) - a_{p,i,j-{\frac{1}{2}}}^{y} \left(C_{p,i,j} -C_{p,i,j-1} \right)
\right].
\nonumber
\end{eqnarray}
\end {subequations}
A 2D second order scheme is given by
\begin {subequations}
\label {eqn4.12}
\begin{eqnarray}
\frac{\partial C_{p,i,j}}{\partial t}
&=& - \frac{1}{\Delta x}
\left( \mathcal{F}_{p,i+{\frac{1}{2}},j} - \mathcal{F}_{p,i-{\frac{1}{2}},j} \right)
 - \frac{1}{\Delta x}
\left( \tilde{\mathcal{F}}_{p,i+{\frac{1}{2}},j} - \tilde{\mathcal{F}}_{p,i-{\frac{1}{2}},j} \right)
\nonumber\\[-1.5ex]
\label{eq:eqn4.12a} \\[-1.5ex]
&-& \frac{1}{\Delta y}
\left( \mathcal {G}_{p,i,j+{\frac{1}{2}}} - \mathcal{G}_{p,i,j-{\frac{1}{2}}} \right)
- \frac{1}{\Delta y}
\left( \tilde{\mathcal {G}}_{p,i,j+{\frac{1}{2}}} - \tilde{\mathcal{G}}_{p,i,j-{\frac{1}{2}}} \right),
\nonumber
\end{eqnarray}
with the dimension-wise higher order corrections
\begin{eqnarray}
\label{eq:eqn4.12b}
\qquad \qquad \tilde{\mathcal{F}}_{p,i-{\frac{1}{2}},j}
&=& \frac {\phi \left(\theta _{p,i-{\frac{1}{2}},j}^{x,+} \right)} {4}
\left[ a_{p,i-{\frac{1}{2}},j}^{x} \left(C_{p,i,j} -C_{p,i-1,j} \right)
+ \left(F_{p,i,j} -F_{p,i-1,j} \right) \right]
\\
& + &
\frac{ \phi \left(\theta _{p,i-{\frac{1}{2}},j}^{x,-} \right)} {4}
\left[ a_{p,i-{\frac{1}{2}},j}^{x} \left(C_{p,i,j} - C_{p,i-1,j}  \right)
-\left(F_{p,i,j} - F_{p,i-1,j} \right) \right],
\nonumber \\
\tilde{\mathcal{G}}_{p,i,j-{\frac{1}{2}}}
&=& \frac {\phi \left(\theta _{p,i,j-{\frac{1}{2}}}^{y,+} \right)} {4}
\left[ a_{p,i,j-{\frac{1}{2}}}^{y} \left(C_{p,i,j} -C_{p,i,j-1} \right)
+ \left(F_{p,i,j} -F_{p,i,j-1} \right) \right]
\nonumber\\
& + &
\frac{ \phi \left(\theta _{p,i,j-{\frac{1}{2}}}^{y,-} \right)} {4}
\left[ a_{p,i,j-{\frac{1}{2}}}^{y} \left(C_{p,i,j} - C_{p,i,j-1}  \right)
-\left(F_{p,i,j} - F_{p,i,j-1} \right) \right],
\nonumber
\end{eqnarray}
and the dimension-wise limiting parameters
\begin{eqnarray}
\theta _{p,i-{\frac{1}{2}},j}^{x,+}
&=& \frac{
\frac  {1 + a_{p,i-{\frac{1}{2}},j}^{x} a_{p,i-{\frac{3}{2}},j}^{x} }
{2a_{p,i-{\frac{3}{2}},j}^{x}}
\left[a_{p,i-{\frac{3}{2}},j}^{x}\left(C_{p,i-1,j} -C_{p,i-2,j}\right)
+ \left( F_{p,i-1,j} -F_{p,i-2,j} \right)  \right]}
{
\frac{1 + a_{p,i-{\frac{1}{2}},j}^{x} a_{p,i-{\frac{1}{2}},j}^{x}}
{2a_{p,i-{\frac{1}{2}},j}^{x}}
\left[ a_{p,i-{\frac{1}{2}},j}^{x} \left(C_{p,i,j} - C_{p,i-1,j} \right)
+ \left( F_{p,i,j} - F_{p,i-1,j} \right) \right]
} ,
\nonumber \\[-1.5ex]
\label{eq:eqn4.12c} \\ [-1.5ex]
\theta _{p,i-{\frac{1}{2}},j}^{x,-}
&=& \frac{
\frac  {1+a_{p,i+{\frac{1}{2}},j}^{x} a_{p,i-{\frac{1}{2}},j}^{x} }
{2a_{p,j+{\frac{1}{2}}}^{x} }
\left[
a_{p,i+{\frac{1}{2}},j}^{x} \left(C_{p,i+1,j} -C_{p,i,j} \right)
-  \left( F_{p,i+1,j} - F_{p,i,j}   \right)
\right]
}
{
\frac{1+a_{p,i-{\frac{1}{2}},j}^{x} a_{p,i-{\frac{1}{2}},j}^{x} }
{2a_{p,i-{\frac{1}{2}},j}^{x} }
\left[ a_{p,i-{\frac{1}{2}},j}\left(C_{p,i,j} - C_{p,i-1,j} \right)
 -  \left( F_{p,i,j} - F_{p,i-1,j}  \right)   \right]
},
\nonumber \\
\theta _{p,i,j-{\frac{1}{2}}}^{y,+}
&=& \frac{
\frac  {1 + a_{p,i,j-{\frac{1}{2}}}^{y} a_{p,i,j-{\frac{3}{2}}}^{y} }
{2a_{p,i,j-{\frac{3}{2}}}^{y} }
\left[a_{p,i,j-{\frac{3}{2}}}^{y} \left(C_{p,i,j-1} - C_{p,i,j-2} \right)
+ \left( F_{p,i,j-1} - F_{p,i,j-2} \right)  \right]}
{
\frac{1 + a_{p,i,j-{\frac{1}{2}}}^{y} a_{p,i,j-{\frac{1}{2}}}^{y} }
{2a_{p,i,j-{\frac{1}{2}}}^{y} }
\left[ a_{p,i,j-{\frac{1}{2}}^{y} } \left(C_{p,i,j} - C_{p,i,j-1} \right)
+ \left( F_{p,i,j} - F_{p,i,j-1} \right) \right]
},
\nonumber \\
\theta _{p,i,j-{\frac{1}{2}}}^{y,-}
&=&  \frac{
\frac  {1+a_{p,i,j+{\frac{1}{2}}}^{y} a_{p,i,j-{\frac{1}{2}}}^{y} }
{2a_{p,i,j+{\frac{1}{2}}}^{y} }
\left[
a_{p,i,j+{\frac{1}{2}}}^{y} \left(C_{p,i,j+1} -C_{p,i,j} \right)
-  \left( F_{p,i,j+1} - F_{p,i,j}  \right)
\right]
}
{
\frac{1+a_{p,j-{\frac{1}{2}}}^{y} a_{p,j-{\frac{1}{2}}}^{y} }
{2a_{p,j-{\frac{1}{2}}}^{y} }
\left[ a_{p,j-{\frac{1}{2}}}^{y}\left(C_{p,j} - C_{p,j-1} \right)
 -  \left( F_{p,j} - F_{p,j-1}  \right)   \right]
}.
\nonumber
\end{eqnarray}
\end {subequations}
The variable relaxed updates with asymmetric, optimal speeds can be obtained similarly by extending the 1D fluxes
(\ref{eq:eqn3.8b}, \ref{eq:eqn3.11b}). But choosing optimal subcharacteristic speeds in multidimensions is nontrivial
because  the subcharacteristic condition for the asymmetric case does not reduce to a simple form as in
(\ref {eq:eqn4.10}). For the 2D case however, one can choose the optimal speeds as
\begin{eqnarray}
\label{eq:eqn4.13}
&a_{p}^{k,+}(x,y,t)&
\ge 2\max \left(\max \limits_{1\le p\le N} \lambda _{p}^k\left( {\bf C} \right) , 0 \right) \\
&a_{p}^{k,-}(x,y,t)&
\ge 2\min \left( \min \limits_{1\le p\le N} \lambda _{p}^k\left( {\bf C} \right) , 0 \right)
\nonumber, \;\; k=x,y.
\end{eqnarray}
The relaxed schemes in this section were based on \textit{donor cell} \cite{lvq02} upwinding on
the charactersitic variables.  Truly multidimensional upwinding methods like
\textit{corner transport upwind} (CTU) \cite{clela} or the  Flat scheme \cite{kmg} can of course also be used.

\section{Numerical Results}
In this section, we demonstrate the results of our new higher order schemes on several problems.
We use a 1D Burgers equation to test convergence: for the constant coeffcient linear
advection problem the variable relaxed schemes reduce to JX scheme (symmetric case) or upwind scheme (optimal case).
Next, we show the results on a weakly hyperbolic, two-phase gas-oil displacement in 1D and 2D. We also present
results for the weakly hyperbolic,  single phase geometric optics problem introduced by
Engquist and Runborg \cite{er}.
While all the relaxation schemes behave well in the presence of weak hyperbolicity, the true advantage of our schemes
becomes apparent when there is a large variation of speeds in the domain.
In all the examples the JX subcharacteristic speeds are set at maximum eigenvalue of the original problem.
The  first order semi-discrete updates (\ref{eqn2.8}) and (\ref{eqn3.8})  are used with
forward Euler time stepping and the second order semi-discrete schemes (\ref{eqn2.9}) and (\ref{eqn3.11}) are used with
a 2-stage TVD Runge-Kutta time-stepping.
For brevity, throughout this section the schemes are abbrevated as: JX for Jin-Xin scheme,
VRS for variable relaxed scheme with symmetric speeds and VRO for variable relaxed scheme with optimal speeds.

\subsection {Burgers Equation}
\begin{figure}[h]
\begin{center}
\hspace{-0.4in}
\subfigure{
\includegraphics[bb=0mm 0mm 208mm 296mm, angle = 90, scale=0.25]{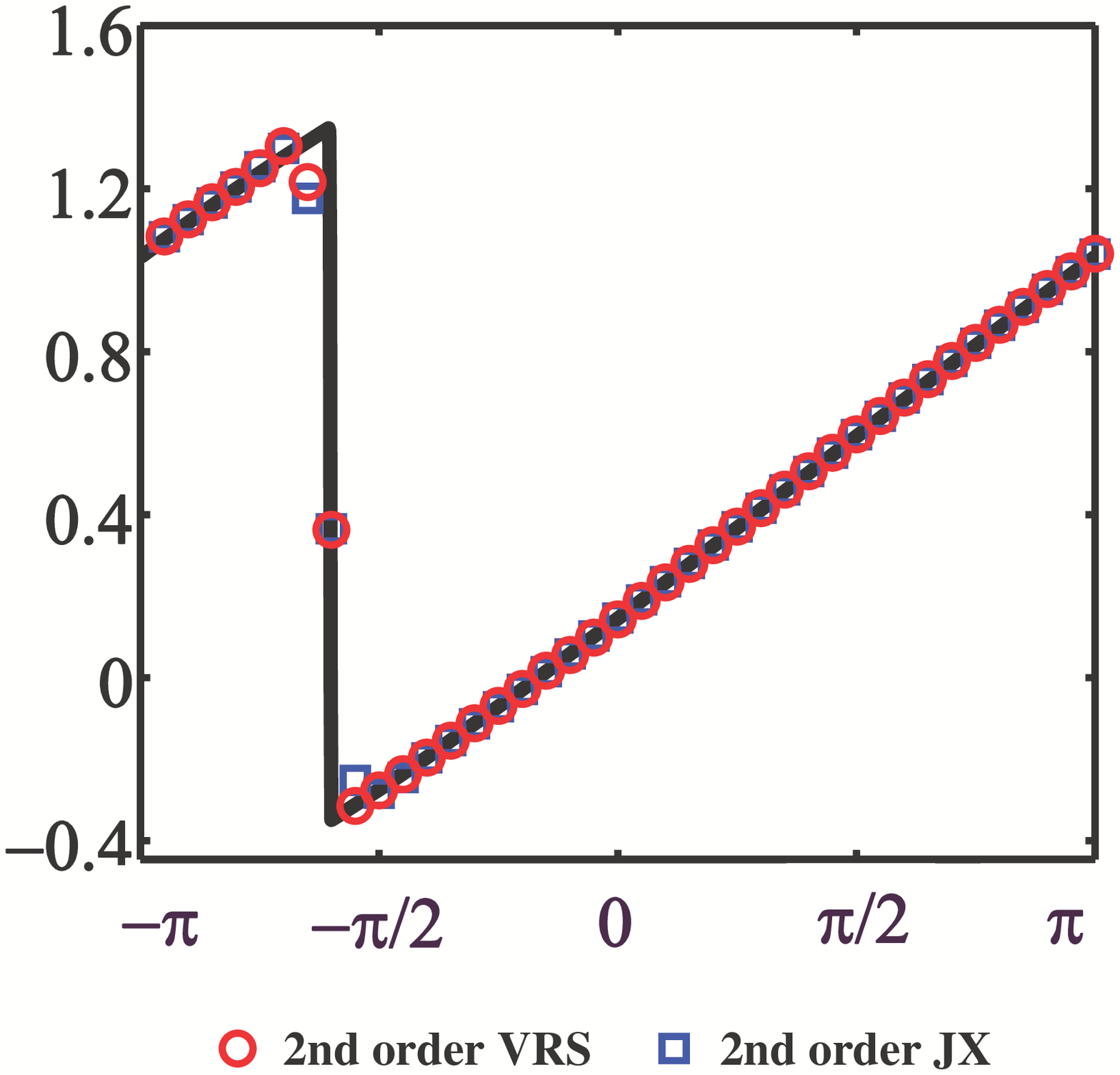}
}
\hspace{-0.75in}
\subfigure{
\includegraphics[bb=0mm 0mm 208mm 296mm, angle = 90, scale=0.25]{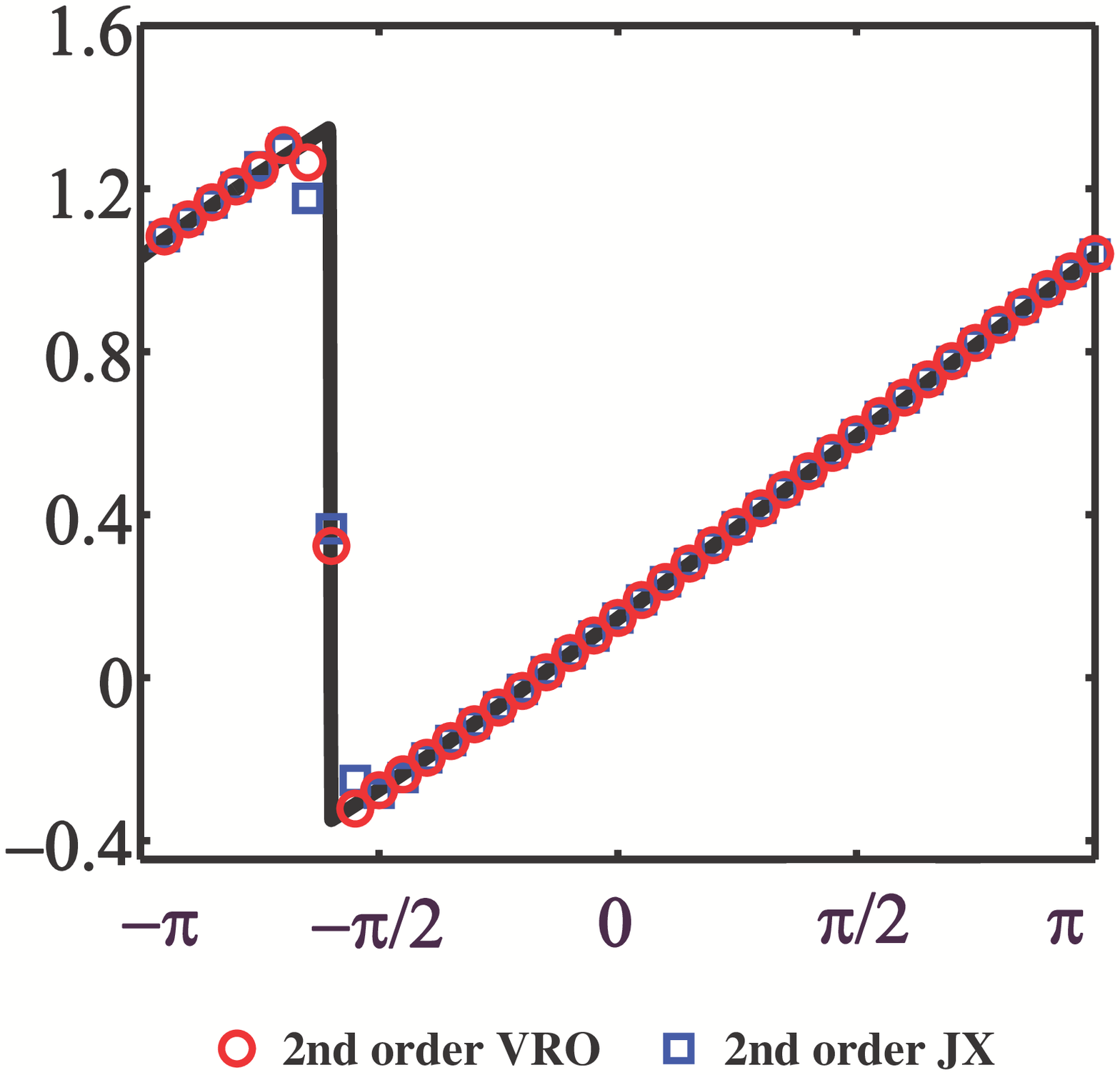}
}
\vspace{-0.5in}
\caption{ Post-shock solutions of 1D Burgers equation (\ref{eq:eqn5.1}) ${T=2.5}$, ${\Delta x = 1/40}$,
JX subcharacteristic speed = Max speed = 1.5, ${\Delta t = 0.5\frac{\Delta x}{1.5}}$,
}
\label{fig5.1}
\end{center}
\end{figure}

We test the order of accuracy of our higher order schemes on the pre-shock solutions of Burgers equation
with periodic initial data
\begin{equation}
\label{eq:eqn5.1}
\centering
C_{t} +\left(\frac{C^{2}}{2} \right)_{x} =0, \qquad C(x,0)= 0.5 + \sin x, \qquad x \in [-\pi,\pi].
\end{equation}
The solution of (\ref{eq:eqn5.1})  develops a shock at ${T_c = 1}$. The
${L_1}$ and ${ L_{\infty }}$ errors and the measured order of accuracy
 are given in table (\ref{table5.1}) at time $T=0.5$ when the solution is still smooth.
\begin{table}[h]
 \caption{${L_1}$ and ${ L_{\infty }}$ error and order for Burgers equation (\ref{eq:eqn5.1}), ${T=0.5}$}
\begin{center} \footnotesize
\begin{tabular}{|p{0.4in}|p{0.6in}|p{0.6in}|p{0.6in}|p{0.6in}|p{0.6in}|p{0.6in}|} \hline
N  & \multicolumn{2}{|p{1.2in}|}{JX ${\rm 2^{nd}}$ order}
& \multicolumn{2}{|p{1.2in}|}{VRS ${\rm 2^{nd}}$ order}
& \multicolumn{2}{|p{1.2in}|}{VRO ${\rm 2^{nd}}$ order}\\ \hline
  & ${L_1}$ error & ${L_1}$ order & ${L_1}$ error & ${L_1}$ order & ${L_1}$ error & ${L_1}$ order \\ \hline
20  & 6.0087e-2  & - & 4.9486e-2   & - & 5.0998e-2   & - \\ \hline
40 & 1.7527e-2 & 1.7778 & 1.3429e-2 & 1.8817 & 1.3439e-2 & 1.9240 \\ \hline
80  & 4.9088e-3  & 1.8358 & 3.9125e-3 & 1.7791 & 3.7378e-3 & 1.8462 \\ \hline
160  & 1.2463e-3 & 1.9777 & 1.0257-3 & 1.9315 & 9.6010-4 & 1.9609 \\ \hline
320  & 3.040e-4 & 2.0357 & 2.645e-4 & 1.9551 & 2.3800e-4 & 2.0124 \\ \hline
N & \multicolumn{2}{|p{1.2in}|}{JX ${\rm 2^{nd}}$ order}
& \multicolumn{2}{|p{1.2in}|}{VRS ${\rm 2^{nd}}$ order}
& \multicolumn{2}{|p{1.2in}|}{VRO ${\rm 2^{nd}}$ order} \\ \hline
& ${L_{\infty}}$ error & ${L_{\infty}}$ order & ${L_{\infty}}$ error & ${L_{\infty}}$ order
& ${L_{\infty}}$ error & ${L_{\infty}}$ order\\ \hline
20  & 3.4612e-2   & - & 3.4623e-2 & - & 3.4559e-2 & - \\ \hline
40  & 1.3580e-2  & 1.3498 & 1.3570e-2 & 1.3512 & 1.3554e-2 & 1.3504 \\ \hline
80  & 5.2178e-3   & 1.3799 & 5.2153e-3 & 1.3797 & 5.2124e-3 & 1.3787 \\ \hline
160  & 1.9775e-3      & 1.3998 & 1.9771e-3 & 1.3993 & 1.9768e-3 & 1.3988 \\ \hline
320  & 7.4200e-4    & 1.4142 & 7.4190e-4 & 1.4141 & 7.4190e-4 & 1.4139 \\ \hline
\end{tabular}
\end{center}
\label{table5.1}
\end{table}
The variable relaxed schemes in general have smaller errors than the JX schemes.
As expected, the relaxation with optimal speeds yields the best result. The post shock solutions are
shown in figure (\ref{fig5.1}), where again the results with variable relaxation are slightly better
than that of the JX scheme.
\subsection {Gas-oil displacements}
\label{gasinjection}

We consider a two-phase (vapor and oil) gas injection system. The vapor and oil phases consist of mixtures of $n_c$ (hydrocarbon) components with overall volume fraction $C_i$. The system is governed by a pressure equation, which determines the phase velocities in the domain, and transport equations (mass balance equations), one for each component. A complete discussion of gas injection processes can be found in \cite{orr}. Here, we consider only gas injection systems with simplified phase behavior and ignore capillary pressure and gravity for ease of presentation. We also assume that the fluid system is incompressible, that is $\nabla . {\bf u}_T=0$, where ${\bf u}_T$ is the total velocity in the multidimensional system.
We use Darcy's law to express the total velocity as a function of the pressure gradient in the domain and phase mobilities. We then obtain
\begin {subequations}
\label{eqn5.2}
\begin{equation}
\label{eq:eqn5.2a}
\phi \frac{\partial C_{i} }{\partial t} +\nabla .\left({\bf u}_{T} F_{i} \right)=0, \qquad  i=1,\cdots,n_c,
\end{equation}
\begin{equation}
\label{eq:eqn5.2b}
\nabla .{\bf u}_{T} =
\nabla .\left( {\bf k}
\left(\frac{k_{rV} \left(S \right)}{\mu _{V}}
 + \frac{k_{rL} \left(S \right)}{\mu _{L}} \right)
\nabla P
 \right)
= 0,
\end{equation}
\end {subequations}
where ${S}$ is the saturation of vapor phase, ${1-S}$ the saturation of oil phase,
 ${\bf k}$ is the permeability tensor and
${\nabla P}$  the phase pressure gradient. The vapor and liquid relative phase permeabilities
${k_{rV} \left(S\right)}$ and ${k_{rL}\left(S\right)}$ are taken to be
\begin{equation}
\label{eq:eqn5.3}
\begin{array}{cc}
&k_{rV}\left(S\right)=0, \;\; k_{rL}\left(S\right)=1, \qquad \qquad \qquad \; S < S_{gc}, \\ \\
&k_{rV}\left(S\right)=\frac{\left(S-S_{gc} \right)^{2} }{\left(1-S_{gc}-S_{or} \right)^{2} }, \;\;\;
k_{rL}\left(S\right)=\frac{\left(1-S-S_{or} \right)^{2}}{\left(1-S_{gc}-S_{or} \right)^{2} },
\qquad S_{gc} < S < 1-S_{or}  \\ \\
&k_{rV}\left(S\right)=1,\;\; k_{rL}\left(S\right)=0, \qquad \qquad \qquad \; S > 1-S_{or}.
\end{array}
\end{equation}
Also, ${\phi}$ is the porosity (volume fraction of the void space).
$C_{i} $ is  the overall volume fraction of component-$i$, which is given by
$C_{i} = c_{iV} S + c_{iL} \left(1-S\right) $,
with  ${c_{iV}}$ and ${c_{iL}}$  being
 the volume fractions of component-${i}$ in vapor and liquid phases.
${F_{i}}$ is the overall fractional volumetric flow of component-$i$ given as
$F_{i} =c_{iV} f + c_{iL} \left(1-f\right) $, and $f$ is the vapor  fractional flow given as
\begin{equation}
\label{eq:eqn5.4}
f=\frac{k_{rV} (S)}{k_{rV} (S)+Mk_{rL} (S)} ,
\end{equation}
where $M$ is the constant viscosity ratio ${ \frac {\mu _{V}}{\mu _{L}}}$.
The phase compositions  ${ c_{iV}}$ and  ${c_{iL}}$ are related as
\begin{equation}
\label{eq:eqn5.5}
K_{i} =\frac{c_{iV} }{c_{iL} }
\end{equation}
where the K-values ${K_{i}}$  are assumed to be constant. When numerically solving, the saturation $S$ and phase
compositions ${ c_{iV}}$ and  ${c_{iL}}$ are usually obtained by performing iterative phase equilibrium
calculations using an equation-of-state (van der Waals \cite{vander} or Peng-Robinson \cite{pr}). Here, for simple
phase behavior, the saturation and phase compositions are obtained by iteratively solving  \cite{orr}
\begin{equation}
\label{eq:eqn5.6}
\sum _{i=1}^{n_{c}} c_{iV}-c_{iL}
=\sum _{i=1}^{n_{c}} \frac{C_i\left(K_i \right)}{1 + S\left(K_i -1 \right)}=0.
\end{equation}
Under the constraint that the volume fractions must sum to 1,
\begin{equation}
\label{eq:eqn5.7}
\sum _{i=1}^{n_{c}}C_{i} =1, \;\;  \sum _{i=1}^{n_{c}}c_{iV} =1, \;\; \textnormal{and}\;\; \sum _{i=1}^{n_{c}}c_{iL} =1,
\end{equation}
the problem (\ref{eqn5.2}) can be expressed in terms of conservation of $n_{c} -1$ components only.
A gas injection problem with two components reduces to a scalar conservation
law (a generalized Buckley-Leverett problem), and
a three-component gas injection system reduces to a 2x2 nonlinear conservation system. In multicomponent systems generally the
lightest components are represented in the system of equations. In this work, we use three-component examples as they are the
simplest multicomponent systems that exhibit weak hyperbolicity.
\begin{figure}[h!]
\begin{center}
\hspace {-0.5in}
\includegraphics[bb=0mm 0mm 208mm 296mm, angle = 90, scale=0.4]{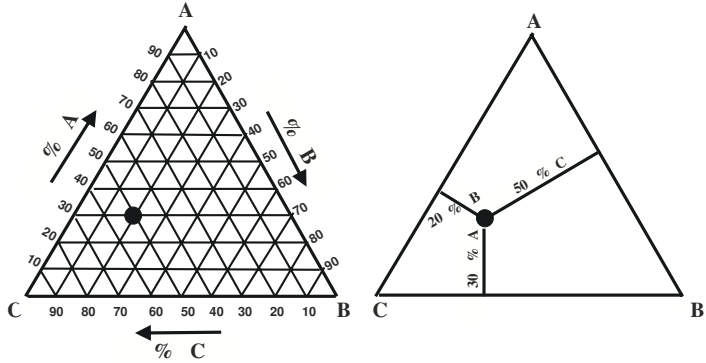}
\vspace{-0.75in}
\caption{ Ternary phase diagram}
\label{fig5.2}
\end{center}
\end{figure}

Ternary systems can be nicely represented by   ternary phase diagrams, which
display phase behavior information at fixed pressure and temperature. An example is given in figure \ref{fig5.2}.
Phase diagrams represent the component concentrations of all possible mixtures of the three components in a two-dimensional space.
Because the volume fractions of the three
components sum to one, the phase compositions can be conveniently represented on an equilateral triangle \cite{lake89,orr}.
Each vertex represents 100 \%
of the component associated with that vertex, and the side opposite 0\%.
Each point within the triangle represents a mixture of the three components; the volume fractions are read from the
perpendicular distance from that point to the three sides of the triangle. For gas/oil systems, the component associated
with the top vertex of the triangle is usually the lightest, and the component associated with the bottom left vertex is
usually the heaviest.

\subsubsection{1D ternary example}
In one dimension,  the transport and the pressure equations (\ref{eqn5.2})  become
\begin {subequations}
\label {eqn5.8}
\begin{eqnarray}
\phi \frac{\partial C_{i}}{\partial t}
+\frac{\partial \left(u_{\textsc{t}} F_{i} \right)}{\partial x} &=& 0, \qquad  \textnormal{for} \; i=1,2,\cdots n_{c} -1,
\label{eq:eqn5.8a} \\
\frac{\partial u_{\textsc{T}}}{\partial x} &=& 0 .
\label{eq:eqn5.8b}
\end{eqnarray}
\end {subequations}
Since the total velocity $u_{\textsc{T}}$ is constant, (\ref{eq:eqn5.8a}) can be conveniently expressed in dimensionless form
\begin{equation}
\label{eq:eqn5.9}
\frac{\partial C_{i} }{\partial {\tilde t}} +\frac{\partial F_{i} }{\partial {\tilde x}} =0
 \qquad \textnormal{for} \; i=1,2,\cdots n_{c} -1 ,
\end{equation}
where ${\tilde t}$  and ${\tilde x}$ are now the dimensionless time and spatial variables given by
\begin{equation}
\nonumber
{\tilde t} = \frac {u_{\textsc{T} t}}{\phi L}, \qquad {\tilde x} = \frac{x}{L},
\end{equation}
where ${L}$  is the length of the domain.
For ease of notation, we will henceforth denote these dimensionless variables
${\tilde t}$  and ${\tilde x}$ by ${t}$  and ${x}$. The analytical theory of 1D diffusion-free,
two phase gas-oil displacements is described in \cite{orr}. An
advection system with constant initial and injection conditions, with simple phase behavior
as given above, can be solved by the method of characteristics (MOC).

Consider a ternary example, described by the conservation of two lightest components $C_1$ and $C_2$
\begin{eqnarray}
\label{eq:eqn5.10}
\frac{\partial C_{1} }{\partial t} &+& \frac{\partial F_{1} }{\partial x} =0,
\\
\frac{\partial C_{2} }{\partial t} &+& \frac{\partial F_{2} }{\partial x} =0,
\nonumber
\end{eqnarray}
where
\begin{eqnarray}
C_{1} = c_{1V} S+c_{1L} (1-S),
\qquad  &F_{1}\left(C_{1}, C_{2}\right) =c_{1V} f(S)+c_{1L} (1-f(S)),
\nonumber  \\
C_{2} = c_{2V} S+c_{2L} (1-S),
\qquad &F_{2}\left(C_{1}, C_{2}\right) =c_{2V} f(S)+c_{2L} (1-f(S)),
\nonumber
\end{eqnarray}
with $f\left(S\right)$ as given by  (\ref{eq:eqn5.4}) and the constant K-values given by
 $ K_1 = 2.5, K_2 = 1.5, K_3 = 0.05$.
This system has two eigenvalues given by
\begin {subequations}
\label {eqn5.11}
\begin{eqnarray}
\label{eq:eqn5.11a}
\lambda _{t} &=& \frac{\partial F_{1} }{\partial C_{1} }
= \left\{\begin{array}{l} {{\frac{df}{dS}}
 {\rm \; \; \; \; \; \; in\; the\; two\; phase\; region}}
 \\
 {1{\rm \; \; \; \; \; \; \; \; in\; the\; single\; phase\; region}} \end{array}\right. ,
 \\  \label{eq:eq5.11b}
\lambda _{nt} &=&\frac{F_{1} +q}{C_{1} +q}, \; \; q=\frac{c_{1L}^{2} }{\gamma \left(K_1, K_2, K_3 \right) },
\end{eqnarray}
\end {subequations}
\begin{figure}[h]
\begin{center}
\begin{minipage}{2.6in}%
\hspace{-0.7in}
\includegraphics[bb=0mm 0mm 208mm 296mm, angle = 90, scale=0.3]{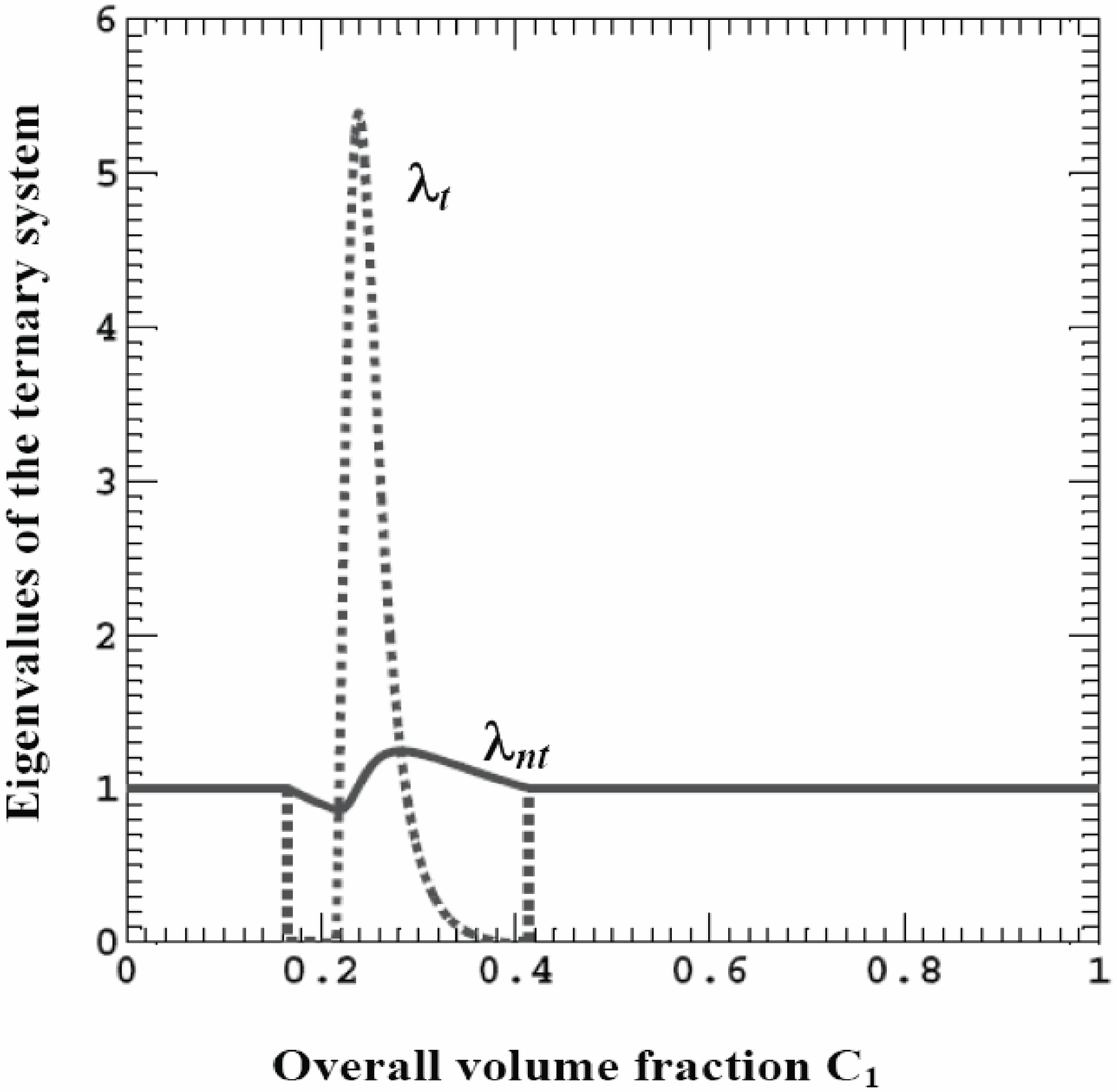}
\caption{Eigenvalues variations  for a \newline ternary system}
\label{fig5.3}
\end {minipage}%
\begin{minipage}{2.5in}%
\hspace{-0.59in}
\includegraphics[bb=0mm 0mm 208mm 296mm, angle = 90, scale=0.3]{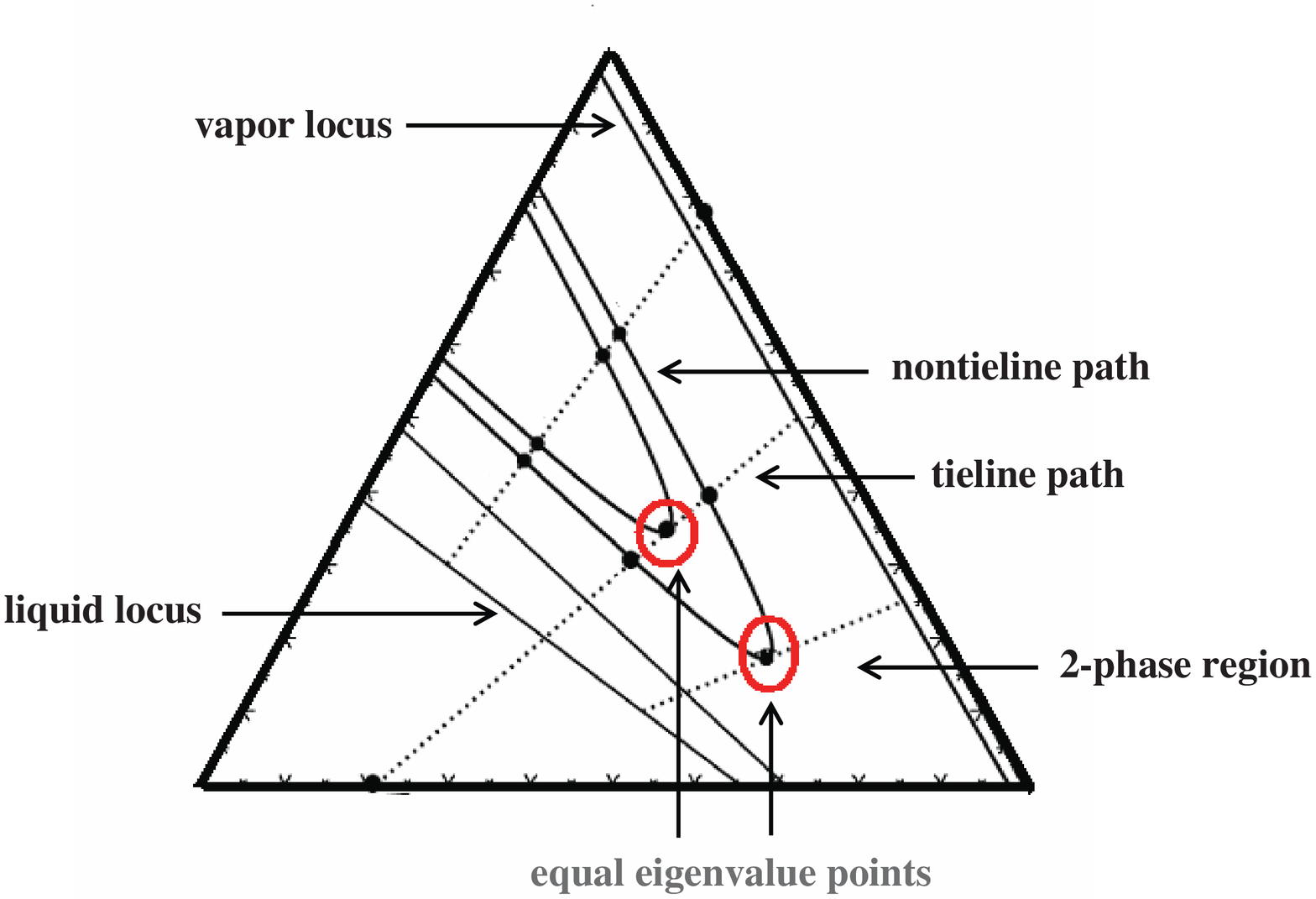}
\caption{Tieline and nontieline paths for a ternary system with constant K-values.}
\label{fig5.4}
\end {minipage}%
\end{center}
\end{figure}
where $\gamma \left(K_1, K_2, K_3 \right)
= \frac{\left(1 - K_3\right) \left(K_2 - 1\right)}{\left(K_1 - K_3\right) \left(K_1 - K_2\right)} $  (see \cite{orr}).
Figure (\ref{fig5.3}) shows the variation of the eigenvalues ${\lambda _{t}}$ and  ${\lambda _{nt}}$
 with $S_{or} =0.1$, $S_{gc} =0.2$ and  ${\frac{\mu _{V}}{\mu _{L}}= \frac{1}{20}} $.
The corresponding eigenvectors are
\begin{equation}
\label{eq:eqn5.12}
\vec{e}_{t} =\left[{\rm 1} \;\;\; {{\rm 0}}\right]^{\rm T}
\; \; \textnormal {and}  \;\;
\vec{e}_{nt} =\left[{\rm 1} \;\;\; {\frac{\lambda _{t} -\lambda _{nt}}
{{\frac{\partial F_{1} }{\partial c_{1L} }} } } \right]^{\rm T} .
\end{equation}
The eigenvectors correspond to the possible paths a solution can trace in  phase space.
The eigenvector $\vec{e}_{t}$ gives the straight line paths in the phase space, known as
the \textit{tie-line} paths (figure \ref{fig5.4}). The eigenvector $\vec{e}_{nt}$ gives the curved paths in the phase
space, known as the \textit{nontie-line} paths. Within the two-phase region only certain specific volume fractions of
liquid and vapor phase ($c_{iL} $ and $c_{iV} $) can be in equilibrium,
and each tie-line connects a pair of equilibrium volume fractions $c_{iL} $ and $c_{iV} $.
Tie-lines also connect the vapor locus and the liquid locus in the phase space. The point at
which a tie-line intersects the vapor locus has $S =1$ and
the point where it intersects the liquid locus has $S = 0$.

When the eigenvalues coincide, $\lambda _{t} =\lambda _{nt} $, so do the eigenvectors,
$\vec{e}_{t} =\vec{e}_{nt} $, and the system has dependent eigenvectors, i.e., the system becomes weakly hyperbolic.
The weak hyperbolicity is not limited only to ternary systems. With every additional component, there will be an additional
tieline-nontieline intersection, and hence an additional point of weak hyperboliciy.
\begin{figure}[h!]
\vspace {-0.2in}
\subfigure{
\begin{minipage}{2.8in}%
\hspace{-0.35in}
\includegraphics[bb=0mm 0mm 208mm 296mm, angle = 90, scale=0.25]{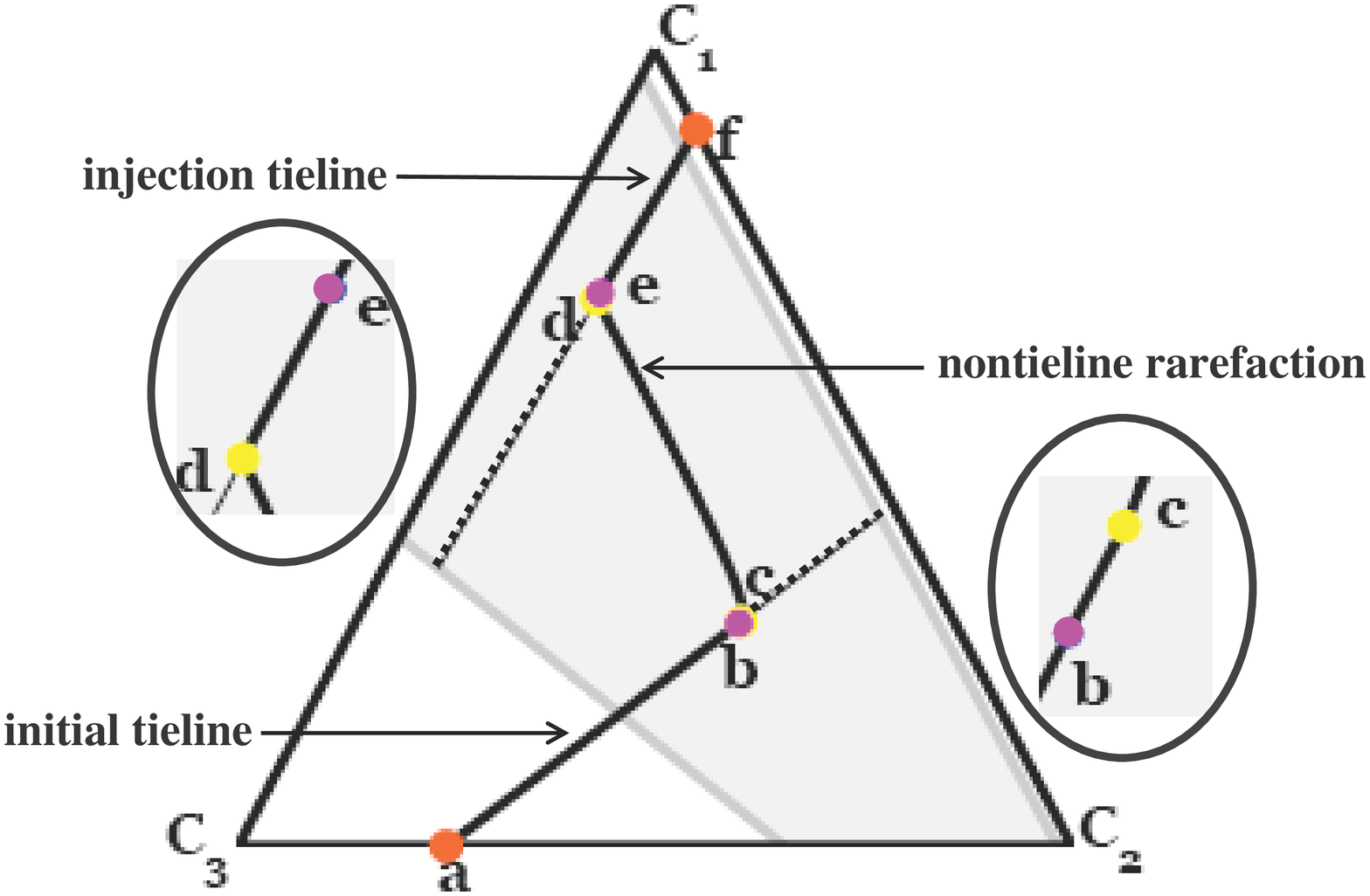}
\end {minipage}%
\begin{minipage}{2.5in}%
\hspace{-0.25in}
\includegraphics[bb=0mm 0mm 208mm 296mm, angle = 90,scale=0.2]{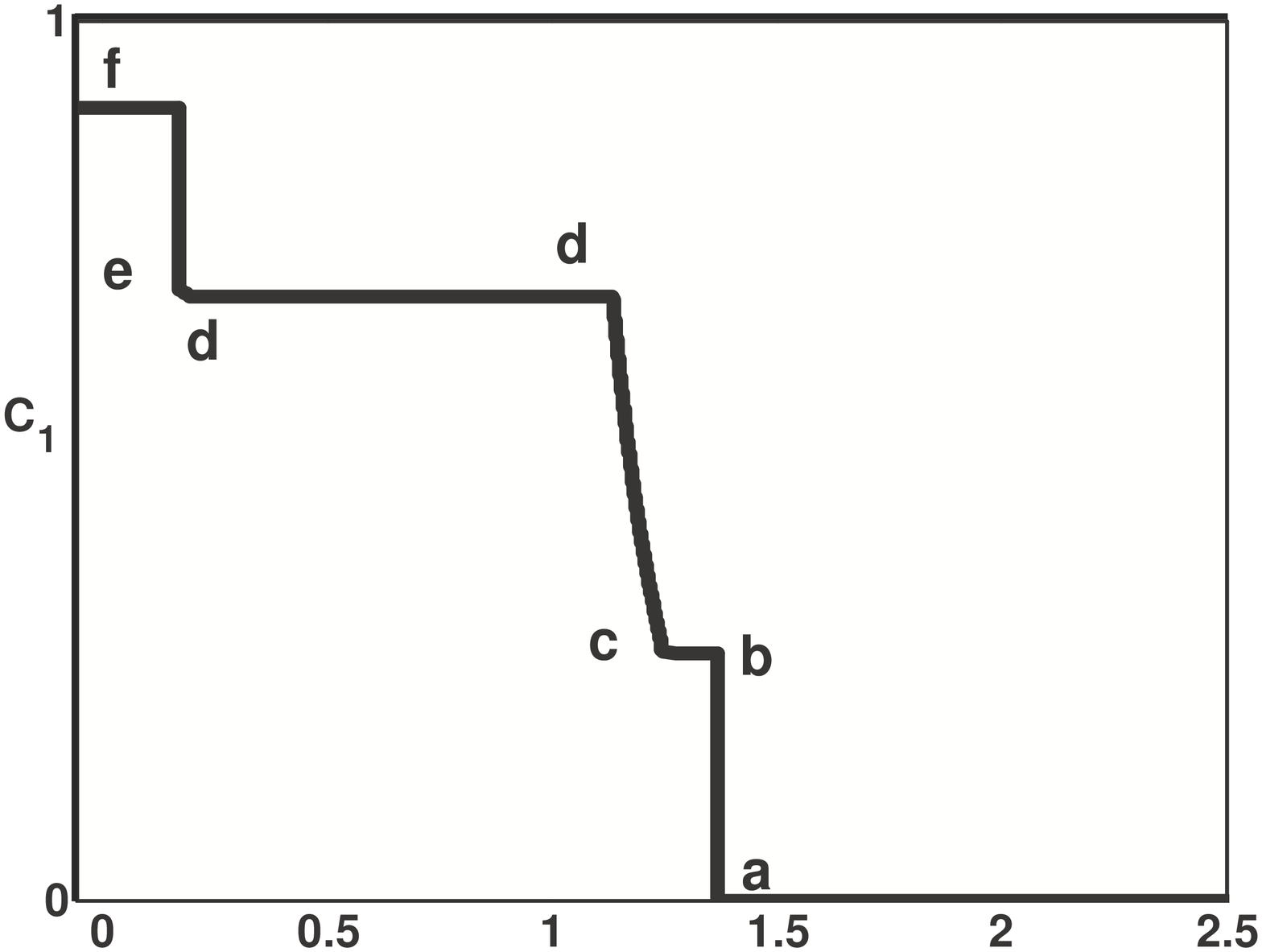}
\end {minipage}
}
\\[-5.5ex]
\subfigure{
\begin{minipage}{2.8in}%
\hspace{-0.17in}
\includegraphics[bb=0mm 0mm 208mm 296mm, angle = 90, scale=0.2]{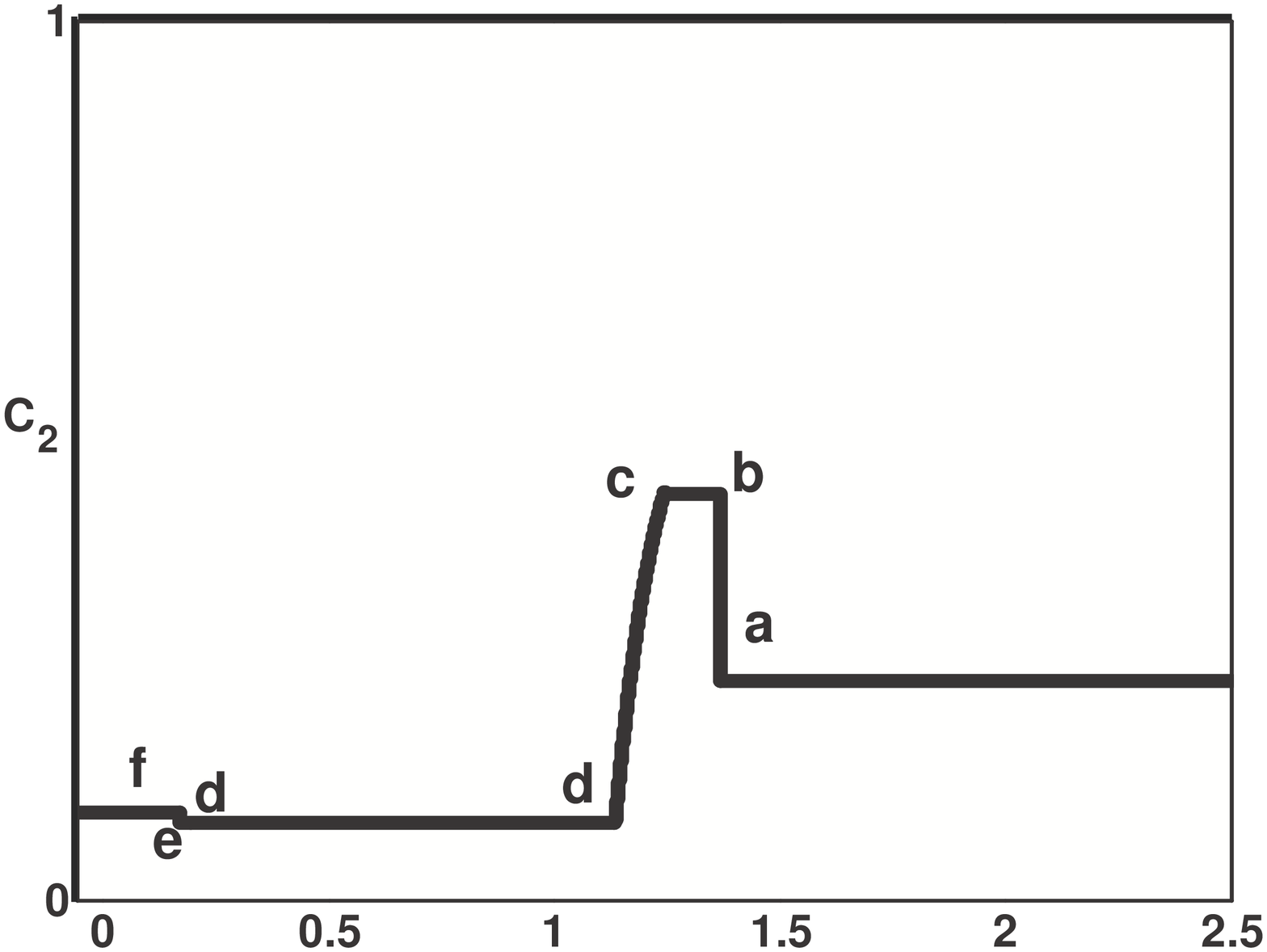}
\end {minipage}%
\begin{minipage}{2.5in}%
\hspace{-0.3in}
\includegraphics[bb=0mm 0mm 208mm 296mm, angle = 90, scale=0.2]{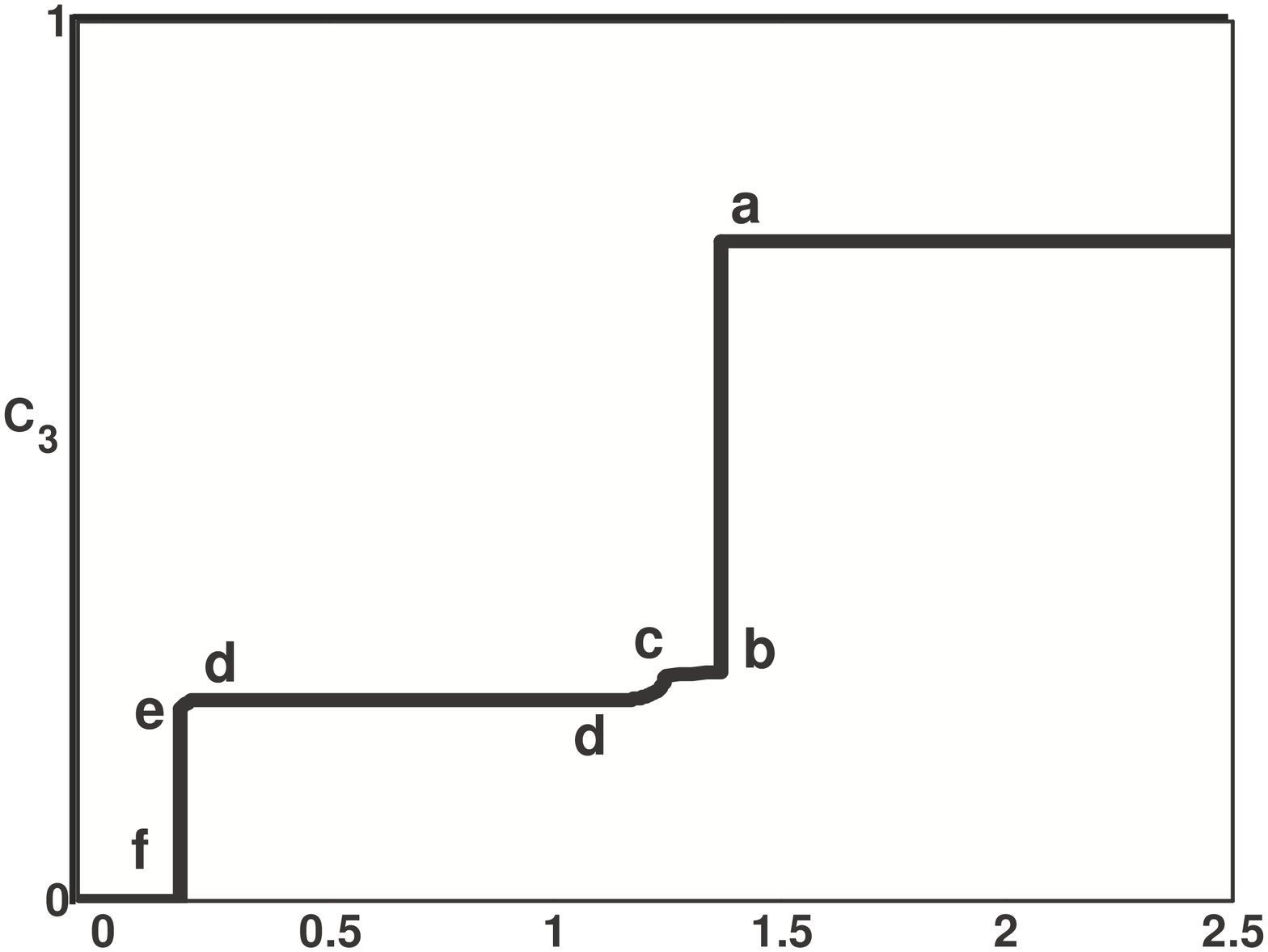}
\end {minipage}%
}
\caption{Solution path in the ternary phase space and composition profiles for 1D ternary example .
Point \textbf{a} is initial oil, \textbf{b} is landing point of leading shock, \textbf{c} is point of weak hyperbolicity,
\textbf{d} is the point where velocity changes from nontie-line eigenvale to tie-line eigenvalue,
\textbf{e} is landing point of trailing shock, and \textbf{f} is injection gas.
}
\label{fig5.5}%
\end{figure}

We will look at the solution profiles for (\ref{eq:eqn5.10}) with the intial conditions
\[\left[C_{1} (x,0){\rm \; \; }C_{2} (x,0)\right]=
\left\{\begin{array}{l} {\left[0.9{\rm \; \; \; \; 0.1}\right],
{\rm \; \; \;  if\; }x<0} \\ {\left[0{\rm .0\; \; 0.25}\right],
{\rm \; \; \;  if\; }x>0} \end{array}\right.
\begin{array}{l} {\left( {\rm amount\; of\; }
C_{1}, {\rm \;  }C_{2} {\rm \; in\; injected\; gas} \right)}
\\ {\left( {\rm amount\; of\; }
C_{1}, {\rm \;  }C_{2} {\rm \; in\; resident\; oil} \right)} \end{array}\] .
Figure (\ref{fig5.5}) illustrates the salient features of this ternary problem. There are two key tie-lines,
one extending through the
initial oil (point \textbf{a}), and another extending through the injection gas (point \textbf{f}).
On each tie-line of the ternary system the solution has a shock and a rarefaction.
As gas is injected, there are two transitions from
the single-phase region to the two-phase region: a leading shock (\textbf{a-b}) on the initial tie-line, and a
trailing shock (\textbf{e-f}) on the injection tie-line. Inside the two-phase region there is a small rarefaction
(\textbf{b-c}), as the composition varies along the initial tie-line. Point (\textbf{c}) is the equal eigenvalue point. Here, the system becomes
weakly hyperbolic. The composition then traces the nontie-line path as a rarefaction (\textbf{c-d}). At point \textbf{d}
the solution encounters the injection tie-line, where the velocity jumps from the nontie-line eigenvalue $\lambda _{nt} $
to the tie-line eigenvalue $\lambda _{t} $. The composition remains constant at \textbf{d} for the entire jump,
forming a zone of constant state. On the injection tie-line there is one more rarefaction (\textbf{d-e})
which connects to the trailing shock (\textbf{e-f}).

\begin{figure}[h!]
\begin{center}
\subfigure{
\includegraphics[bb=0mm 0mm 208mm 296mm, angle = 90, scale=0.27]{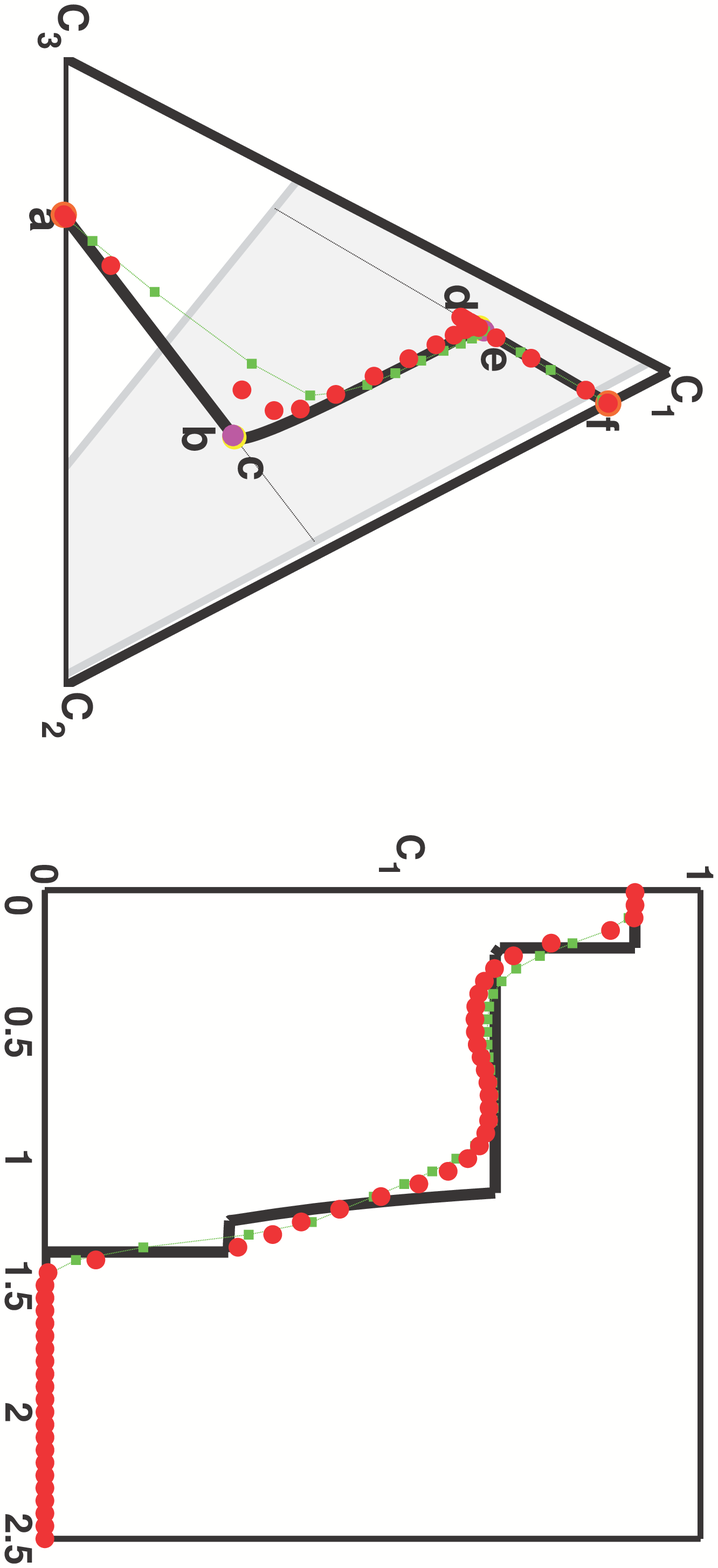}
}
\\
\vspace{-0.9in}
\subfigure{
\includegraphics[bb=0mm 0mm 208mm 296mm, angle = 90, scale=0.27]{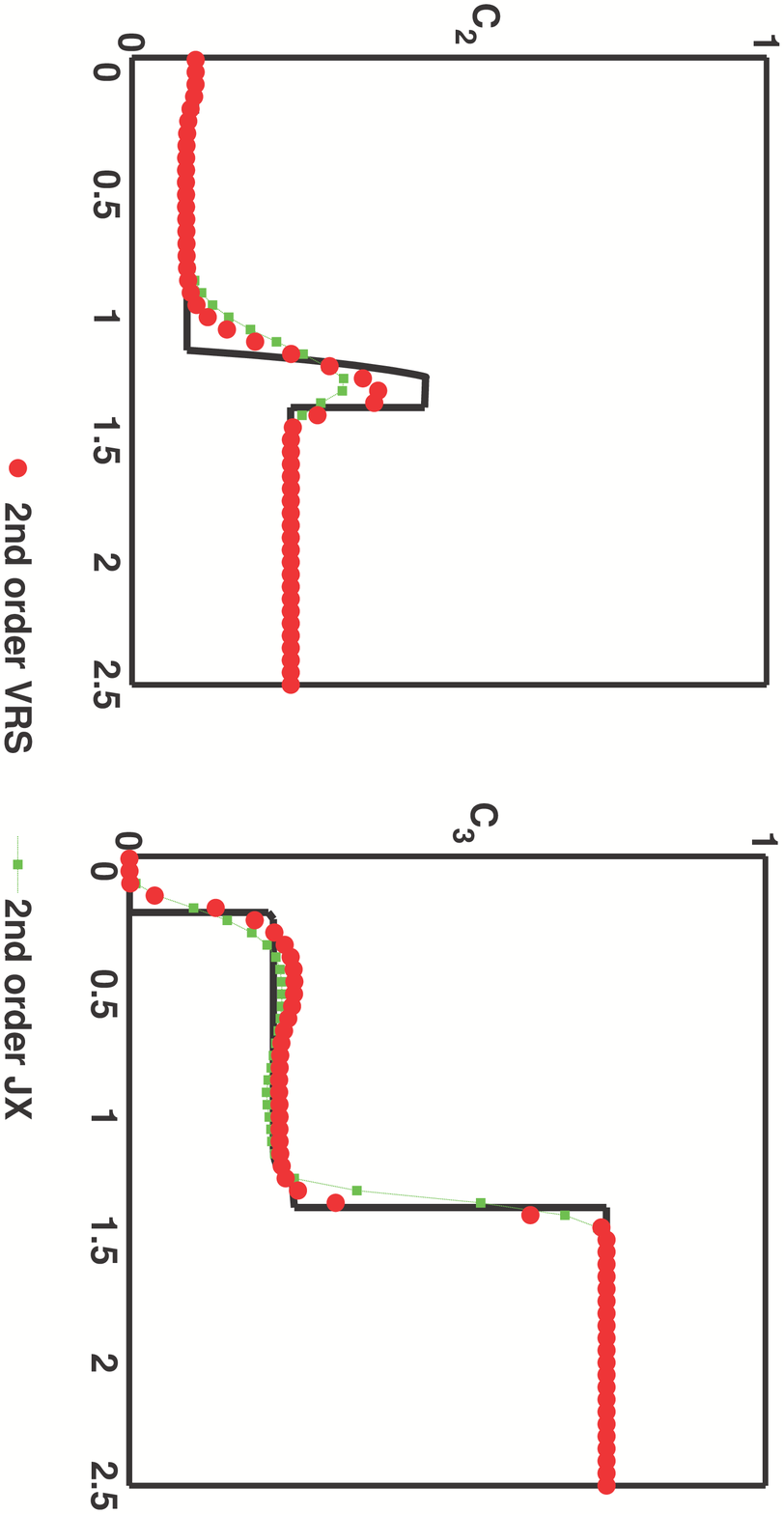}
}
\end{center}
\vspace{-0.6in}
\caption{Solution of ternary displacement with second order JX and VRS schemes,
${T=1}$,  N=50, ${\Delta x = 2.5/N}$, Max speed = 5.4, ${\Delta t = 0.5 \frac{\Delta x}{5.4}}$ .}
\label{fig5.6}
\end{figure}

\begin{figure}[h!]
\begin{center}
\subfigure{
\includegraphics[bb=0mm 0mm 208mm 296mm, angle = 90, scale=0.27]{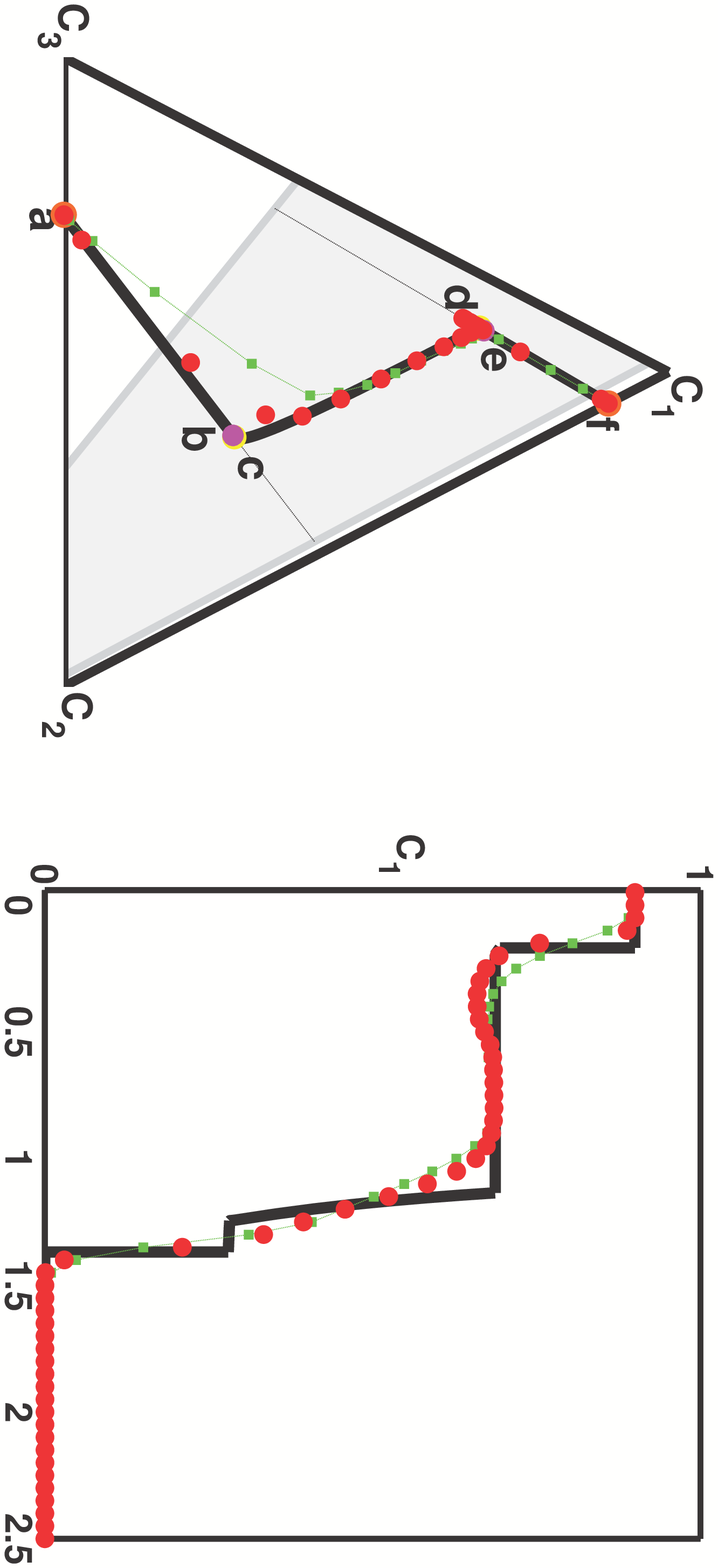}
}
\\
\vspace{-0.9in}
\subfigure{
\includegraphics[bb=0mm 0mm 208mm 296mm, angle = 90, scale=0.27]{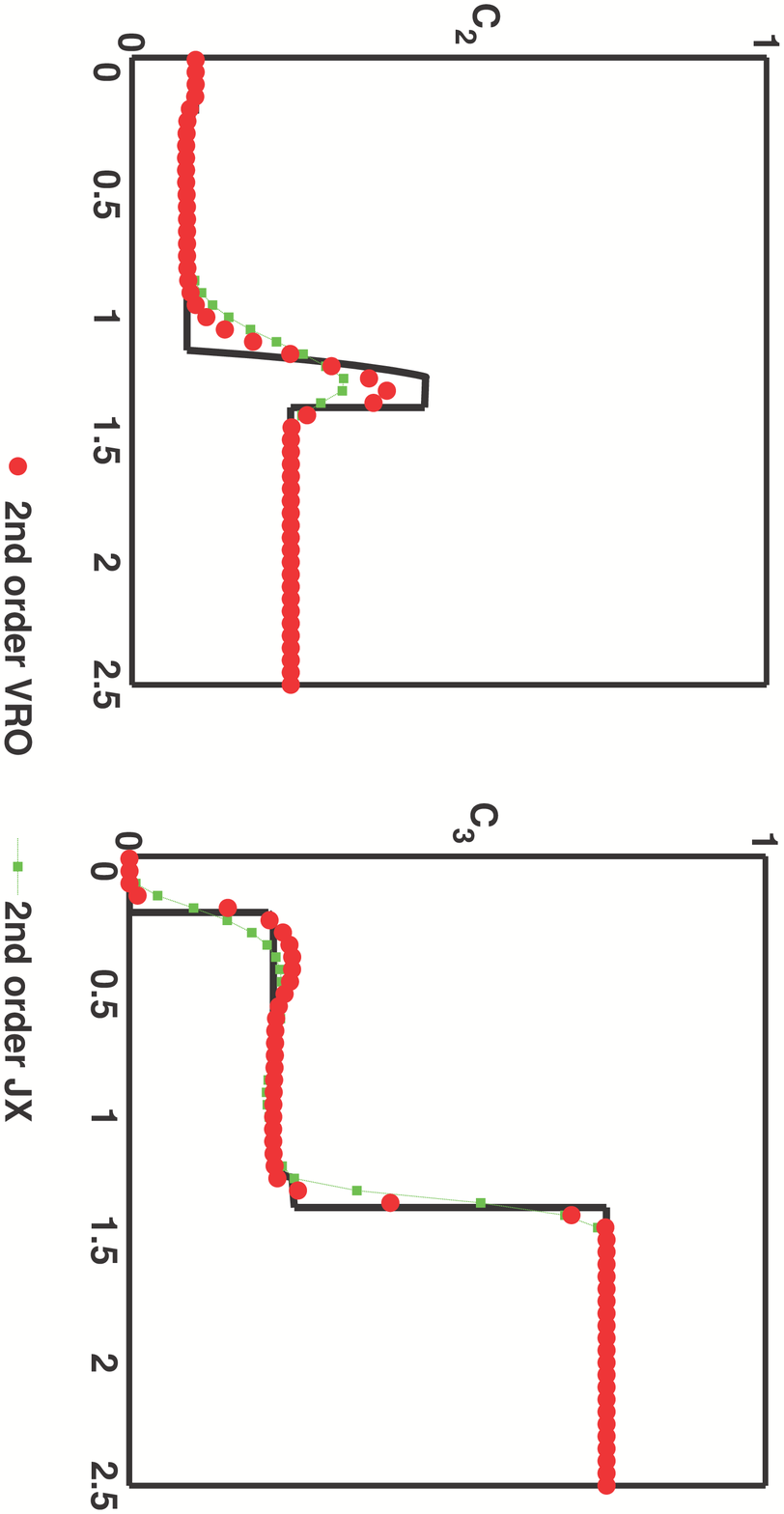}
}
\end{center}
\vspace{-0.6in}
\caption{Solutions of ternary displacement with second order JX and VRO schemes,
${T=1}$,  N=50, ${\Delta x = 2.5/N}$, Max speed = 5.4, ${\Delta t = 0.5 \frac{\Delta x}{5.4}}$ .}
\label{fig5.7}
\end{figure}

The results of second order JX, VRS and VRO schemes are shown in figures (\ref{fig5.6}) and (\ref{fig5.7}). While all
three schemes resolve the $C_1$ and $C_3$ profiles reasonably well, the difference in accuracy can be observed in the
$C_2$ profile, where JX scheme smoothes out the $C_2$ bank. This is also noticeable in the ternary phase diagram, where
the JX path is further from the actual solution path.
\begin{figure}[h!]
\begin{center}
\subfigure{
\includegraphics[bb=0mm 0mm 208mm 296mm, angle = 90, scale=0.27]{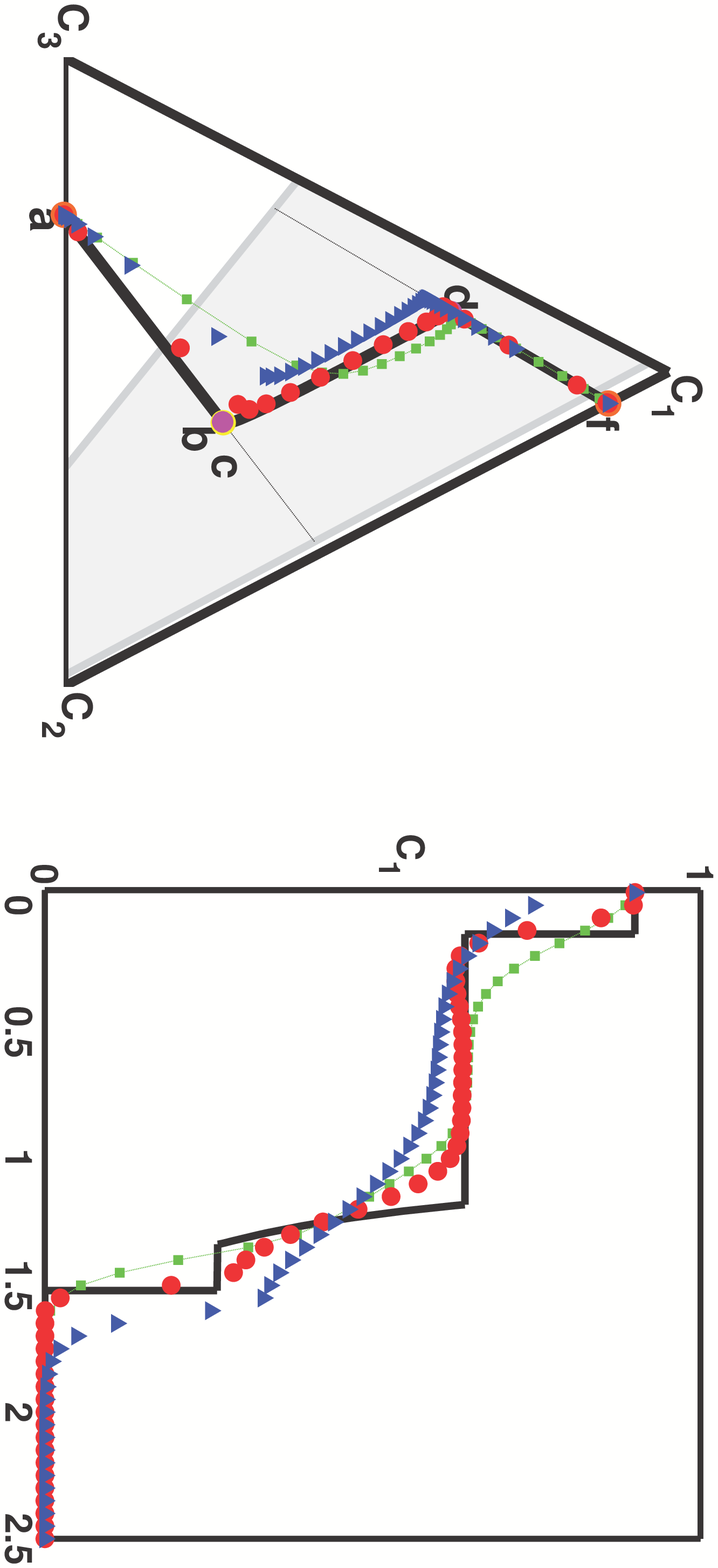}
}
\\
\vspace{-0.9in}
\subfigure{
\includegraphics[bb=0mm 0mm 208mm 296mm, angle = 90, scale=0.27]{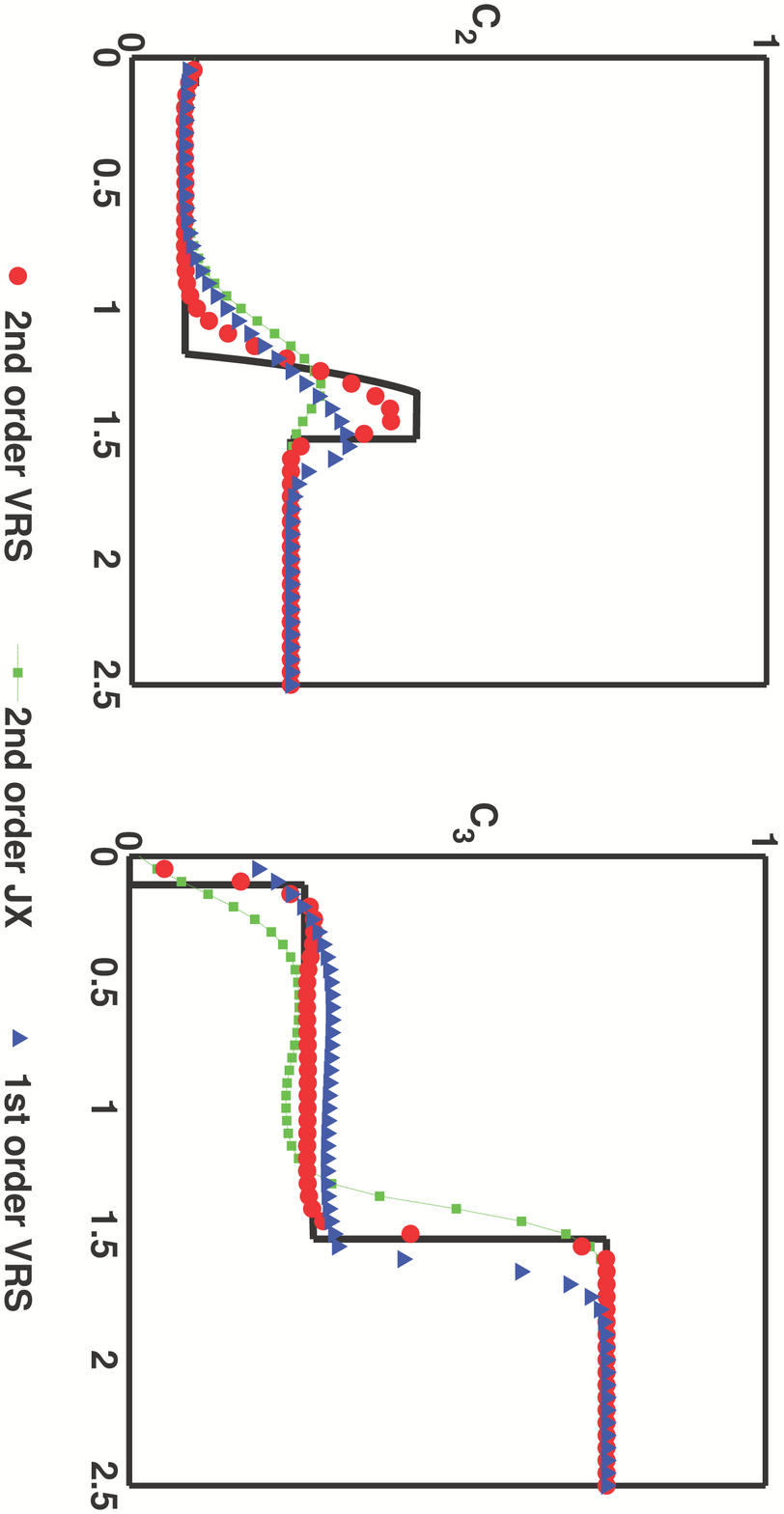}
}
\end{center}
\vspace{-0.6in}
\caption{Solutions of ternary displacement with JX scheme and first and second order VRS schemes,
${T=1}$,  N=50, ${\Delta x = 2.5/N}$, Max speed = 25.5, ${\Delta t = 0.5 \frac{\Delta x}{25.5}}$ .}
\label{fig5.8}
\end{figure}
\begin{figure}[h!]
\begin{center}
\subfigure{
\includegraphics[bb=0mm 0mm 208mm 296mm, angle = 90, scale=0.27]{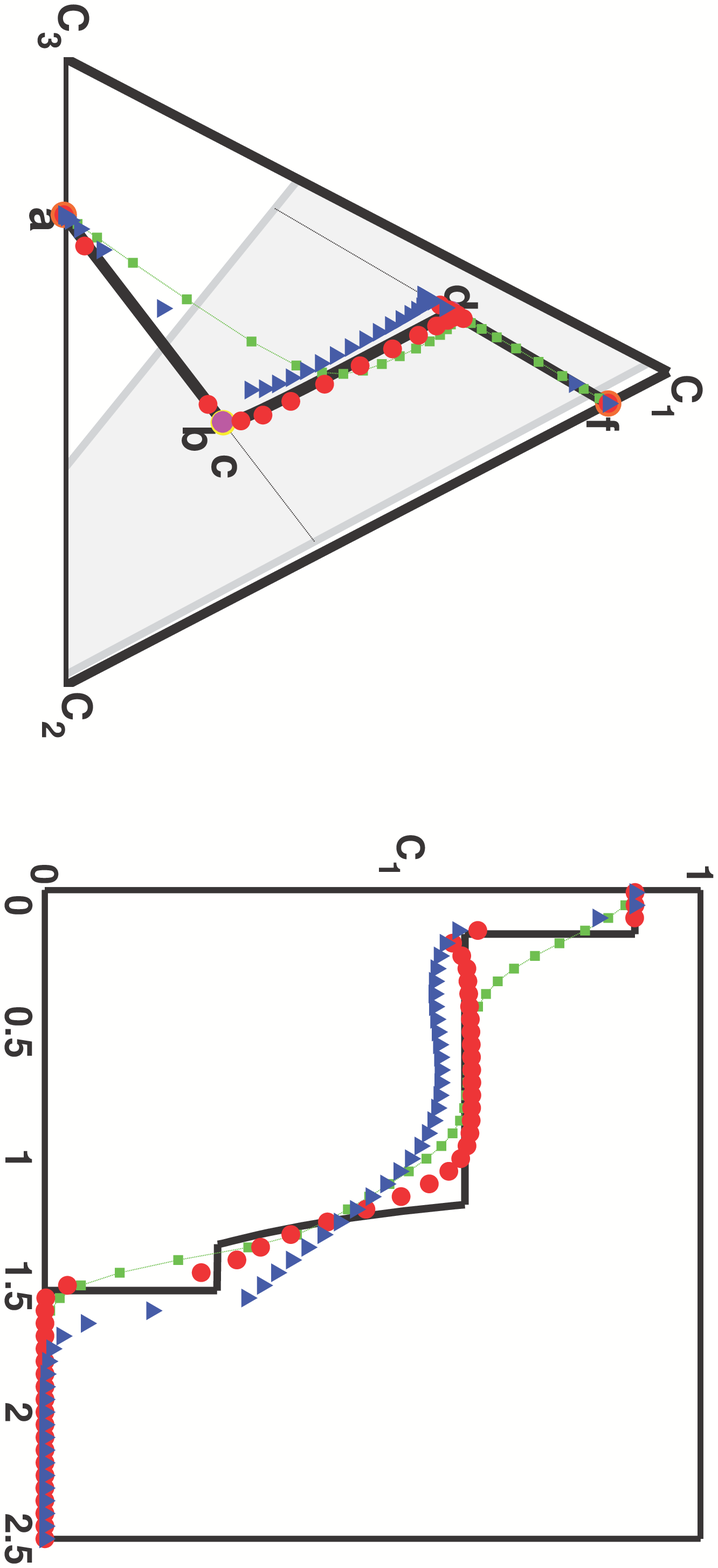}
}
\\
\vspace{-0.9in}
\subfigure{
\includegraphics[bb=0mm 0mm 208mm 296mm, angle = 90, scale=0.27]{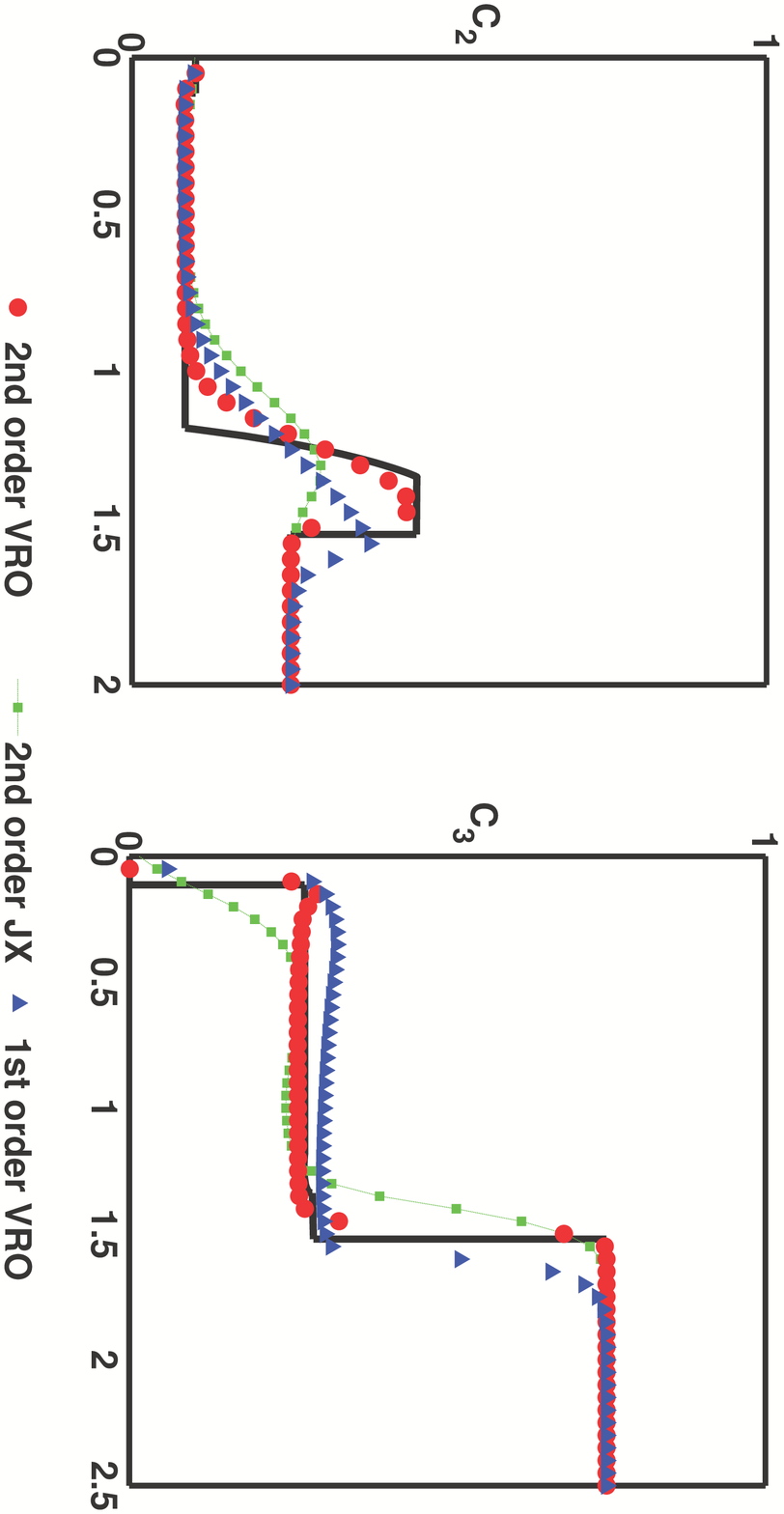}
}
\end{center}
\vspace{-0.6in}
\caption{Solutions of ternary displacement with JX scheme and first and second order VRO schemes,
${T=1}$,  N=50, ${\Delta x = 2.5/N}$, Max speed = 25.5, ${\Delta t = 0.5 \frac{\Delta x}{25.5}}$ .}
\label{fig5.9}
\end{figure}
\begin{figure}[h!]
\begin{center}
\begin{minipage}{2.5in}%
\hspace{-0.7in}
\includegraphics[bb=0mm 0mm 208mm 296mm, angle = 90, scale=0.3]{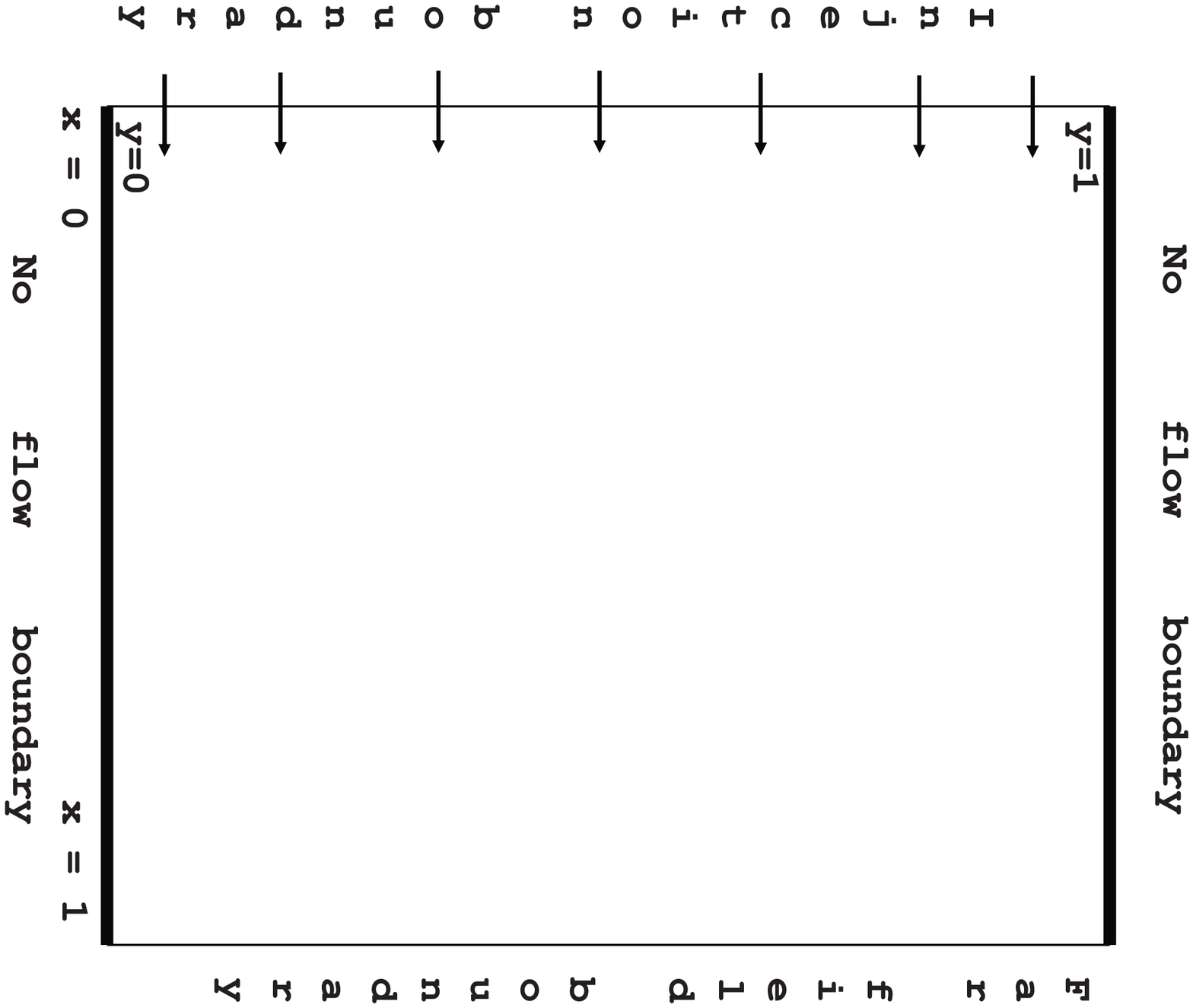}
\caption{$2$D simulation domain}
\label{fig5.10}
\end {minipage}%
\begin{minipage}{2.8in}%
\hspace{-0.59in}
\includegraphics[bb=0mm 0mm 208mm 296mm, angle = 90, scale=0.3]{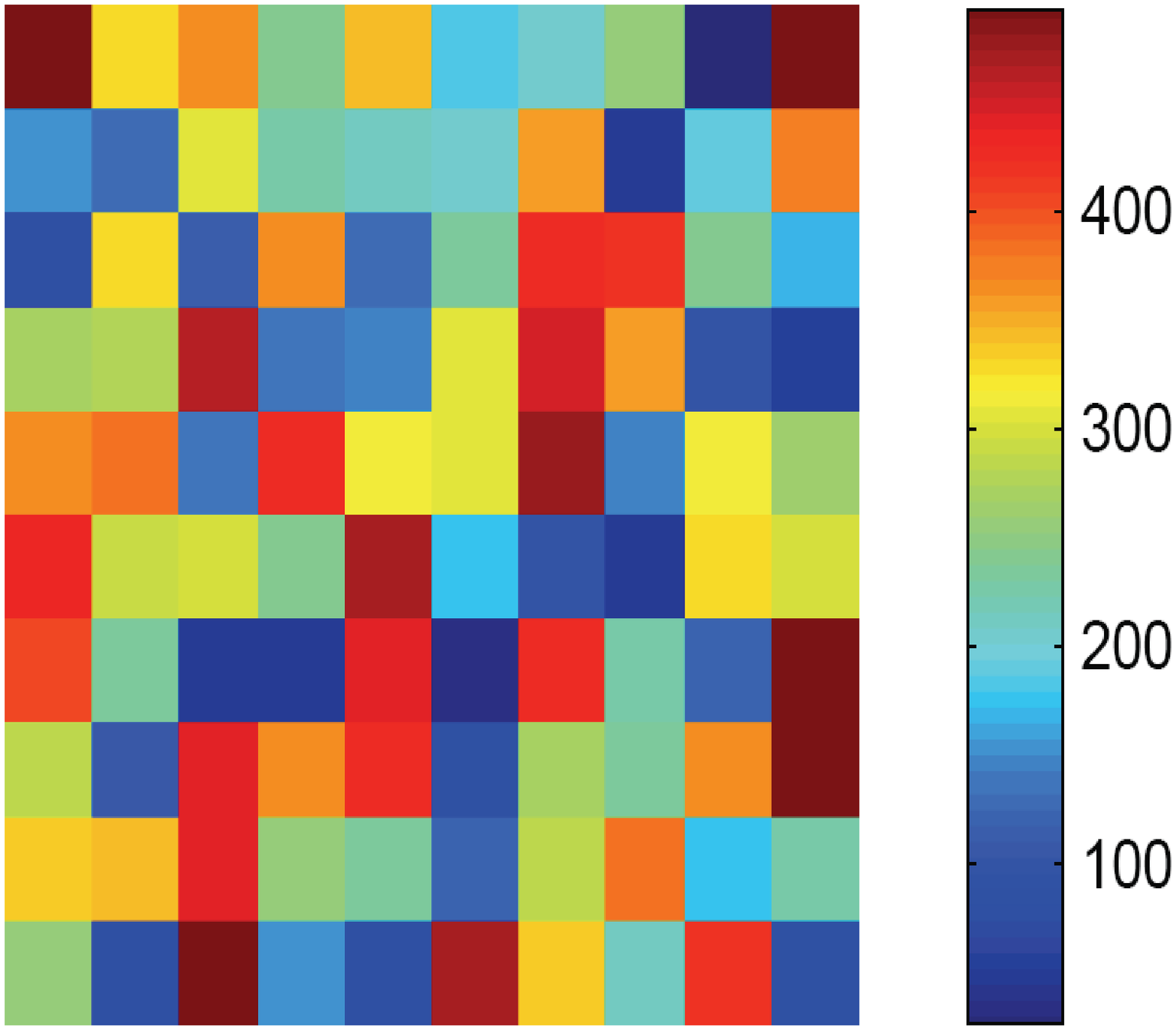}
\caption{Synthetic permeability field for a ${40 \times 40}$ grid }
\label{fig5.11}
\end {minipage}%
\end{center}
\end{figure}

In the above example, the maximum eigenvalue is 5.4 and the average  is 1.
When the fluid properties are changed in this example such that the maximum eigenvalue is 25.5 with the average
 still being 1
(by setting $S_{or} =0.1$, $S_{gc} =0.3$,  ${\frac{\mu _{V}}{\mu _{L}}= \frac{1}{20}} $), the
difference in the accuracy of the VRS/VRO and the JX scheme becomes more apparent (figures (\ref{fig5.8}) and (\ref{fig5.9})).
In fact, in the phase space and also in the $C_2$ profile,  the first order VRS and VRO schemes
show better resolution than the second order JX scheme.
This becomes even more important in multidimensional problems,
where, because of heterogeneous permeability fields or presence of wells (sources and sinks) the global maximum speeds can be
much greater than the average speeds in the domain.

\subsubsection{2D example}
We consider the 2D ternary gas-oil displacement described by,
\begin {subequations}
\label{eqn5.13}
\begin{eqnarray}
\frac{\partial C_{1} }{\partial t}
+\frac{\partial \left(u_{T}^{x} F_{1} \right)}{\partial x}
+\frac{\partial \left(u_{T}^{y} F_{1} \right)}{\partial y} =0,
\nonumber\\
\frac{\partial C_{2} }{\partial t}
+\frac{\partial \left(u_{T}^{x} F_{2} \right)}{\partial x}
+\frac{\partial \left(u_{T}^{y} F_{2} \right)}{\partial y} =0,
\label{eq:eqn5.13a} \\
\nabla . {\bf u}_{T} =
\nabla . \left( \left( {\bf k}\frac{k_{r\textsc{v}}\left(S\right)}{\mu _{\textsc{v}}}
+  {\bf k}\frac{k_{r\textsc{l}}\left(S\right) }{\mu _{\textsc{l}}} \right)
\left(\nabla P \right)
 \right) = 0,
\label{eq:eqn5.13b}
\end{eqnarray}

where  ${\bf u}_{T} =\left[u_{T}^{x} {\rm \; \; }u_{T}^{y} \right]^{T} $ is the vector of total Darcy velocities
given by
\begin{eqnarray}
{\bf u}_{T} =
\left( \left( {\bf k}\frac{k_{r\textsc{v}}\left(S\right)}{\mu _{\textsc{v}}}
+  {\bf k}\frac{k_{r\textsc{l}}\left(S\right) }{\mu _{\textsc{l}}} \right)
\left(\nabla P \right)
 \right),
\label{eq:eqn5.13c}
\end{eqnarray}
\end {subequations}
and $ \bf{k} = \left[\begin{array}{cc} k^x & 0 \\0  & k^y \end{array} \right]$  is the domain permeability tensor.
The fluxes $F_1$ and $F_2$ in the above system are the same as in 1D:
 $F_{1,2} =c_{1,2 \textsc{v}} f + c_{1,2 \textsc{l}} \left(1-f\right) $. Hence, the 2D system also exhibits strong nonlinear coupling and the solutions are composed of compound waves.
The property of weak hyperbolicity also gets carried over to 2D since the linear combination of the Jacobians can
still have an incomplete set of eigenvectors.

We solve the coupled system (\ref{eqn5.13}) sequentially. In each time step, (\ref{eq:eqn5.13b}) is solved for pressure
using the component fractions from the previous time step. A finite difference discretization \cite{as} is applied on
(\ref{eq:eqn5.13b}) on the 2-D Cartesian grid and a direct sparse solver is used to solve the resulting algebraic system of
equations.
The velocities are then calculated from the newly computed pressure field as given in equation (\ref{eq:eqn5.13c})
and used in the transport system (\ref{eq:eqn5.13a}) to advance the component fractions (for details see \cite{as}, chapters 3 and 5).
\begin{figure}[h!]
\begin{center}
\vspace {-0.17in}
\subfigure{
\hspace {-0.0in}
\includegraphics[bb=0mm 0mm 208mm 296mm, angle = 90, scale=0.2]{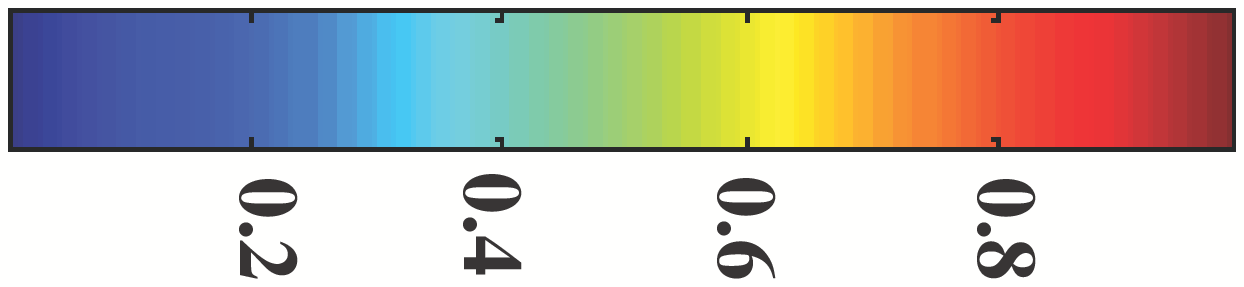}
}
\subfigure{
\hspace{-2in}
\includegraphics[bb=0mm 0mm 208mm 296mm, angle = 90,scale=0.2]{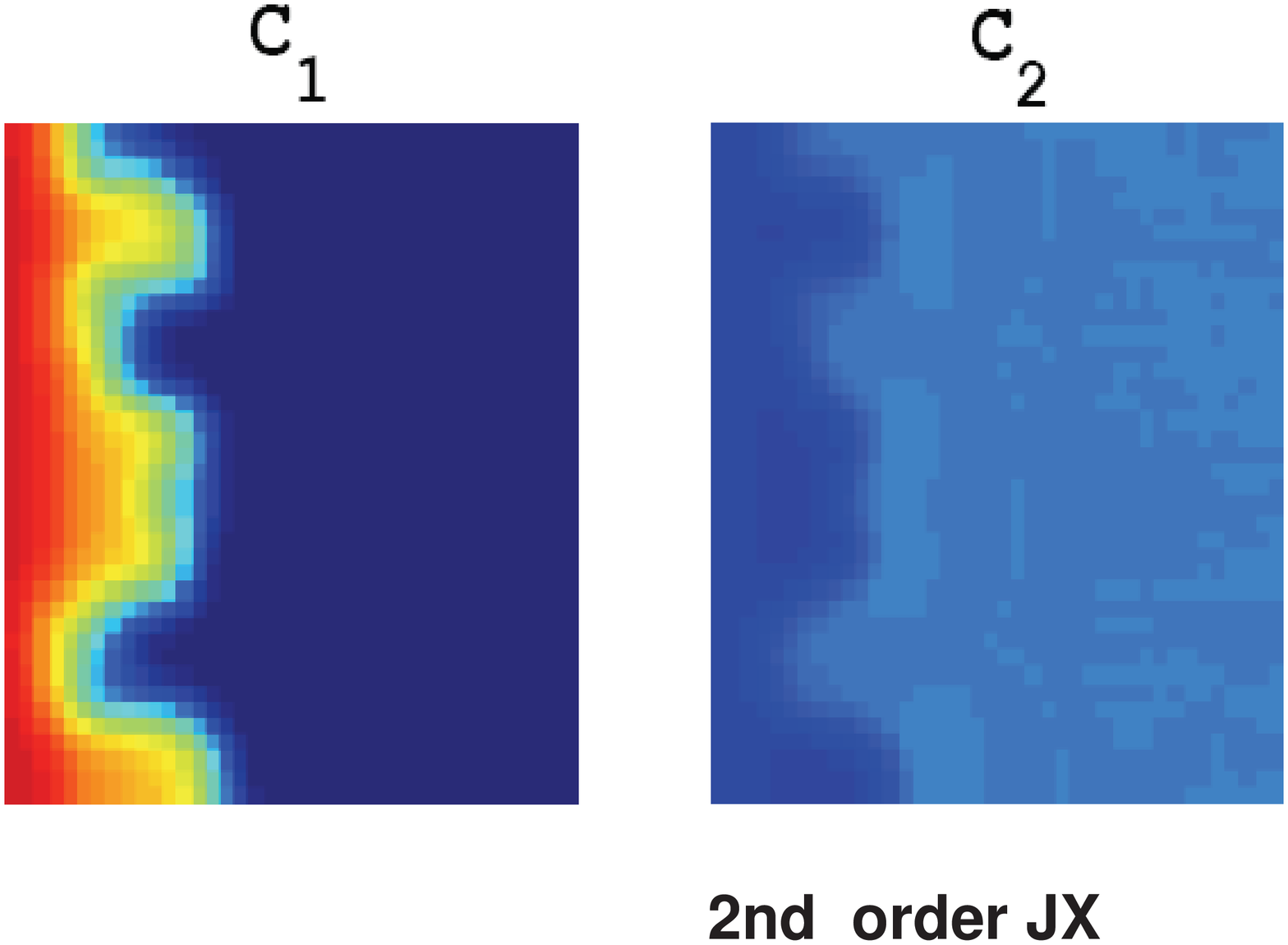}
}
\subfigure{
\hspace{-0.5in}
\includegraphics[bb=0mm 0mm 208mm 296mm, angle = 90, scale=0.2]{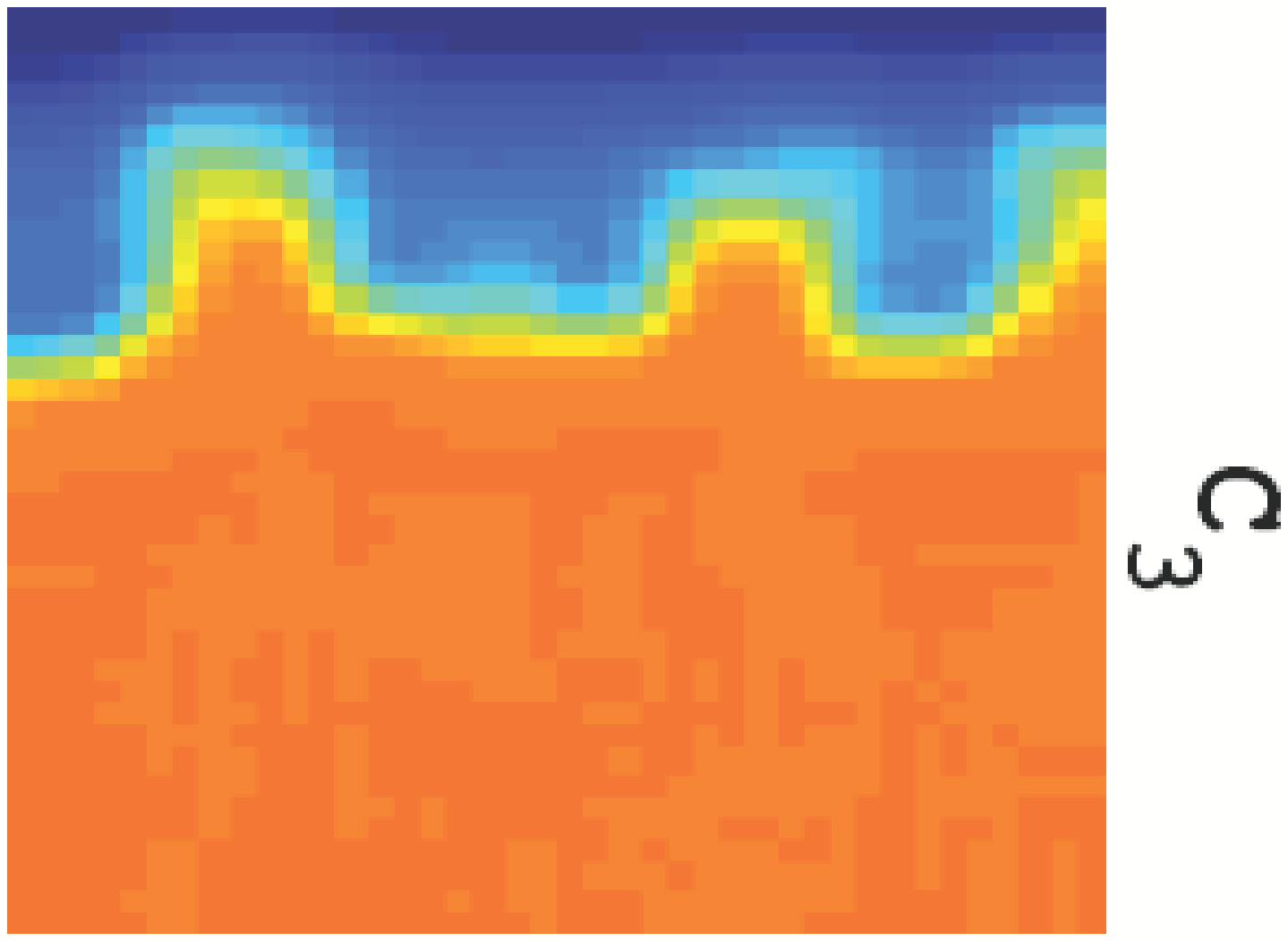}
}
\end{center}
\vspace{-0.3in}
\caption{ Solution profiles for the $2$D ternary displacement at a time ${T=0.2}$
with the second order JX scheme on a ${40 \times 40}$ grid,
JX subcharacteristic speed in either dimensions = 5.8, ${\Delta t = 0.5 \frac{\Delta x}{5.8}}$
.}
\label{fig5.12}
\end{figure}
\begin{figure}[h!]
\begin{center}
\vspace {-0.2in}
\subfigure{
\hspace {-0.0in}
\includegraphics[bb=0mm 0mm 208mm 296mm, angle = 90, scale=0.2]{figures/colorbar.eps}
}
\subfigure{
\hspace{-2.0in}
\includegraphics[bb=0mm 0mm 208mm 296mm, angle = 90,scale=0.2]{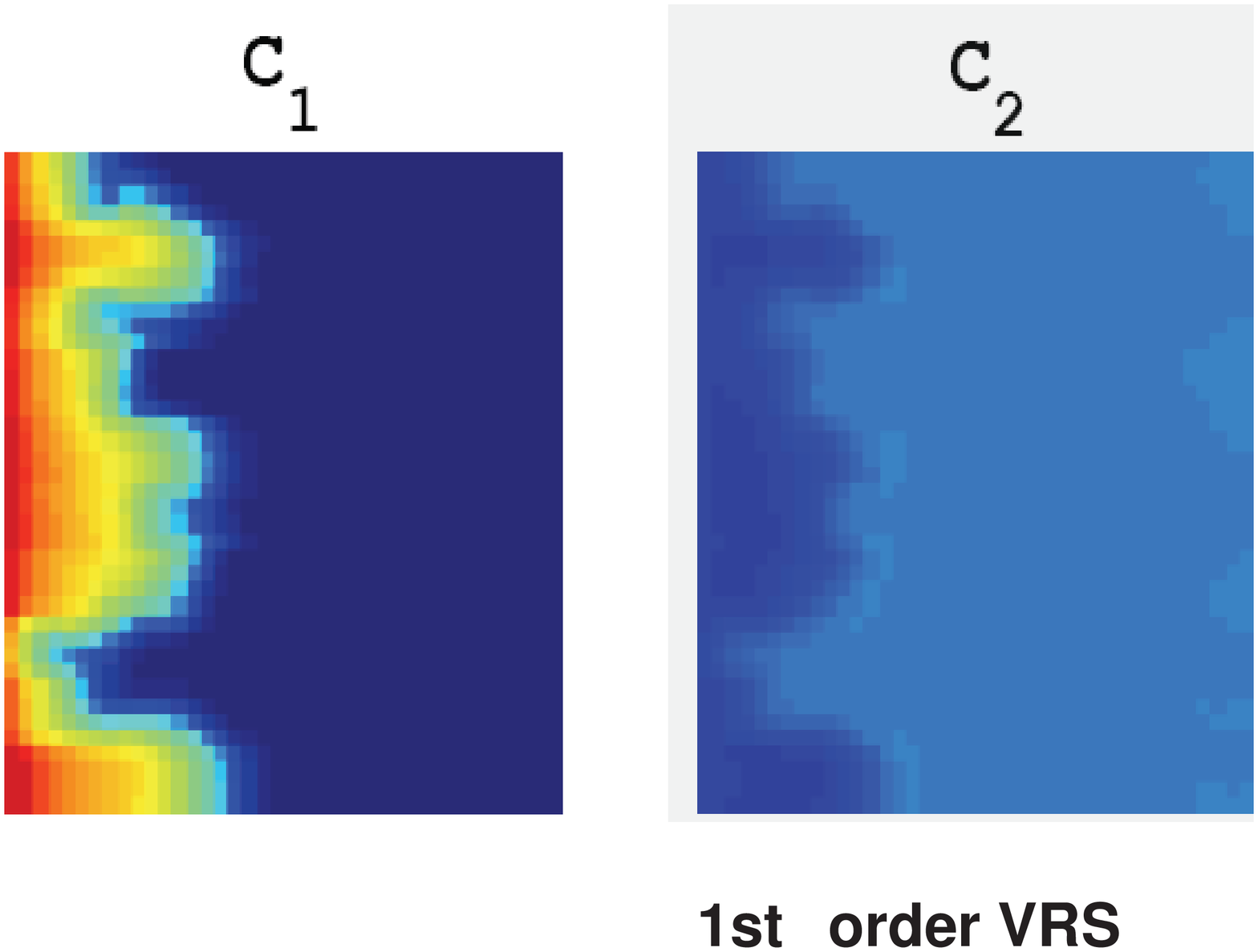}
}
\subfigure{
\hspace{-0.5in}
\includegraphics[bb=0mm 0mm 208mm 296mm, angle = 90, scale=0.2]{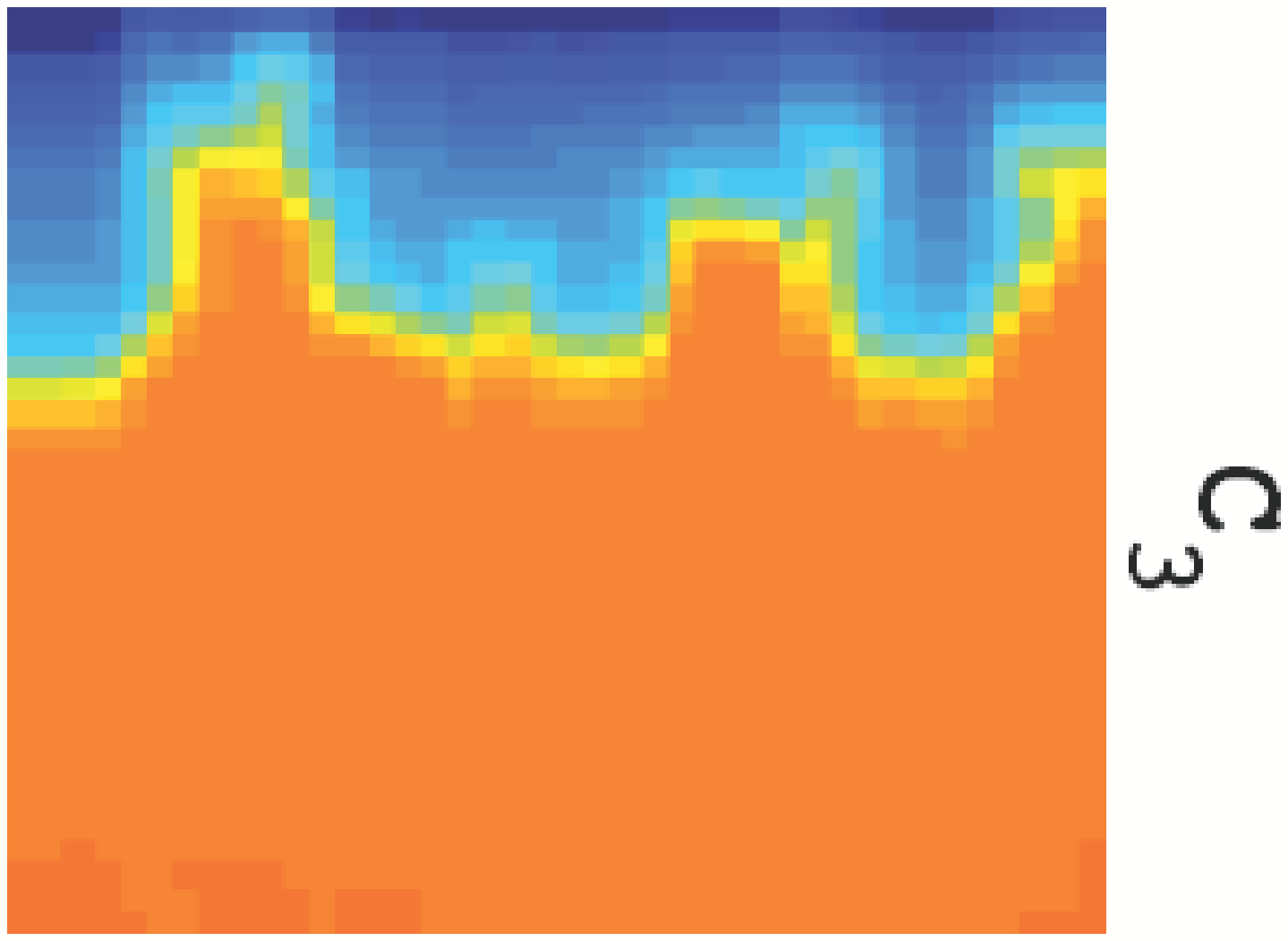}
}
\\
\vspace {-0.4in}
\subfigure{
\hspace {-0.0in}
\includegraphics[bb=0mm 0mm 208mm 296mm, angle = 90, scale=0.2]{figures/colorbar.eps}
}
\subfigure{
\hspace{-2.0in}
\includegraphics[bb=0mm 0mm 208mm 296mm, angle = 90,scale=0.2]{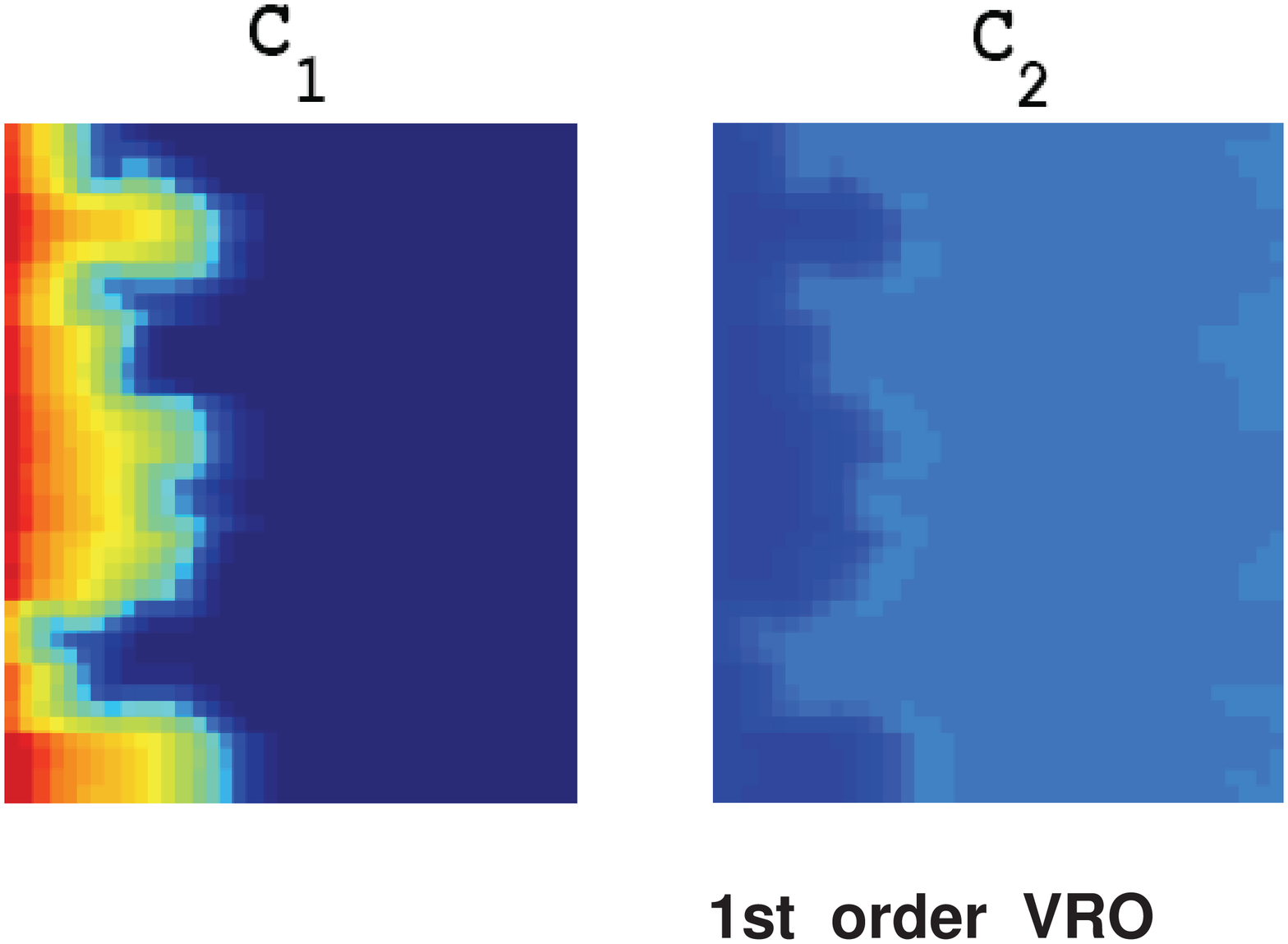}
}
\subfigure{
\hspace{-0.5in}
\includegraphics[bb=0mm 0mm 208mm 296mm, angle = 90, scale=0.2]{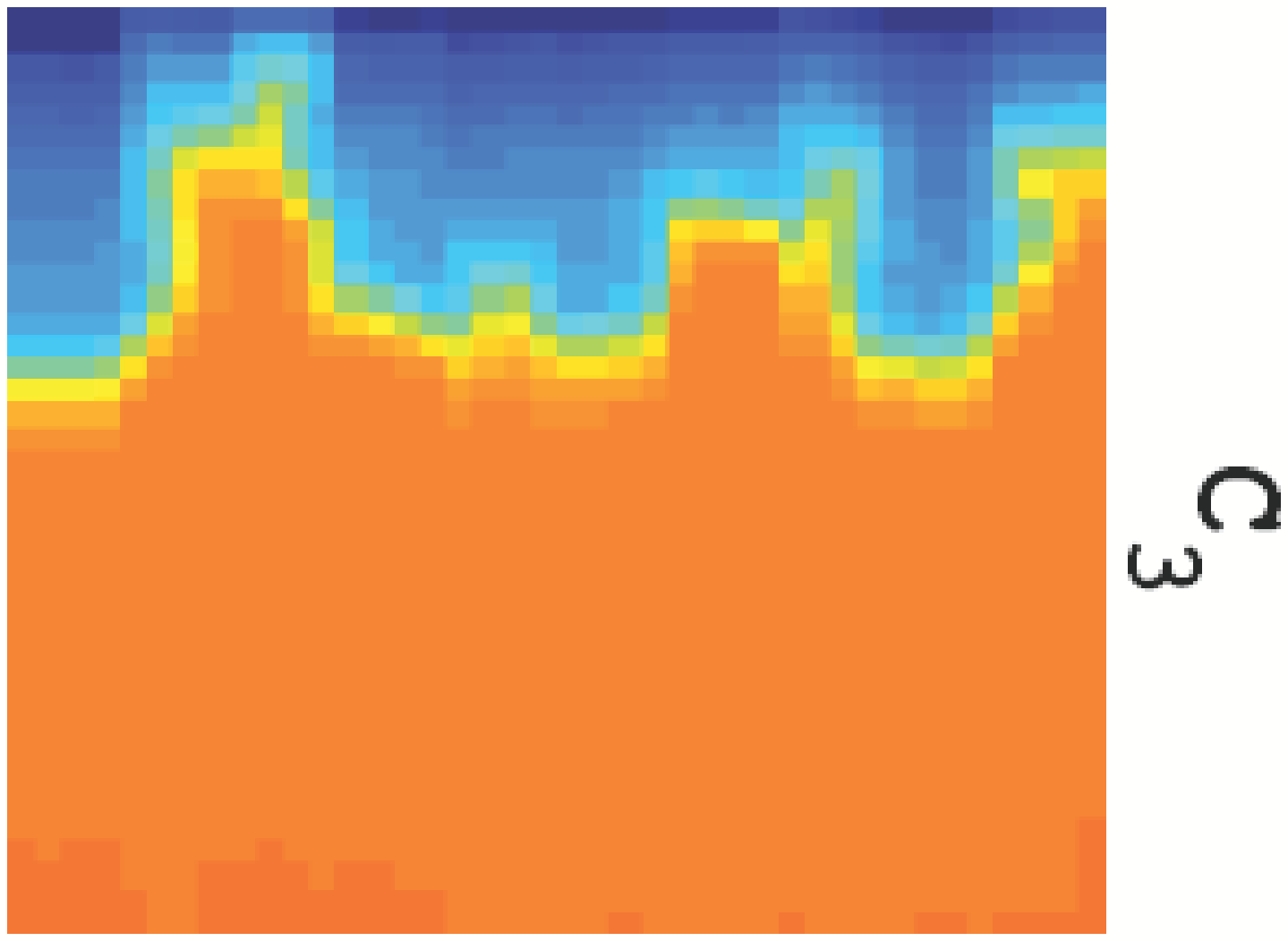}
}
\end{center}
\vspace{-0.3in}
\caption{ Solution profiles for $2$D ternary displacement at ${T=0.2}$
with first order optimal variable relaxed scheme on a ${40 \times 40}$ grid,
max subcharacteristic speed = 5.8, ${\Delta t = 0.5 \frac{\Delta x}{5.8}}$
.}
\label{fig5.14}
\end{figure}

With the initial and injection conditions given by
\[\left[C_{1} (x,0){\rm \; \; }C_{2} (x,0)\right]=
\left\{\begin{array}{l} {\left[0.9{\rm \; \; \; \; 0.1}\right],
{\rm \; \; \;  if\; }x<0} \\ {\left[0{\rm .0\; \; 0.25}\right],
{\rm \; \; \;  if\; }x>0} \end{array}\right.
\begin{array}{l} {\left( {\rm amount\; of\; }
C_{1}, {\rm \;  }C_{2} {\rm \; in\; injected\; gas} \right)}
\\ {\left( {\rm amount\; of\; }
C_{1}, {\rm \;  }C_{2} {\rm \; in\; resident\; oil} \right),} \end{array}\]
the system (\ref{eqn5.13}) is solved on the domain shown in figure (\ref{fig5.10}), with permeability
$ k^x$ and $ k^y$ both given by the synthetic heterogeneous field of figure (\ref{fig5.11}). The fluid properties are taken to be
 $ K_1 = 2.5, K_2 = 1.5, K_3 = 0.05, S_{or} =0.1$, $S_{gc} =0.05$ and
 ${\frac{\mu _{\textsc{v}}}{\mu _{\textsc{l}}}= \frac{1}{2}} $.
 Gas is injected at the left boundary  at a constant nondimensional rate 1.
 For this problem,   the velocities of the
 components in $x$ and $y$ directions are given by the eigenvalues  (\ref{eqn5.11}) multiplied by  $u_{T}^{x}$
 and $u_{T}^{y}$ respectively. The  $x$-velocities are
 always nonnegative, but the  $y$-velocities can have either sign because of the heterogeneity and the boundary conditions.

\begin{figure}[h!]
\begin{center}
\vspace {-0.3in}
\subfigure{
\hspace {-0.0in}
\includegraphics[bb=0mm 0mm 208mm 296mm, angle = 90, scale=0.2]{figures/colorbar.eps}
}
\subfigure{
\hspace{-2.0in}
\includegraphics[bb=0mm 0mm 208mm 296mm, angle = 90,scale=0.2]{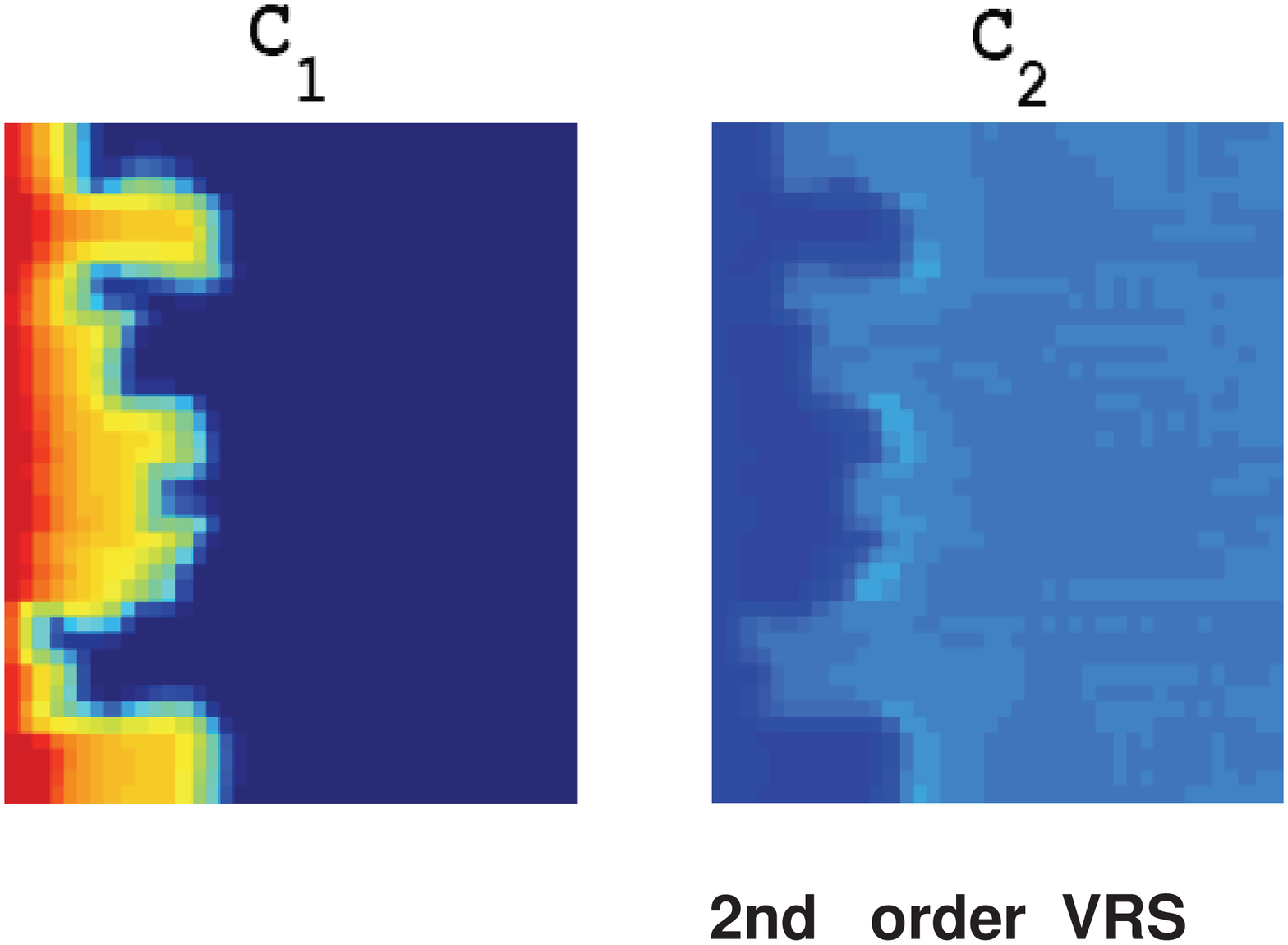}
}
\subfigure{
\hspace{-0.5in}
\includegraphics[bb=0mm 0mm 208mm 296mm, angle = 90, scale=0.2]{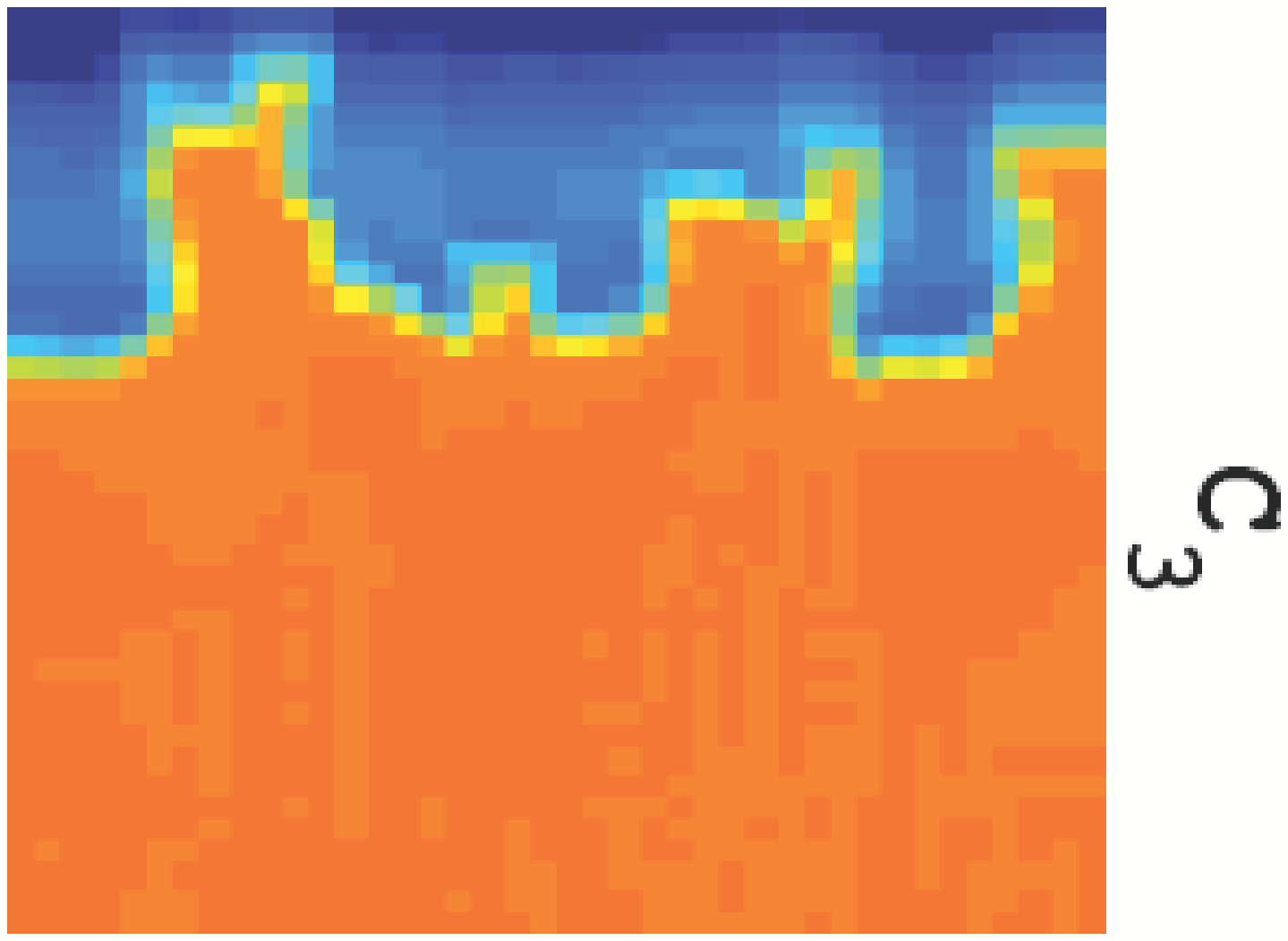}
}
\\
\vspace {-0.4in}
\subfigure{
\hspace {-0.0in}
\includegraphics[bb=0mm 0mm 208mm 296mm, angle = 90, scale=0.2]{figures/colorbar.eps}
}
\subfigure{
\hspace{-2.0in}
\includegraphics[bb=0mm 0mm 208mm 296mm, angle = 90,scale=0.2]{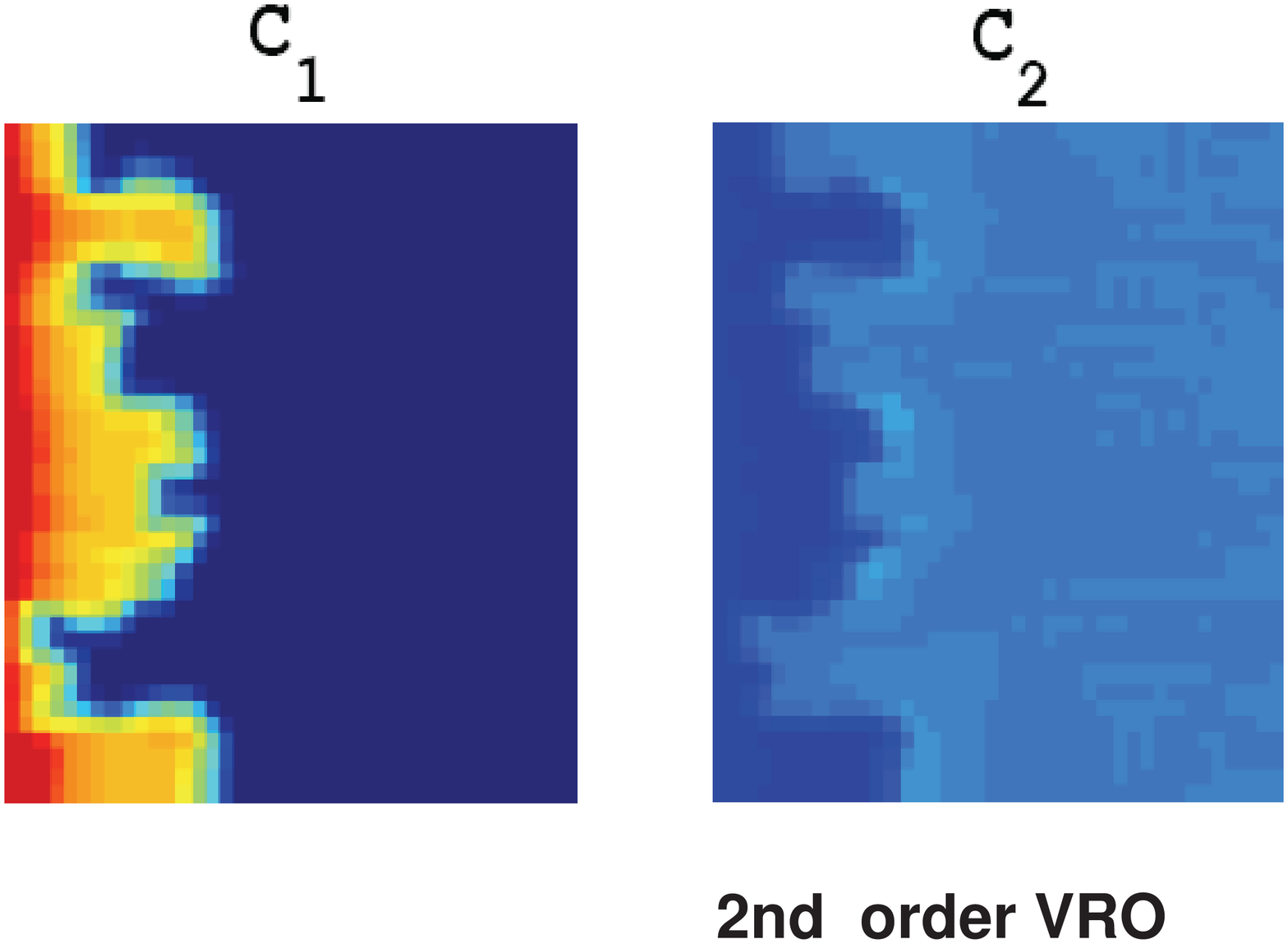}
}
\subfigure{
\hspace{-0.5in}
\includegraphics[bb=0mm 0mm 208mm 296mm, angle = 90, scale=0.2]{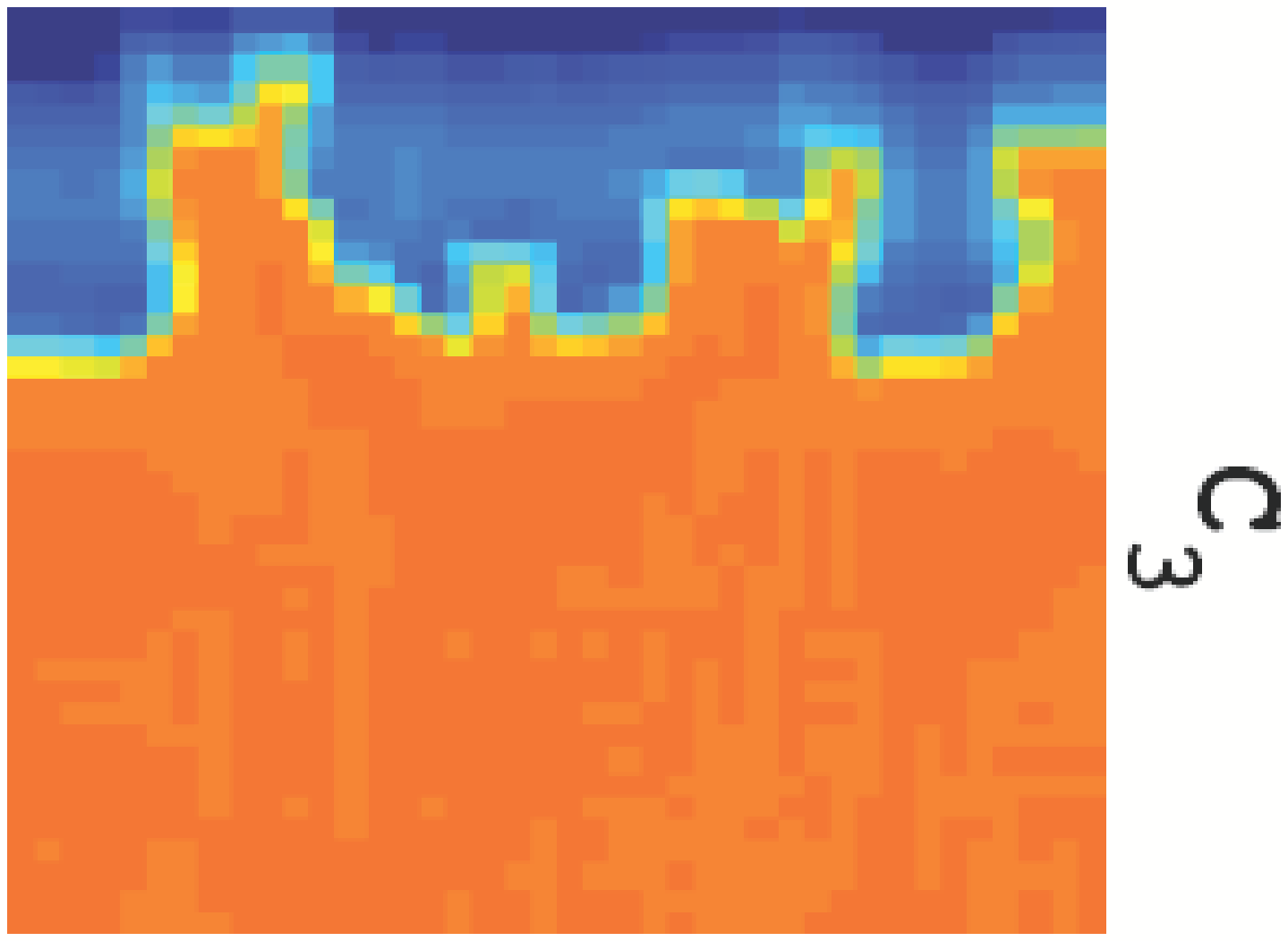}
}
\end{center}
\vspace{-0.3in}
\caption{ Solution profiles for $2$D ternary displacement at ${T=0.2}$
with second order VRO on a ${40 \times 40}$ grid,
max subcharacteristic speed = 5.8, ${\Delta t = 0.5 \frac{\Delta x}{5.8}}$
.}
\label{fig5.16}
\end{figure}
\begin {table}[h]
\caption{Choice of subcharacteristic speeds for the 2D ternary displacement example}
\begin {center}
\begin{tabular}{|p{1.5in}|p{1.6in}|p{1.6in}|} \hline
 & $a^{x} $ & $a^{y} $ \\ \hline
Minimum $a^{y} $  & 9.4  (4 times average) & 3.4  (10 times average) \\ \hline
Minimum $a^{x} $ & 5.3 (2.3 times average) & 8.5  (24 times average) \\ \hline
$a=\sqrt{\left(\lambda _{\max }^{x} \right)^{2} +\left(\lambda _{\max }^{y} \right)^{2} } $
& 5.8 (2.4 times average) & 5.8  (17 times average) \\ \hline
\end{tabular}
\end {center}
\label{table5.2}
\end {table}
\begin{figure}[h!]
\begin{center}
\vspace {-0.3in}
\subfigure{
\hspace {-0.0in}
\includegraphics[bb=0mm 0mm 208mm 296mm, angle = 90, scale=0.2]{figures/colorbar.eps}
}
\subfigure{
\hspace{-2.0in}
\includegraphics[bb=0mm 0mm 208mm 296mm, angle = 90,scale=0.2]{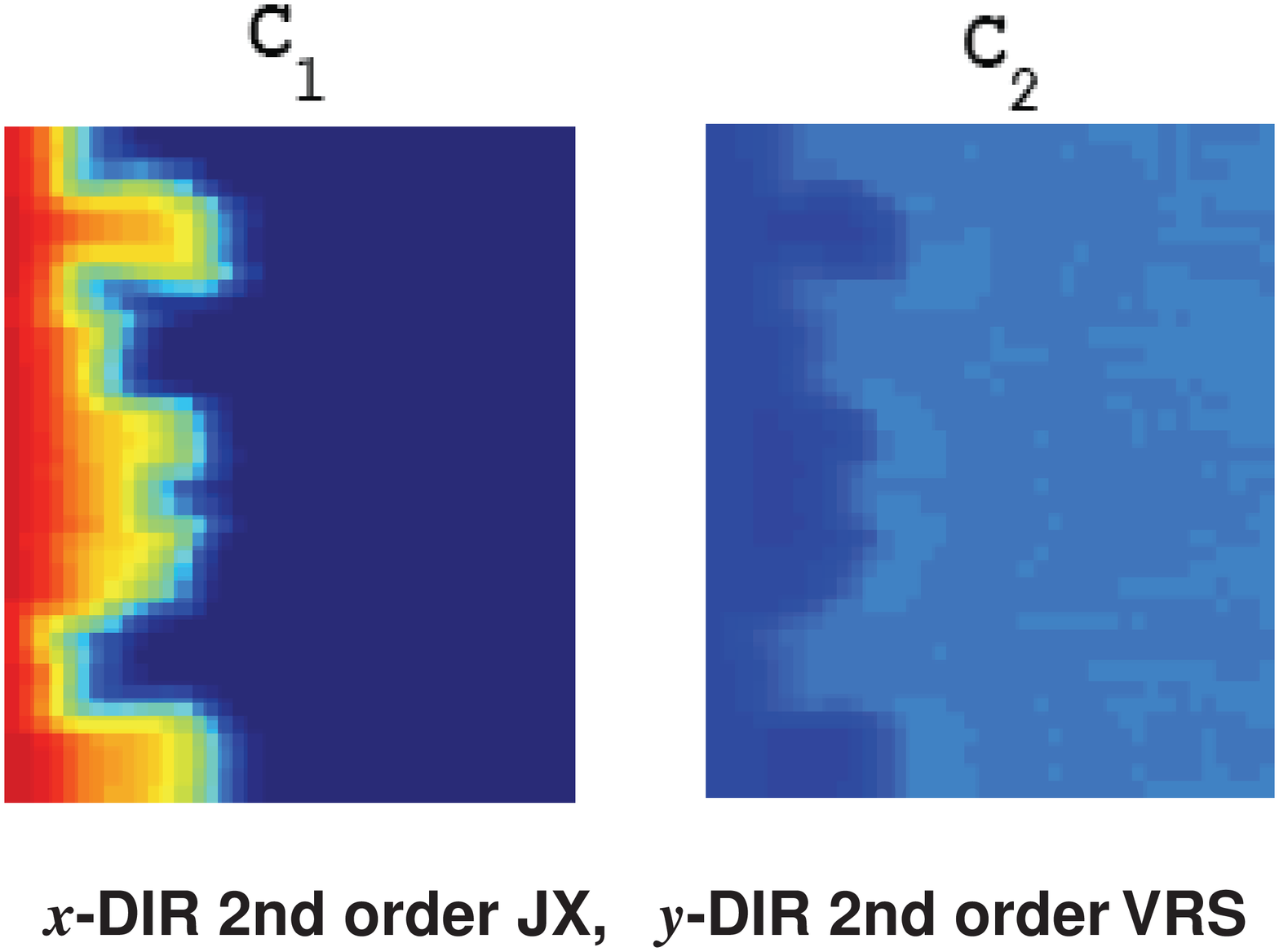}
}
\subfigure{
\hspace{-0.5in}
\includegraphics[bb=0mm 0mm 208mm 296mm, angle = 90, scale=0.2]{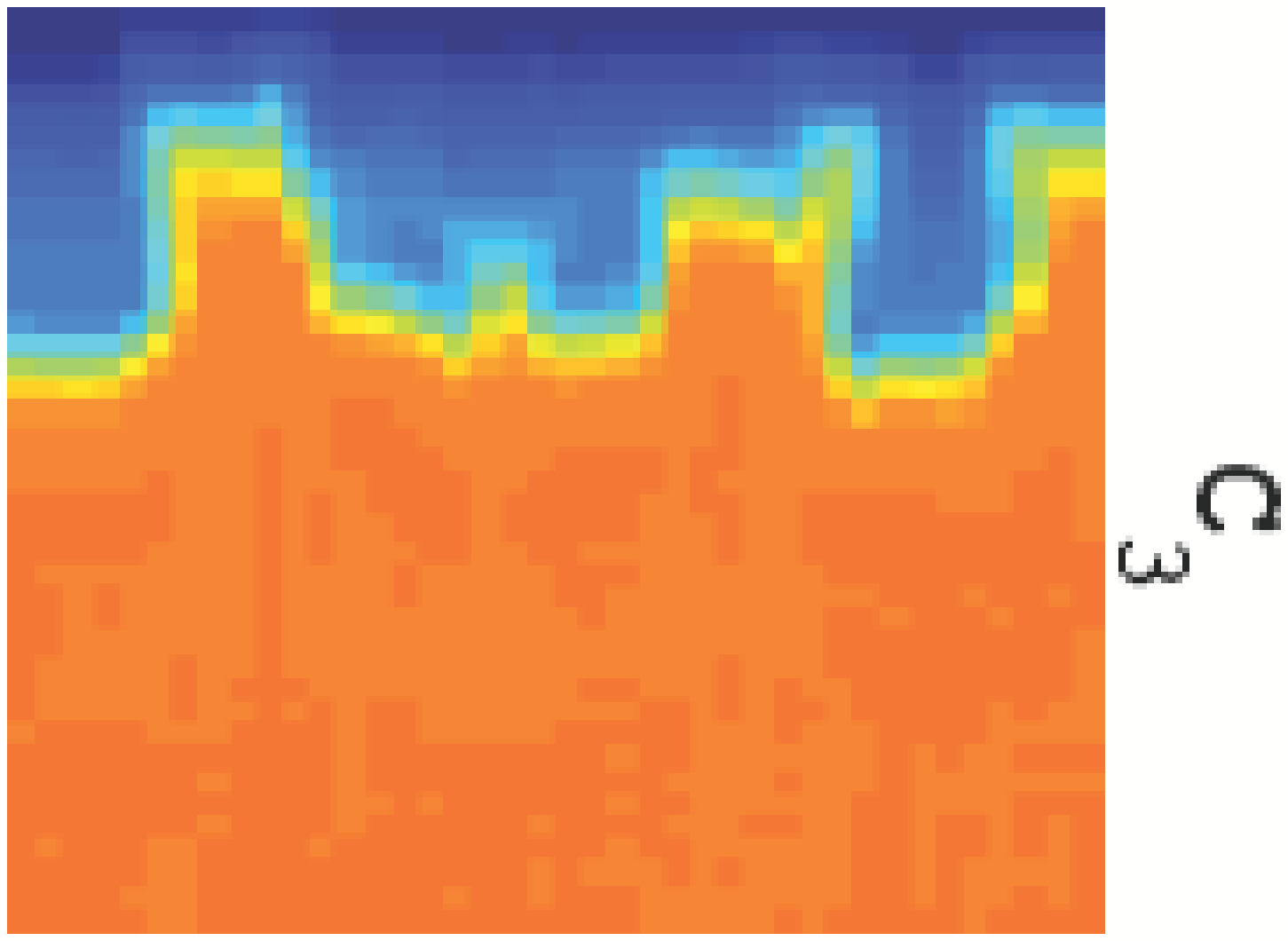}
}
\end{center}
\vspace{-0.3in}
\caption{ Solution profiles for $2$D ternary displacement at ${T=0.2}$
with the second order JX scheme in ${x}$-dimension
and second order optimal variable relaxed scheme ${y}$-dimension on a ${40 \times 40}$ grid,
JX subcharacteristic speed in either dimensions = 5.8, ${\Delta t = 0.5 \frac{\Delta x}{5.8}}$
.}
\label{fig5.17}
\end{figure}

Figure (\ref {fig5.12}) demonstrates the increase in numerical diffusion for the JX scheme in the
presence of heterogeneity that results in high local velocities.
The second order variable relaxed
schemes produce solution profiles that are noticeably sharper than that of the second order JX schemes. Even
the first order variable relaxed schemes (figure (\ref{fig5.14}))
resolve solutions profiles sharper than the second order JX scheme.

The extreme smearing seen with the second order JX scheme is caused by the disparity between the JX subcharacteristic
speeds, which is dependent on the global maximum speeds, and the average speeds. This is particularly strong in the $y$-dimension for this example.
The maximum and average velocities in the $y$-dimension
are 2.87 and  0.35 and those in the $x$-dimension are about 5 and 2.4. The JX subcharacteristic is
governed by maximal velocity and the sub-characteristic condition as
\begin{eqnarray}
\frac{\left(\lambda _{\max }^{x} \right)^{2} }{\left(a^{x} \right)^{2} }
+\frac{\left(\lambda _{\max }^{y} \right)^{2} }{\left(a^{y} \right)^{2} } \le 1,
\nonumber
\end{eqnarray}
which imposes a much severe restriction on the possible values of $a^{x}$ and  $a^{y}$.  Only certain pairs of
$a^{x}$ and  $a^{y}$ obey the above condition.
Table (\ref{table5.2}) lists three possible pairs of $a^{x}$ and $a^{y}$: (i) the pair which has minimum possible $a^{x}$,
(ii) the pair with minimum possible $a^{y}$,
and (iii) the pair for which $a^{x} = a^{y}$.
The table shows that the sub-characteristic speed in ${y}$-dimension will be at least 10 times the average speed.
This causes the solution to be heavily smeared. Indeed,
if we apply a variable relaxed scheme in the ${y}$-dimension, while still retaining the JX discretization in
the ${x}$-dimension, we can observe an immediate improvement (see figure (\ref{fig5.17})).

\subsection {Single phase geometric optics problem}
Consider the 2x2 system
\begin{eqnarray}
\label{eq:eqn5.15}
\hspace{-0.2in}
&{\bf C}_{t}  + {\bf{F} \left(\bf{C} \right)}_{x} + {\bf{G} \left(\bf{C} \right)}_{y} = {\bf 0}, &
\\ \nonumber \\[-1.5ex]
\hspace{-0.7in}\textnormal{where}\; \;\;
&{\bf C}=\left[\begin{array}{c} {C_1} \\ {C_2} \end{array}\right] , \;
{\bf F}=\left[\begin{array}{c} {\frac {C_1^2} {\sqrt{C_1^2+ C_2^2}} } \\ {\frac{C_1C_2}{\sqrt{C_1^2 + C_2^2}} }
\end{array}\right] \;\;
\textnormal{and}\; \;
{\bf G}=\left[\begin{array}{c} {\frac {C_1C_2} {\sqrt{C_1^2+ C_2^2}} } \\ {\frac{C_2^2}{\sqrt{C_1^2 + C_2^2}} }
\end{array}\right] . &
\nonumber
\end{eqnarray}
This system was introduced by Engquist and Runborg \cite{er} in a study of multiphase modeling of geometric optics.
The system (\ref{eq:eqn5.15}) represents a single phase wave equation traveling through vacuum. Its solution
is a single ray of strength ${g(r,t) = \sqrt{C_1^2 + C_2^2}}$,  a distance ${r \equiv r(x,y) }$ at an angle
${\theta = arctan (\frac{C_2}{C_1})}$.

The Jacobians of this system are given by
\begin{eqnarray}
\label{eq:eqn5.16}
\hspace{0.2in}
{\bf F'({\bf C})} &=& {\bf R}\left[\begin{array}{cc} \textnormal{cos}\;\theta& - \textnormal{sin}\;\theta \\
0&  \textnormal{cos}\;\theta \end{array}\right]{\bf R^{-1}}, \; \;\;
{\bf G'({\bf C})} = {\bf R}\left[\begin{array}{cc}  \textnormal{sin}\;\theta&  \textnormal{cos}\;\theta \\
0&  \textnormal{sin}\;\theta \end{array}\right]{\bf R^{-1}},
\\ \nonumber \\[-1.5ex]
\hspace{-0.7in}\textnormal{where}\; \;\;
{\bf R} &=& \left[\begin{array}{cc}  \textnormal{cos}\;\theta&  -\textnormal{sin}\;\theta \\
 \textnormal{sin}\;\theta&  \textnormal{cos}\;\theta \end{array}\right] .
\nonumber
\end{eqnarray}
The Jacobians  ${\bf F'({\bf C})}$ and ${\bf G'({\bf C})}$, and any linear combination of the Jacobians, have an incomplete set of eigenvectors. So the system is weakly hyperbolic everwhere and hence forms an interesting test case.

The system (\ref{eq:eqn5.15}) is  solved over the rectangle ${ 0 \le x \le 1 ,\; 0 \le y \le 2 }$ with
initial conditions $C_1 = 0, C_2 = 0$. At time $t = 0$ the system is activated by a point source
located at $\left( -0.2, 1 \right)$. The exact solution to this problem
$ { g = \textnormal{max}\frac {\left( 0, t-r\right)^3}{r}}$ is used as a Dirichlet condition on all boundaries.

\begin{figure}[h!]
\begin{center}
\vspace {-0.3in}
\subfigure{
\includegraphics[bb=0mm 0mm 208mm 296mm, angle = 90, scale=0.3]{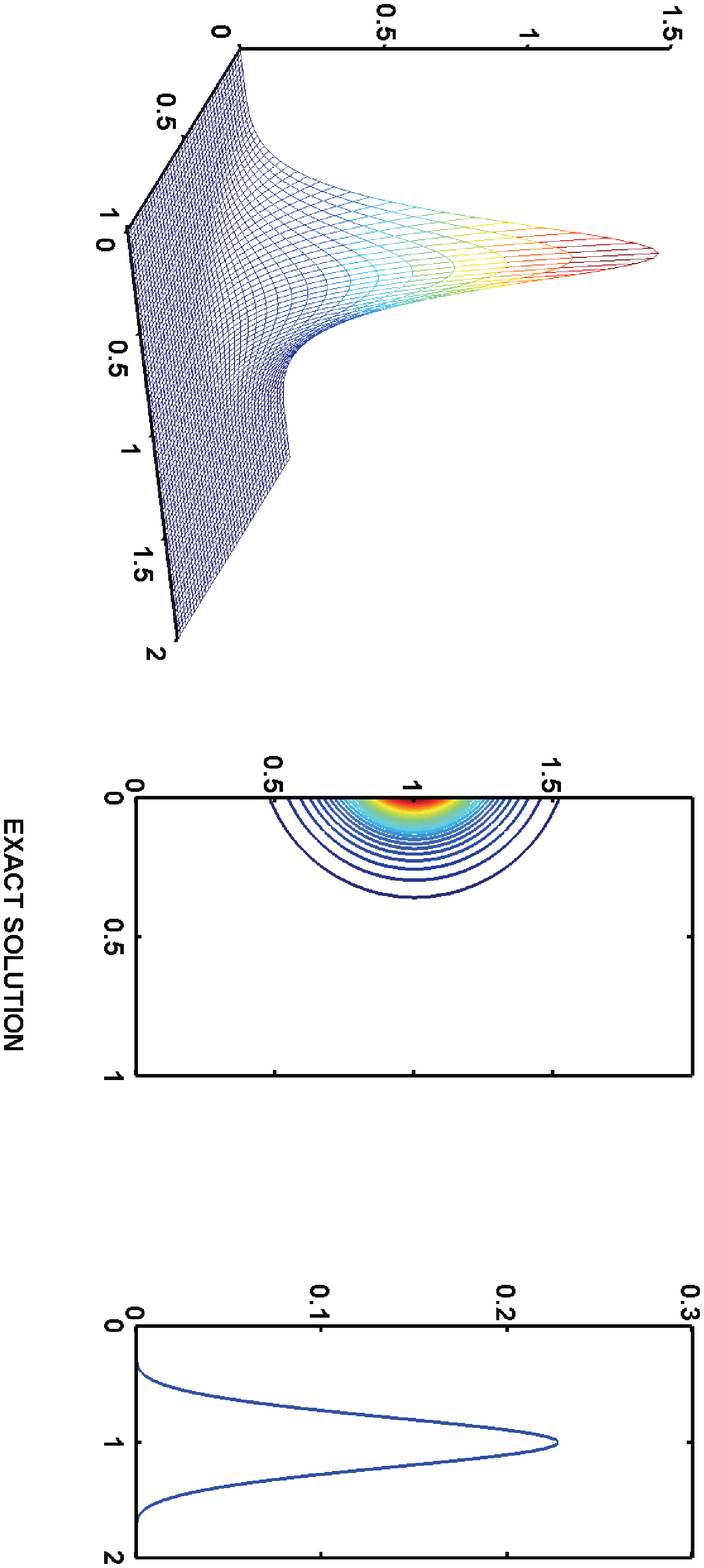}
}
\\
\vspace{-0.4in}
\subfigure{
\includegraphics[bb=0mm 0mm 208mm 296mm, angle = 90,scale=0.3]{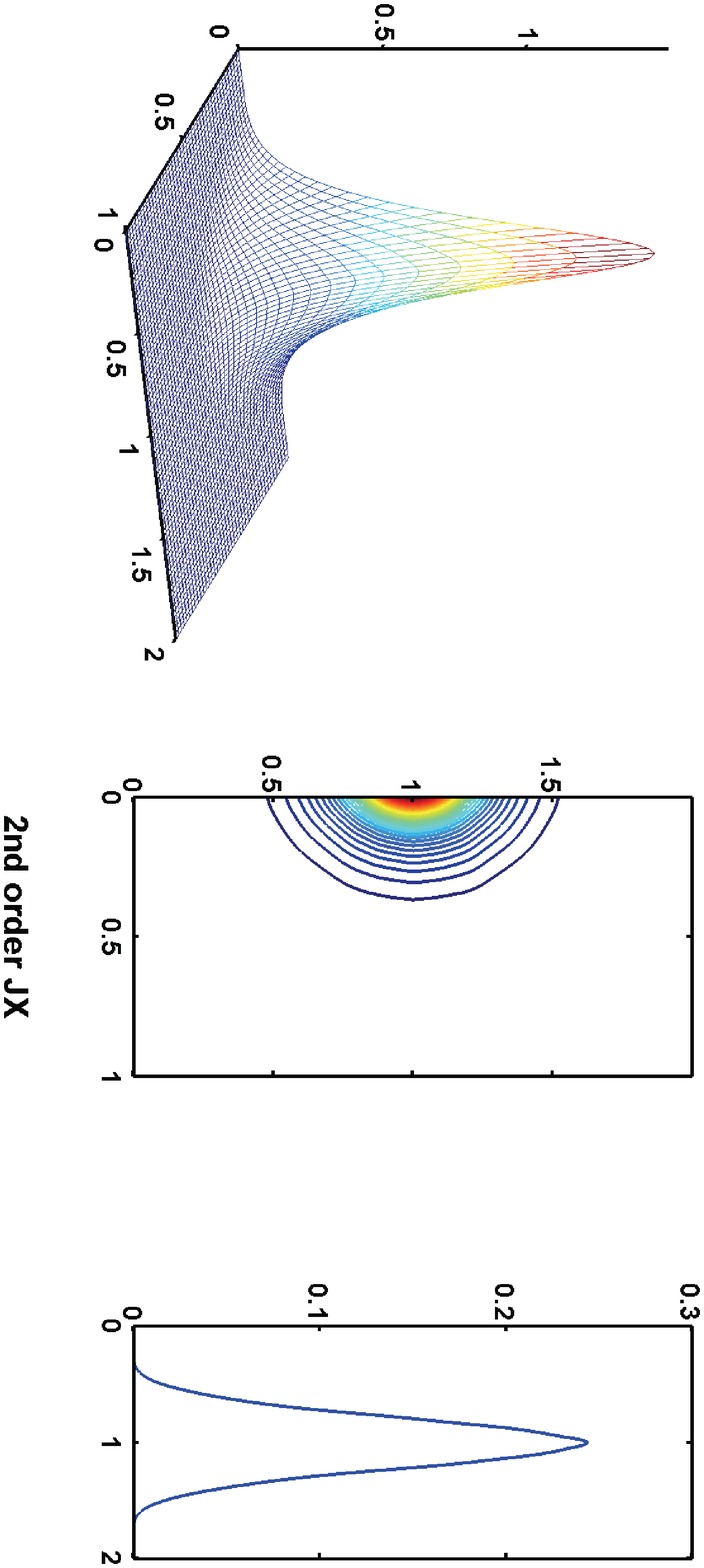}
}
\\
\vspace{-0.5in}
\subfigure{
\includegraphics[bb=0mm 0mm 208mm 296mm, angle = 90,scale=0.3]{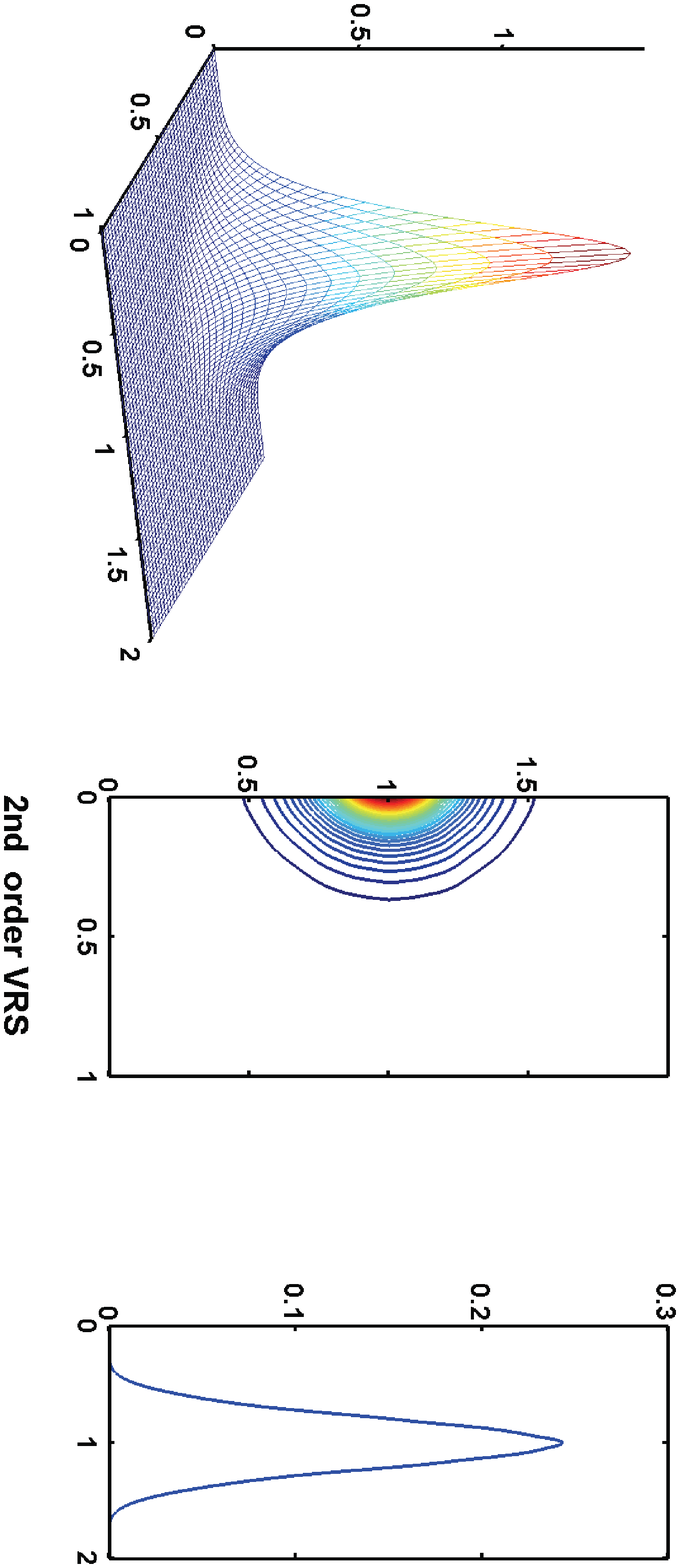}
}
\end{center}
\vspace{-0.1in}
\caption{Solution profiles for $2$D Engquist-Runborg problem ${T=0.85}$,
on a ${40 \times 80}$ grid, JX subcharacteristic speed = 1.4142, ${\Delta t = 0.5 \frac{\Delta x}{1.4142}}$.
The leftmost column shows the ray strengths, contour plots in the middle column and vertical cuts of solution
at $x=0.2$
in the rightmost column. The top row shows the exact solution, the middle row is that of the second order JX scheme and the last row
shows the result of the second order VRS scheme.}
\label{fig5.18}
\end{figure}

\begin{figure}[h!]
\begin{center}
\vspace {-0.3in}
\subfigure{
\includegraphics[bb=0mm 0mm 208mm 296mm, angle = 90, scale=0.3]{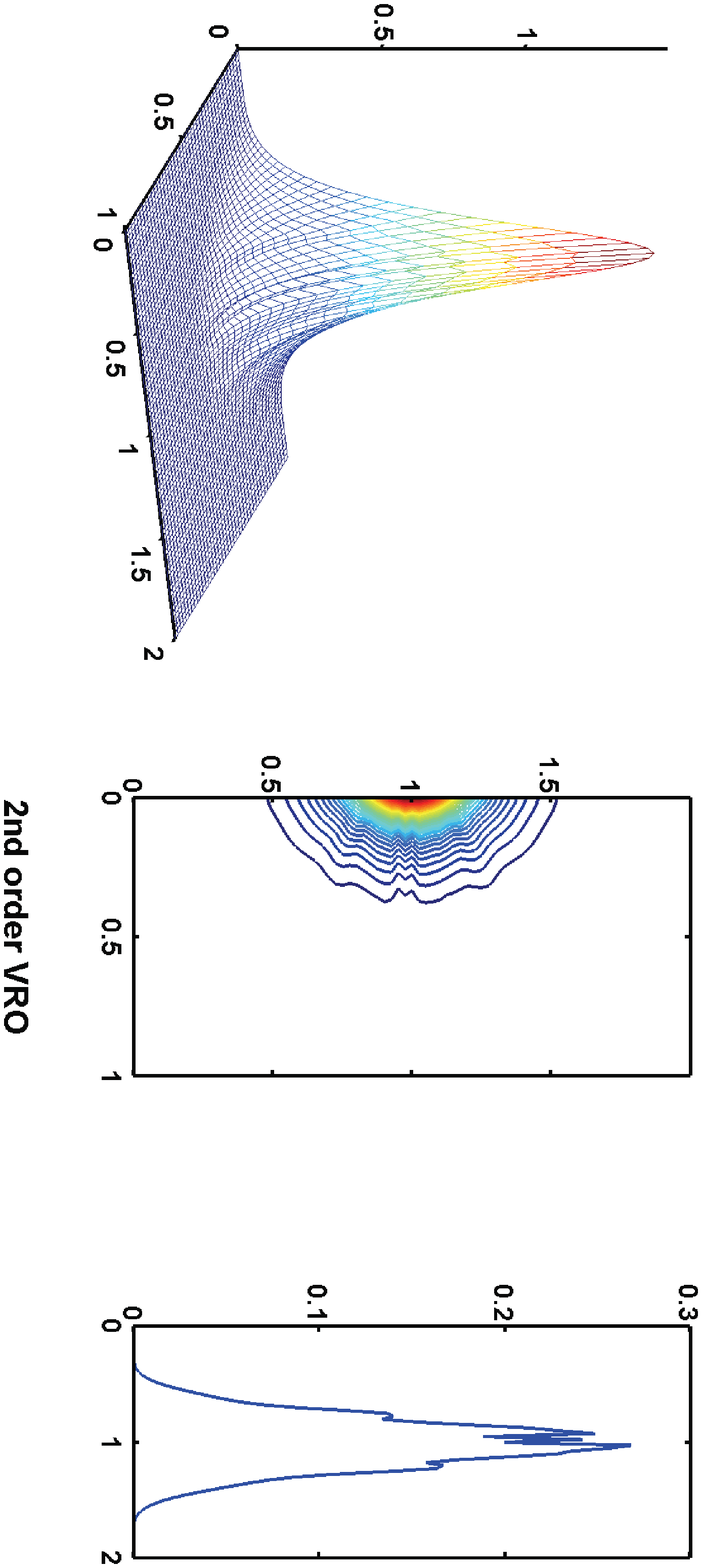}
}
\\
\vspace{-0.7in}
\subfigure{
\includegraphics[bb=0mm 0mm 208mm 296mm, angle = 90, scale=0.3]{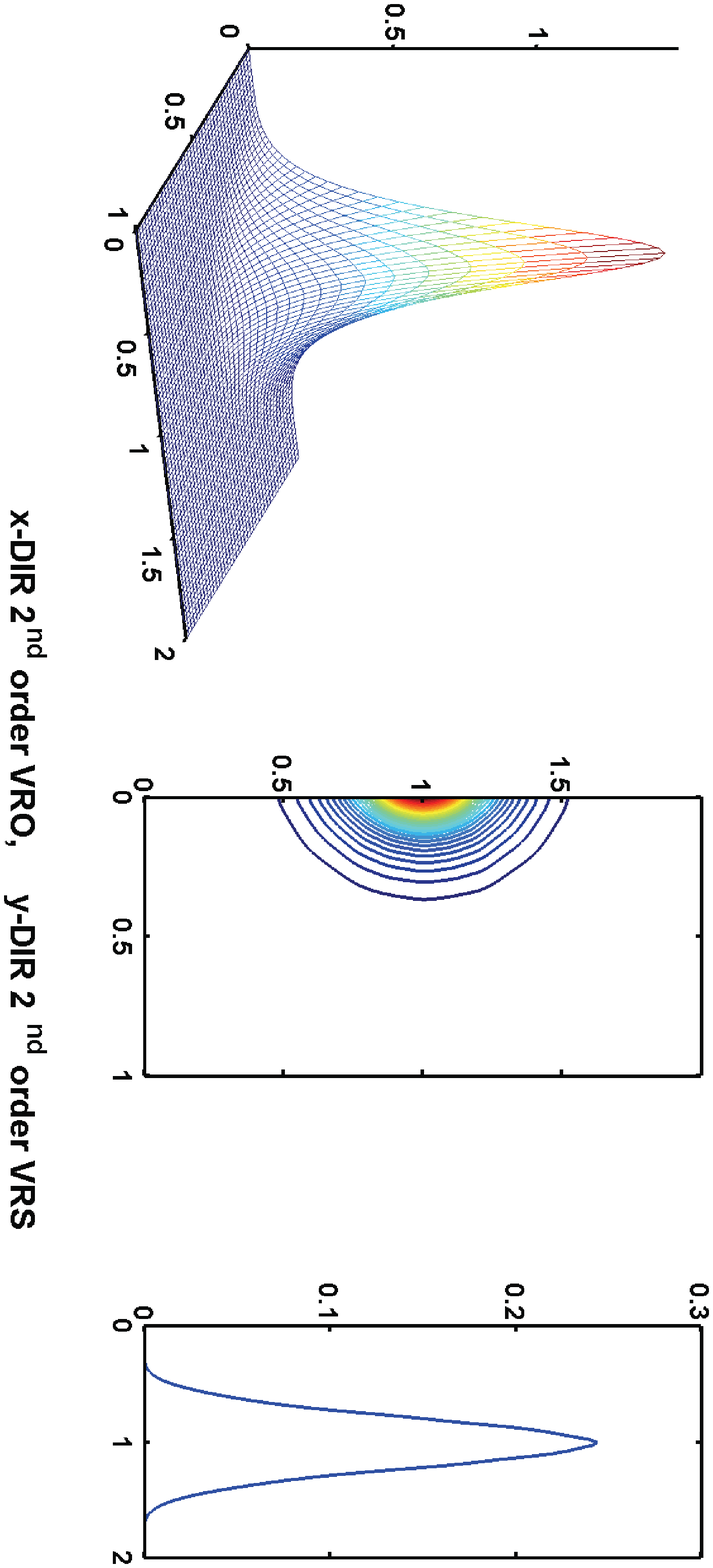}
}
\end{center}
\vspace{-0.1in}
\caption{Ray strengths, contour plots and vertical cuts of solution at  $x=0.2$,
for $2$D Engquist-Runborg problem ${T=0.85}$
with the second order VRO scheme on a ${40 \times 80}$ grid,
maximum subcharacteristic speed = 2,
${\Delta t = 0.5 \frac{\Delta x}{2}}$.}
\label{fig5.19}
\end{figure}

Both the second order JX and VRS schemes, which use symmetric subcharacteristics, produce resonably accurate profiles
(figure (\ref{fig5.18})).
However, the VRO scheme, which is similar to the traditional second order upwind scheme, leads to oscillations in the solution profile.
A similar phenomenon was observed for these upwind schemes in \cite{er}. The problem arises primarily in the
$y$-dimension. Indeed, if JX or VRS  scheme is applied in the $y$-dimension, while retaining the
VRO in $x$-dimension, there is an immediate improvement (figure (\ref{fig5.19})).
\newpage
\newpage
\section{Computational effort} In this section we explain in detail why  gas injection processes
require lesser phase equilibrium calculations in a second order relaxation framework than in a second order
central/central upwind framework.

Consider the scalar 1D semi-discrete KT scheme (equations 4.3,4.4 and 4.5 in \cite{kt}):
\begin {subequations}
\label{eqn6.1}
\begin{eqnarray}
\frac{d}{dt}{ C_{j}(t)}
= -\frac{{\mathcal H}_{j+{\frac{1}{2}}} (t) - {\mathcal H}_{j-{\frac{1}{2}}} (t)}
{\Delta x}
\label{eq:eqn6.1a}
\end{eqnarray}
where the numerical flux
\begin{eqnarray}
{\mathcal H}_{j+{\frac{1}{2}}} (t) :=
\frac{F(C_{j+{\frac{1}{2}} }^{+}(t)) + F(C_{j+{\frac{1}{2}} }^{-}(t))  }
{2}
- \frac{ a_{j+{\frac{1}{2}}}(t) }{2}
[C_{j+{\frac{1}{2}} }^{+}(t) - C_{j+{\frac{1}{2}} }^{-}(t)]
\label{eq:eqn6.1b}
\end{eqnarray}
where the intermediate values ${C_{j+{\frac{1}{2}} }^{\pm}}$ are given by
\begin{eqnarray}
C_{j+{\frac{1}{2}} }^{+}(t) :=
C_{j+1 }(t)  - \frac{\Delta x}{2}(C_x)_{j+1}(t), \;\;\;
C_{j+{\frac{1}{2}} }^{-}(t) :=
C_{j }(t)  + \frac{\Delta x}{2}(C_x)_{j}(t)
\label{eq:eqn6.1c}
\end{eqnarray}
\end {subequations}
For a N-cell grid, evaluating the numerical flux at the N+1 interfaces requires 2N flux evaluations at the points
\begin{eqnarray}
\nonumber
\{C_{1}(t)  + \frac{\Delta x}{2}(C_x)_{1}(t), \quad
C_{2}(t)  - \frac{\Delta x}{2}(C_x)_{2}(t), \quad
C_{2}(t)  + \frac{\Delta x}{2}(C_x)_{2}(t),  \quad
\cdots
\\
\nonumber
C_{j }(t)  + \frac{\Delta x}{2}(C_x)_{j}(t),  \quad
C_{j+1 }(t)  - \frac{\Delta x}{2}(C_x)_{j+1}(t),  \quad
\cdots
\\
\nonumber
C_{N}(t)  - \frac{\Delta x}{2}(C_x)_{N}(t), \quad
C_{N }(t)  + \frac{\Delta x}{2}(C_x)_{N}(t), \quad
C_{N+1 }(t)  - \frac{\Delta x}{2}(C_x)_{N+1}(t) \}
\end{eqnarray}
With the second order relaxed framework, the numerical flux at the N+1 interfaces (equations \ref{eqn2.9},  \ref{eqn3.11} )
requires only N+1 flux evaluations at points
\begin{eqnarray}
\nonumber
\{C_{1}(t), \quad
C_{2}(t), \quad
C_{3}(t), \quad
\cdots , \quad
C_{j }(t),  \quad
\cdots , \quad
C_{N-1}(t), \quad
C_{N}(t), \quad
C_{N+1}(t) \}
\end{eqnarray}
In general, for most systems, this difference in the number of flux computations adds only a very small overhead. However
with gas injection processes this leads to significant increase in computational effort. For example,
consider the 1-D ternary gas injection problem of equation(\ref{eq:eqn5.10}), whose fluxes are given by
\begin{eqnarray}
F_{1}\left(C_{1}, C_{2}\right) =c_{1V} f(S)+c_{1L} (1-f(S)),  \qquad
F_{2}\left(C_{1}, C_{2}\right) =c_{2V} f(S)+c_{2L} (1-f(S)),
\nonumber
\end{eqnarray}
where ${f(S)}$ is given by  equation(\ref{eq:eqn5.4}). Here, the fluxes ${F_1}$  and ${F_2}$ cannot be represented directly in terms of components
${C_1}$  and ${C_2}$, and are represented via the fractional flow curve ${f(S)}$. Evaluation of flux at composition points ${\{C_{1,2;j}(t)\} }$
or at the intermediate pair of composition points
\begin{eqnarray}
\nonumber
\{C_{1,2;j}(t)  - \frac{\Delta x}{2}(C_{1,2;x})_{j}(t), C_{1,2;j }(t)  + \frac{\Delta x}{2}(C_{1,2;x})_{j}(t)\}
\end{eqnarray}
necessitates phase equilibrium calculations in order to determine the saturation, and hence  ${f(S)}$ at those points.
For a 1D N-cell grid, each semi-discrete step in a central framework requires N-1 more flux evaluations (and hence phase equilibrium calculations) than the
semi-discrete step in a relaxation framework. When used with a 2-stage RK timestepping, this difference will double in 1D and
will be quadruple in 2D. A similar arguement applies for second order central upwind framework.

Note that for first order schemes, the intermediate pair of points
\begin{eqnarray}
\nonumber
\{C_{1,2;j }(t)  - \frac{\Delta x}{2}(C_{1,2;x})_{j}(t), C_{1,2;j }(t)  + \frac{\Delta x}{2}(C_{1,2;x})_{j}(t)\}
\end{eqnarray}
 collapse into one point${\{C_{1,2;j}(t)\} }$. Thus the computation effort will be same for first order central and relaxed schemes.
This has to be expected since, as we saw in section 3,  the first order central/central upwind schemes
are the same as the first order relaxed schemes.
\section{Discussion and conclusions} In this work we presented improved variable relaxation schemes for multidimensional
hyperbolic conservation laws. The motivation for our work is the weakly hyperbolic gas-injection displacements, the equations
of which are governed by strong nonlinear coupling and require costly thermodynamic equilibrium calculations every timestep.
For these problems, the traditional use of upwind schemes is problematic. Simulating these processes in central/central-upwind framework, which are
eigenstructure independent, requires more number of expensive thermodynamic equilibrium calculations in everytime step.
Jin and Xin's relaxation scheme, while providing a useful alternative to central schemes, still exhibits excessive
numerical diffusion in the presence of high contrasts in local velocities, as shown in our numerical experiments.
Our variable relaxation schemes retain the simplicity of the Jin-Xin relaxation schemes, but
improve the resolution significantly by using local subcharacteristic speeds, both in one and two spatial dimensions.
We solve for the relaxed systems themselves using traditional first and second order (TVD) upwind schemes.
We prove that the fully discrete one-dimensional schemes are monotone (first order) and  TVD (Total Variation Diminishing).
We proposed two types of subcharacteristic  speeds: optimal (VRO) and symmetric (VRS). Both work well for the gas-injection cases tested,
which have isolated points of weak hyperbolicity.
The optimal scheme, which adapts itself to become fully upwind in the presence of one-sided fluxes
can exhibit problems when the system to be solved is weakly hyperbolic everywhere in the domain.
Similar behavior was observed for traditional upwind schemes and because of the close relation of VRO with such schemes this behavior is not surprising.
The symmetric variable relaxed scheme VRS can be successfully used irrespective of the eigenstructure of
the problem. Both optimal and symmetric variable relaxation can be used on strongly hyperbolic problems,
with the advantage of avoiding characteristic decomposition or nonlinear Riemann solutions, and thus being faster.
\renewcommand{\theequation}{A-\arabic{equation}}
\setcounter{equation}{0}  
\section*{APPENDIX: Proofs for theorems 3.1 and 3.2}
${}$
\renewcommand{\thetheorem}{3.1}
\begin{theorem}
The first order, fully discrete,  variable relaxed scheme
\begin{eqnarray}
{C}_{j}^{n+1} = {C}_{j}^{n}
 - \frac{\Delta t}{\Delta x } \left({\mathcal F}_{j+{\frac{1}{2}}}^{n} - {\mathcal F}_{j-{\frac{1}{2}}}^{n} \right),
\label {eq:eqnC.1}
\end{eqnarray}
with symmetric speeds, where ${\mathcal F}_{j\pm{\frac{1}{2}}}^{n}$ is given by (\ref {eq:eqn3.13}),
is monotonic under the local subcharacteristic condition
\begin{eqnarray}
\nonumber
a_{j-{\frac{1}{2}} } \ge \left|F'\left(C\right)_{j-1}^{n}\right| \;\; \textnormal{and} \;\;
a_{j-{\frac{1}{2}}} \ge \left|F'\left(C\right)_{j}^{n}\right|,
\end{eqnarray}
and the time step restriction $\frac{\Delta t}{\Delta x} a_{\max } \le 1$, where $a_{\max }$ is the
maximum subcharacteristic speed.
The first order scheme (\ref {eq:eqn3.16})  with optimal speeds,
where ${\mathcal F}_{j\pm{\frac{1}{2}}}^{n}$ is given by (\ref {eq:eqn3.8b}),  is monotonic under
the local subcharacteristic condition
\begin{eqnarray}
\nonumber
a_{j-{\frac{1}{2}} }^{-} \le    \min \left( F'(C)_{j-1}^{n}, 0\right), \;\; \textnormal{and} \;\;
a_{j-{\frac{1}{2}} }^{+} \ge    \max \left( F'(C)_{j}^{n}, 0\right),
\end{eqnarray}
and the time step restriction $\frac{\Delta t}{\Delta x} a_{\max } \le \frac {1}{2}$.
\end {theorem}
\vspace{0.1in}
\begin {proof}
Representing the above update as $C_{j}^{n+1} ={\rm H} (C_{}^{n} ;j)$, the scheme is monotone if
$\frac{\partial {\rm H} (C_{}^{n} ;j)}{\partial C_{i}^{n} } \ge 0$ $\forall {\rm \; \; }i,j,C_{}^{n} $ (see \cite{lvq02}).

The update  $C_{j}^{n+1} ={\rm H} (C_{}^{n} ;j)$ depends only on variables from cells $j-1$, $j$, $j+1$.
So $\frac{\partial {\rm H} (C_{}^{n} ;j)}{\partial C_{i}^{n} } $ is nonzero only w.r.t
$C_{j-1}^{n} $, $C_{j+1}^{n} $, $C_{j}^{n} $.

\noindent Differential of update (\ref{eq:eqnC.1}) w.r.t $C_{j-1}^{n} $:
\begin{eqnarray}
\nonumber
\frac{\partial {\rm H} (C_{}^{n} ;i)}{\partial C_{j-1}^{n} }
&=& 0 -
\frac{\Delta t}{\Delta x} \left({\rm \; 0} + 0 - F'(C)_{j-1}^{n} +
\frac{\left(-a_{j-{\frac{1}{2}} }^{-} \right)}
{a_{j-{\frac{1}{2}} }^{+} + \left(-a_{j-{\frac{1}{2}} }^{-} \right)}
F'(C)_{j-1}^{n}
- \frac{ \left(-a_{j-{\frac{1}{2}} }^{-} \right)  a_{j-{\frac{1}{2}} }^{+} }
{a_{j-{\frac{1}{2}} }^{+} + \left(-a_{j-{\frac{1}{2}} }^{-} \right)}
\right)
\\
&=&\frac{\Delta t}{\Delta x} a_{j-{\frac{1}{2}} }^{+}
\frac{{\rm \; }F'(C)_{j-1}^{n} + \left(-a_{j-{\frac{1}{2}} }^{-} \right)}
{a_{j-{\frac{1}{2}} }^{+} + \left(-a_{j-{\frac{1}{2}} }^{-} \right)} .
\nonumber
\end{eqnarray}
\begin{eqnarray}
\nonumber
\hspace{-2.75in} \textnormal{If} \; \;
F'(C)_{j-1}^{n} \ge 0
\;\; \textnormal{then} \;\;
\frac{\partial {\rm H} (C_{}^{n} ;i)}{\partial C_{j-1}^{n} } \ge 0 .
\end{eqnarray}
\begin{eqnarray}
\nonumber
\hspace{-1in}
\textnormal{If} \; \;
F'(C)_{j-1}^{n} < 0
\;\; \textnormal{then} \;\;
\frac{\partial {\rm H} (C_{}^{n} ;i)}{\partial C_{j-1}^{n} } \ge 0
\;\; \textnormal{only if} \;\;
\left|a_{j-{\frac{1}{2}} }^{-} {\rm \; }\right| \ge \left| F'(C)_{j-1}^{n}\right|.
\end{eqnarray}
For symmetric case this requirement becomes,
\begin {subequations}
\label{eqC.2}
\begin{equation}
\label{eq:eqC.2a}
\left|a_{j-{\frac{1}{2}} }^{-} {\rm \; }\right|=a_{j-{\frac{1}{2}} }  \ge \left| F'(C)_{j-1}^{n}\right|,
\end{equation}
and for optimal case
\begin{equation}
\label{eq:eqC.2b}
a_{j-{\frac{1}{2}} }^{-} \le    \min \left( F'(C)_{j-1}^{n}, 0\right).
\end{equation}
\end {subequations}

\noindent Differential of update (\ref{eq:eqnC.1}) w.r.t $C_{j+1}^{n} $
\begin{eqnarray}
\nonumber
\frac{\partial {\rm H} (C_{}^{n} ;i)}{\partial C_{j+1}^{n} }
&=& 0 -
\frac{\Delta t}{\Delta x} \left( 0
+ \frac{\left(-a_{j+{\frac{1}{2}}}^{-} \right)}
{a_{j+{\frac{1}{2}} }^{+} +\left( - a_{j+{\frac{1}{2}} }^{-} \right)}
 \left[F'(C)_{j+1}^{n} - a_{j+{\frac{1}{2}} }^{+} \right] - 0 - 0
\right)
\\
&=& \frac{\Delta t}{\Delta x}
\left(-a_{j+{\frac{1}{2}} }^{-} \right)
\frac{a_{j+{\frac{1}{2}} }^{+} -F'(C)_{j+1}^{n} }{a_{j+{\frac{1}{2}} }^{+} +\left(-a_{j+{\frac{1}{2}} }^{-} \right)}.
\nonumber
\end{eqnarray}
\begin{eqnarray}
\nonumber
\hspace{-2.75in} \textnormal{If} \; \;
F'(C)_{j+1}^{n} \le 0
\;\; \textnormal{then} \;\;
\frac{\partial {\rm H} (C_{}^{n} ;i)}{\partial C_{j-1}^{n} } \ge 0 .
\end{eqnarray}
\begin{eqnarray}
\nonumber
\hspace{-1in}
\textnormal{If} \; \;
F'(C)_{j+1}^{n} > 0
\;\; \textnormal{then} \;\;
\frac{\partial {\rm H} (C_{}^{n} ;i)}{\partial C_{j-1}^{n} } \ge 0
\;\; \textnormal{only if} \;\;
a_{j+{\frac{1}{2}} }^{+}  \ge \left| F'(C)_{j+1}^{n}\right|.
\end{eqnarray}
For symmetric case this requirement becomes,
\begin {subequations}
\label{eqC.3}
\begin{equation}
\label{eq:eqC.3a}
a_{j+{\frac{1}{2}} }^{+} = a_{j+{\frac{1}{2}} }  \ge \left| F'(C)_{j+1}^{n}\right|,
\end{equation}
and for optimal case
\begin{equation}
\label{eq:eqC.3b}
a_{j+{\frac{1}{2}} }^{+} \ge    \max \left( F'(C)_{j+1}^{n}, 0\right).
\end{equation}
\end {subequations}
\noindent Differential of update  (\ref{eq:eqnC.1}) w.r.t $C_{j}^{n} $
\begin{eqnarray}
\label{eq:eqC.4}
\frac{\partial {\rm H} (C_{}^{n} ;i)}{\partial C_{j}^{n} }
&=& 1 -
\frac{\Delta t}{\Delta x} \left(\frac{\left|a_{j-{\frac{1}{2}} }^{-} \right| a_{j-{\frac{1}{2}} }^{+} }
{a_{j-{\frac{1}{2}} }^{+} + \left|a_{j-{\frac{1}{2}} }^{-} \right|}
+\frac{\left|a_{j+{\frac{1}{2}} }^{-} \right|a_{j+{\frac{1}{2}} }^{+} }
{a_{j+{\frac{1}{2}} }^{+} +\left|a_{j+{\frac{1}{2}} }^{-} \right|} \right)
\nonumber
\\
&-&
\frac{\Delta t}{\Delta x} F'(C)_{j}^{n} \left(\frac{a_{j-{\frac{1}{2}} }^{+} }
{a_{j-{\frac{1}{2}} }^{+} +\left|a_{j-{\frac{1}{2}} }^{-} \right|} -
\frac{\left|a_{j+{\frac{1}{2}} }^{-} \right|}{a_{j+{\frac{1}{2}} }^{+}
+\left|a_{j+{\frac{1}{2}} }^{-} \right|} \right)
\nonumber
\\
&=& 1- \frac{\Delta t}{\Delta x} a_{\max } \frac{1}{a_{\max}}
\left(\frac{\left|a_{j-{\frac{1}{2}} }^{-} \right|a_{j-{\frac{1}{2}} }^{+} }
{a_{j-{\frac{1}{2}} }^{+} +\left|a_{j-{\frac{1}{2}} }^{-} \right|} +
\frac{\left|a_{j+{\frac{1}{2}} }^{-} \right|a_{j+{\frac{1}{2}} }^{+} }
{a_{j+{\frac{1}{2}} }^{+} +\left|a_{j+{\frac{1}{2}} }^{-} \right|} \right)
\\
\nonumber
&-&
\frac{\Delta t}{\Delta x} a_{\max } \frac{F'(C)_{j}^{n} }{a_{\max }}
\left(\frac{a_{j-{\frac{1}{2}} }^{+}}{a_{j-{\frac{1}{2}} }^{+} +\left|a_{j-{\frac{1}{2}} }^{-} \right|} -
\frac{\left|a_{j+{\frac{1}{2}} }^{-} \right|}{a_{j+{\frac{1}{2}} }^{+} +\left|a_{j+{\frac{1}{2}} }^{-} \right|} \right).
\nonumber
\end{eqnarray}
\noindent In the first term of the RHS of the (\ref{eq:eqC.4}),
\begin{eqnarray}
\nonumber
\frac{\max \left(\left|a_{j-{\frac{1}{2}} }^{-} \right|, a_{j-{\frac{1}{2}} }^{+} \right)}{a_{\max } } \le 1,
\;\; \textnormal{and} \;\;
\frac{\min \left(\left|a_{j-{\frac{1}{2}} }^{-} \right|, a_{j-{\frac{1}{2}} }^{+} \right)}{a_{j-{\frac{1}{2}} }^{+}
+\left|a_{j-{\frac{1}{2}} }^{-} \right|} \le \frac{1}{2},
\end{eqnarray}
\noindent which implies
\begin{eqnarray}
\nonumber
\frac{1}{a_{\max } } \frac{\left|a_{j-{\frac{1}{2}} }^{-} \right|a_{j-{\frac{1}{2}} }^{+} }
{a_{j-{\frac{1}{2}} }^{+} +\left|a_{j-{\frac{1}{2}} }^{-} \right|}
=\frac{\max \left(\left|a_{j-{\frac{1}{2}} }^{-} \right|,a_{j-{\frac{1}{2}} }^{+} \right)}{a_{\max } }
\frac{\min \left(\left|a_{j-{\frac{1}{2}} }^{-} \right|,a_{j-{\frac{1}{2}} }^{+} \right)}{a_{j-{\frac{1}{2}} }^{+}
+\left|a_{j-{\frac{1}{2}} }^{-} \right|} \le \frac{1}{2}.
\end{eqnarray}
\noindent Similarly
\begin{eqnarray}
\nonumber
\frac{1}{a_{\max }}
\frac{\left|a_{j+{\frac{1}{2}} }^{-} \right|a_{j+{\frac{1}{2}} }^{+}}
{a_{j+{\frac{1}{2}} }^{+} +\left|a_{j+{\frac{1}{2}} }^{-} \right|}
=\frac{\max \left(\left|a_{j+{\frac{1}{2}} }^{-} \right|{\rm \; \; ,\; \; }a_{j+{\frac{1}{2}} }^{+} \right)}{a_{\max } }
\frac{\min \left(\left|a_{j+{\frac{1}{2}} }^{-} \right|{\rm \; \; ,\; \; }a_{j+{\frac{1}{2}} }^{+} \right){\rm \; \; }}
{a_{j+{\frac{1}{2}} }^{+} +\left|a_{j+{\frac{1}{2}} }^{-} \right|} \le \frac{1}{2}.
\end{eqnarray}
\noindent Combining the two inequalities we have,
\begin{eqnarray}
\nonumber
\frac{1}{a_{\max } } \left(\frac{\left|a_{j-{\frac{1}{2}} }^{-} \right|a_{j-{\frac{1}{2}} }^{+} {\rm \; \; }}
{a_{j-{\frac{1}{2}} }^{+} +\left|a_{j-{\frac{1}{2}} }^{-} \right|}
+\frac{\left|a_{j+{\frac{1}{2}} }^{-} \right|a_{j+{\frac{1}{2}} }^{+} {\rm \; }}{a_{j+{\frac{1}{2}} }^{+}
+\left|a_{j+{\frac{1}{2}} }^{-} \right|} \right)\le 1.
\end{eqnarray}
\noindent In the second term of the RHS of the  (\ref{eq:eqC.4}),
\begin{eqnarray}
\nonumber
\frac{a_{j-{\frac{1}{2}} }^{+} {\rm \; }}{a_{j-{\frac{1}{2}} }^{+}
+\left|a_{j-{\frac{1}{2}} }^{-} \right|} \le 1,
\;\;
\frac{\left|a_{j+{\frac{1}{2}} }^{-} \right|{\rm \; }}{a_{j+{\frac{1}{2}} }^{+}
+\left|a_{j+{\frac{1}{2}} }^{-} \right|} \le 1
\;\; \textnormal{so},\;\;
\left(\frac{a_{j-{\frac{1}{2}} }^{+} {\rm \; }}{a_{j-{\frac{1}{2}} }^{+} +\left|a_{j-{\frac{1}{2}} }^{-} \right|} {\rm \; }
-\frac{\left|a_{j+{\frac{1}{2}} }^{-} \right|{\rm \; }}{a_{j+{\frac{1}{2}} }^{+}
+\left|a_{j+{\frac{1}{2}} }^{-} \right|} \right)\le 1.
\end{eqnarray}
\noindent So, for a time-step restriction $\frac{\Delta t}{\Delta x} a_{\max } \le \frac{1}{2} $
\begin{eqnarray}
\nonumber
\frac{\Delta t}{\Delta x} a_{\max }
\frac{1}{a_{\max } } \left(\frac{\left|a_{j-{\frac{1}{2}} }^{-} \right|a_{j-{\frac{1}{2}} }^{+} {\rm \; \; }}
{a_{j-{\frac{1}{2}} }^{+} +\left|a_{j-{\frac{1}{2}} }^{-} \right|}
+\frac{\left|a_{j+{\frac{1}{2}} }^{-} \right|a_{j+{\frac{1}{2}} }^{+} {\rm \; }}
{a_{j+{\frac{1}{2}} }^{+} +\left|a_{j+{\frac{1}{2}} }^{-} \right|} \right) \le \frac{1}{2}
\;\; \textnormal{and}
\\ \nonumber
\left| \frac{\Delta t}{\Delta x} a_{\max } \frac{F'(C)_{j}^{n} }{a_{\max } } \left(\frac{a_{j-{\frac{1}{2}} }^{+} {\rm \; }}
{a_{j-{\frac{1}{2}} }^{+} +\left|a_{j-{\frac{1}{2}} }^{-} \right|} {\rm \; } -
\frac{\left|a_{j+{\frac{1}{2}} }^{-} \right|{\rm \; }}{a_{j+{\frac{1}{2}} }^{+}
+\left|a_{j+{\frac{1}{2}} }^{-} \right|} \right) \right| \le \frac{1}{2}.
\nonumber
\end{eqnarray}
\noindent Therefore,  $\frac{\partial {\rm H} (C_{}^{n} ;i)}{\partial C_{j}^{n} } \ge 0$ and the scheme
(\ref{eq:eqnC.1}) is monotone.

\noindent For the symmetric choice of speeds
 $a_{j-{\frac{1}{2}} }^{+} =\left|a_{j-{\frac{1}{2}} }^{-} \right| = a_{j-{\frac{1}{2}} }$ and
$a_{j+{\frac{1}{2}} }^{+} =\left|a_{j+{\frac{1}{2}} }^{-} \right| = a_{j+{\frac{1}{2}} } $
\begin{eqnarray}
\nonumber
\frac{\partial {\rm H} (C_{}^{n} ;i)}{\partial C_{j}^{n} } =
1-\frac{\Delta t}{\Delta x} a_{\max } \frac{1}{a_{\max } } \left(\frac{a_{j-{\frac{1}{2}} } {\rm \; \; }}{2}
+\frac{a_{j+{\frac{1}{2}} } {\rm \; }}{2} \right)\ge 0,
\end{eqnarray}
with a less restrictive time-step restriction $\frac{\Delta t}{\Delta x} a_{\max } \le 1$.
\end {proof}
\newpage
\renewcommand{\thetheorem}{3.2}
\begin{theorem}
The variable relaxed scheme with second order spatial discretization and forward Euler time stepping
\begin{equation}
\label{eq:eqC.5}
C_{j}^{n+1} =C_{j}^{n} -  \frac{1}{\Delta x}
\left( \mathcal{F}_{j+{\frac{1}{2}}}^{n} - \mathcal{F}_{j-{\frac{1}{2}}}^{n} \right)
- \frac{1}{\Delta x}
\left( \tilde{\mathcal {F}}_{j+{\frac{1}{2}}}^{n} -\tilde{\mathcal{F}}_{j-{\frac{1}{2}}}^{n} \right),
\end{equation}
with $\mathcal{F}_{j \pm {\frac{1}{2}}}^{n}$  given by (\ref{eq:eqn3.8b}) and
$\tilde{\mathcal {F}}_{j \pm {\frac{1}{2}}}^{n}$  given by (\ref{eq:eqn3.11b}),
is TVD  under the CFL condition $\frac{\Delta t}{\Delta x} a_{\max } \le \frac{1}{2} $, and the local subcharacteristic
condition
\begin{eqnarray}
a_{j-{\frac{1}{2}} }
\ge \left| \frac{ F\left(C\right)_{j}^{n} -  F\left(C\right)_{j-1}^{n} }{C_{j}^{n} -C_{j-1}^{n} } \right|,
\nonumber
\end{eqnarray}
for the symmetric case and,
\begin{eqnarray}
\nonumber
\hspace {0.4in}
a_{j-{\frac{1}{2}} }^{-} \le
\min \left( \frac{ F\left(C\right)_{j}^{n} -  F\left(C\right)_{j-1}^{n} }{C_{j}^{n} -C_{j-1}^{n} }, 0 \right),
a_{j-{\frac{1}{2}} }^{+} \ge \max
 \left( \frac{ F\left(C\right)_{j}^{n} -  F\left(C\right)_{j-1}^{n} }{C_{j}^{n} -C_{j-1}^{n} }, 0 \right),
\end{eqnarray}
for the optimal case.
\end{theorem}
\vspace{0.1in}
\begin{proof}
The proof is done along the lines of Harten's Theorem \cite{h83}, which states that, a scheme that
is in the form
\begin{equation}
\label{eq:eqC.6}
C_{j}^{n+1} =C_{j}^{n} - \kappa {1}_{j-1} \left(C_{j}^{n} -C_{j-1}^{n} \right)
+ \kappa {2}_{j} \left( C_{j+1}^{n} - C_{j}^{n} \right)
\end{equation}
is TVD, if $\kappa {1}_{j-1} \ge 0$, $\kappa {2}_{j} \ge 0$  and $\kappa {1}_{j} +\kappa {2}_{j} \le 1$, $\forall {j}$. \\
Equation  (\ref{eq:eqC.5}) can be rewritten in the form  (\ref{eq:eqC.6}) by setting
\begin{eqnarray}
&\kappa {1}_{j-1} = \hspace{4.625in}
\nonumber \\
&\frac{\Delta t}{\Delta x}
\frac
{\left( -a_{j-{\frac{1}{2}}}^{-} + \frac{F_{j}^{n} -F_{j-1}^{n}}{C_{j}^{n} -C_{j-1}^{n}} \right)}
{ a_{j-{\frac{1}{2}}}^{+} - a_{j-{\frac{1}{2}}}^{-}}
\left[
a_{j-{\frac{1}{2}}}^{+} \left(1- \frac{\phi \left(\theta _{j-{\frac{1}{2}}}^{+} \right)}{2} \right)
+a_{j+{\frac{1}{2}}}^{+}
\frac
{ \left(1+a_{j+{\frac{1}{2}}}^{+} a_{j-{\frac{1}{2}}}^{+} \right) }
{ \left(1+a_{j+{\frac{1}{2}}}^{+} a_{j+{\frac{1}{2}}}^{+} \right) }
\frac{ \phi \left( \theta _{j+{\frac{1}{2}}}^{+} \right) }{ 2\theta _{j+{\frac{1}{2}}}^{+} }
\right]&
\nonumber \\ \nonumber \\
&\kappa {2}_{j} =  \hspace{5in}
\nonumber \\
&\frac{\Delta t}{\Delta x}
\frac
{\left(a_{j+{\frac{1}{2}}}^{+} - \frac{F_{j+1}^{n} -F_{j}^{n}}{C_{j+1}^{n} -C_{j}^{n} } \right) }
{ a_{j+{\frac{1}{2}}}^{+} - a_{j+{\frac{1}{2}}}^{-} }
\left[
-a_{j+{\frac{1}{2}} }^{-}\left( 1- \frac{\phi \left(\theta _{j+{\frac{1}{2}}}^{-} \right)}{2}\right)
-a_{j-{\frac{1}{2}} }^{-}
\frac
{\left(1+a_{j+{\frac{1}{2}}}^{-} a_{j-{\frac{1}{2}}}^{-} \right)}
{\left(1+a_{j-{\frac{1}{2}}}^{-} a_{j-{\frac{1}{2}}}^{-} \right)}
\frac{\phi \left(\theta _{j-{\frac{1}{2}}}^{-} \right)}{2\theta _{j-{\frac{1}{2}}}^{-} }
\right] &
\nonumber
\end{eqnarray}
\noindent
\textbf {To prove that  $\kappa {1}_{j-1} \ge 0\;\;$,  $\kappa {2}_{j} \ge 0\;\;$  $\forall {j}$:}
\\In the expression for $\kappa {1}_{j-1}$, the following inequalities hold
\begin{eqnarray}
\frac{1}{a_{j-{\frac{1}{2}}}^{+} - a_{j-{\frac{1}{2}}}^{-}} > 0, \;
\left[
a_{j-{\frac{1}{2}}}^{+} \left(1- \frac{\phi \left(\theta _{j-{\frac{1}{2}}}^{+} \right)}{2} \right)
+a_{j+{\frac{1}{2}}}^{+}
\frac
{1+a_{j+{\frac{1}{2}}}^{+} a_{j-{\frac{1}{2}}}^{+}}
{1+a_{j+{\frac{1}{2}}}^{+} a_{j+{\frac{1}{2}}}^{+}}
\frac
{\phi \left(\theta _{j+{\frac{1}{2}}}^{+} \right) }
{2\theta _{j+{\frac{1}{2}}}^{+}}
\right]>0
\nonumber
\end{eqnarray}
by the choice of subcharacteristic speeds and because van Leer limiter obeys the bound
$\left(1-\frac{1}{2} \phi \left(\theta _{j} \right)\right)\ge 0$   and
$\frac{\phi \left(\theta _{j} \right)}{2\theta _{j}} \ge 0$.
\begin{eqnarray}
\nonumber
\hspace{-2.6in} \textnormal{If} \; \;
\frac{F_{j}^{n} -F_{j-1}^{n}}{C_{j}^{n} - C_{j-1}^{n}} \ge 0
\;\; \textnormal{then} \;\;
\kappa {1}_{j-1} \ge 0, \; \forall {j}.
\end{eqnarray}
\begin{eqnarray}
\nonumber
\hspace{-0.8in}
\textnormal{If} \; \;
\frac{F_{j}^{n} -F_{j-1}^{n}}{C_{j}^{n} - C_{j-1}^{n}} < 0
\;\; \textnormal{then} \;\;
\kappa {1}_{j-1} \ge 0, \; \forall {j}
\;\; \textnormal{only if} \;\;
\left|a_{j-{\frac{1}{2}}}^{-} \right|\ge \left|\frac{F_{j}^{n} -F_{j-1}^{n} }{C_{j}^{n} -C_{j-1}^{n} } \right|.
\end{eqnarray}
For symmetric case this requirement becomes,
\begin {subequations}
\label{eqC.7}
\begin{equation}
\label{eq:eqC.7a}
\left|a_{j-{\frac{1}{2}} }^{-} {\rm \; }\right|=a_{j-{\frac{1}{2}} }  \ge
\left| \frac{F_{j}^{n} -F_{j-1}^{n} }{C_{j}^{n} -C_{j-1}^{n} }  \right| \;\; \forall {j},
\end{equation}
and for optimal case
\begin{equation}
\label{eq:eqC.7b}
a_{j-{\frac{1}{2}} }^{-} \le    \min \left( \frac{F_{j}^{n} -F_{j-1}^{n} }{C_{j}^{n} -C_{j-1}^{n} }, 0\right)
\;\; \forall {j}.
\end{equation}
\end {subequations}
By similar reasoning, we can see that
\begin{eqnarray}
\nonumber
\hspace{-2.75in} \textnormal{if} \; \;
\frac{F_{j+1}^{n} -F_{j}^{n}}{C_{j+1}^{n} -C_{j}^{n} } \le 0
\;\; \textnormal{then} \;\;
\kappa {2}_{j} \ge 0,
\end{eqnarray}
\begin{eqnarray}
\nonumber
\hspace{-1in}
\textnormal{ and if} \; \;
\frac{F_{j+1}^{n} -F_{j}^{n}}{C_{j+1}^{n} -C_{j}^{n} } > 0
\;\; \textnormal{then} \;\;
\kappa {2}_{j} \ge 0
\;\; \textnormal{only if} \;\;
a_{j+{\frac{1}{2}} }^{+}  \ge \left|\frac{F_{j+1}^{n} -F_{j}^{n}}{C_{j+1}^{n} -C_{j}^{n} }\right|.
\end{eqnarray}
For symmetric case this requirement becomes,
\begin {subequations}
\label{eqC.8}
\begin{equation}
\label{eq:eqC.8a}
a_{j+{\frac{1}{2}} }^{+} = a_{j+{\frac{1}{2}} }  \ge \left| \frac{F_{j+1}^{n} -F_{j}^{n}}{C_{j+1}^{n} -C_{j}^{n} }\right|
\;\forall {j},
\end{equation}
and for optimal case
\begin{equation}
\label{eq:eqC.8b}
a_{j+{\frac{1}{2}} }^{+} \ge    \max \left( \frac{F_{j+1}^{n} -F_{j}^{n}}{C_{j+1}^{n} -C_{j}^{n} }, 0\right)
\; \forall {j}.
\end{equation}
\end {subequations}
To prove that $\; \kappa {1}_{j} +\kappa {2}_{j} \le 1 \;\;$, $\forall {j}$, consider,
\begin{eqnarray}
&\kappa {1}_{j}& +\kappa {2}_{j} \hspace{8in}
\nonumber \\
&=& \frac {\Delta t}{\Delta x}
\frac
{ \left( -a_{j+{\frac{1}{2}}}^{-} + \frac{F_{j+1}^{n}-F_{j}^{n}}{C_{j+1}^{n} -C_{j}^{n}} \right)}
{ a_{j+{\frac{1}{2}}}^{+} - a_{j+{\frac{1}{2}}}^{-} }
\left[
a_{j+{\frac{1}{2}}}^{+} - a_{j+{\frac{1}{2}}}^{+}\frac{\phi \left(\theta _{j+{\frac{1}{2}} }^{+} \right)}{2}
+ a_{j+{\frac{3}{2}}}^{+}
\frac
{\left(1+a_{j+{\frac{3}{2}}}^{+} a_{j+{\frac{1}{2}}}^{+} \right)}
{\left(1+a_{j+{\frac{3}{2}} }^{+} a_{j+{\frac{3}{2}} }^{+} \right)}
\frac
{\phi \left(\theta _{j+{\frac{3}{2}} }^{+} \right)} {2\theta _{j+{\frac{3}{2}}}^{+} }
\right]
\nonumber \\
&+& \frac{\Delta t}{\Delta x}
\frac
{ \left(a_{j+{\frac{1}{2}} }^{+} -\frac{F_{j+1}^{n} -F_{j}^{n}}{C_{j+1}^{n} -C_{j}^{n}} \right)}
{a_{j+{\frac{1}{2}} }^{+} - a_{j+{\frac{1}{2}} }^{-} }
\left[
- a_{j+{\frac{1}{2}}}^{-} + a_{j+{\frac{1}{2}}}^{-} \frac{\phi \left(\theta _{j+{\frac{1}{2}}}^{-} \right)}{2}
-a_{j-{\frac{1}{2}} }^{-} \frac
{\left(1+ a_{j+{\frac{1}{2}}}^{-} a_{j-{\frac{1}{2}}}^{-}\right)}
{\left(1+ a_{j-{\frac{1}{2}} }^{-} a_{j-{\frac{1}{2}} }^{-} \right)}
\frac{\phi \left(\theta _{j-{\frac{1}{2}} }^{-} \right)}{2\theta _{j-{\frac{1}{2}} }^{-} }
\right]
\nonumber
\end {eqnarray}
\begin{eqnarray}
&=&\frac {\Delta t}{\Delta x}
\frac
{ \left( \left | a_{j+{\frac{1}{2}}}^{-} \right | + \frac{F_{j+1}^{n}-F_{j}^{n}}{C_{j+1}^{n} -C_{j}^{n}} \right)}
{ a_{j+{\frac{1}{2}}}^{+} + \left | a_{j+{\frac{1}{2}}}^{-} \right | }
\left[
a_{j+{\frac{1}{2}}}^{+}
 +  \frac{1}{2}
\frac
{\left(a_{j+{\frac{3}{2}}}^{+} \right)^{2} a_{j+{\frac{1}{2}}}^{+}
\left(
\frac {\phi \left(\theta _{j+{\frac{3}{2}}}^{+}\right)}{\theta _{j+{\frac{3}{2}}}^{+}}
-\phi \left(\theta _{j+{\frac{1}{2}}}^{+} \right)
\right)       }
{\left(1+a_{j+{\frac{3}{2}} }^{+} a_{j+{\frac{3}{2}} }^{+} \right)}
\hspace{1in} \right .  \nonumber \\
& +& \left .  \frac{1}{2}
\frac
{
a_{j+{\frac{3}{2}}}^{+}
\left(\frac{\phi \left(\theta _{j+{\frac{3}{2}}}^{+} \right)}{\theta _{j+{\frac{3}{2}}}^{+} } \right)
 -a_{j+{\frac{1}{2}}}^{+}
\phi \left(\theta _{j+{\frac{1}{2}} }^{+} \right)}
{\left(1+a_{j+{\frac{3}{2}} }^{+} a_{j+{\frac{3}{2}} }^{+} \right)}
\right]
\nonumber \\
&+& \frac{\Delta t}{\Delta x}
\frac
{ \left(a_{j+{\frac{1}{2}} }^{+} -\frac{F_{j+1}^{n} -F_{j}^{n}}{C_{j+1}^{n} -C_{j}^{n}} \right)}
{a_{j+{\frac{1}{2}} }^{+} + \left | a_{j+{\frac{1}{2}}}^{-} \right |}
\left[
\left | a_{j+{\frac{1}{2}}}^{-} \right |
+\frac{1}{2}
\frac
{\left | a_{j+{\frac{1}{2}}}^{-} \right | \left(a_{j-{\frac{1}{2}}}^{-} \right)^{2}
\left(
\frac{\phi \left(\theta _{j-{\frac{1}{2}}}^{-} \right)}{\theta _{j-{\frac{1}{2}}}^{-}}
-\phi \left(\theta _{j+{\frac{1}{2}}}^{-} \right)
\right)}
{\left(1+ a_{j-{\frac{1}{2}} }^{-} a_{j-{\frac{1}{2}} }^{-} \right)}
\hspace{1in} \right . \nonumber \\
&+& \left . \frac{1}{2}
\frac
{
\left | a_{j-{\frac{1}{2}}}^{-} \right |
\left(\frac{\phi \left(\theta _{j-{\frac{1}{2}}}^{-} \right)}{\theta _{j-{\frac{1}{2}}}^{-}} \right)
- \left | a_{j+{\frac{1}{2}}}^{-} \right |
\phi \left(\theta _{j+{\frac{1}{2}}}^{-}
\right) }
{\left(1+ a_{j-{\frac{1}{2}} }^{-} a_{j-{\frac{1}{2}} }^{-} \right)}
\right]
\nonumber
\end {eqnarray}
\begin{eqnarray}
& \le &
\frac {\Delta t}{\Delta x}
\frac
{ \left( \left | a_{j+{\frac{1}{2}}}^{-} \right | + \frac{F_{j+1}^{n}-F_{j}^{n}}{C_{j+1}^{n} -C_{j}^{n}} \right)}
{ a_{j+{\frac{1}{2}}}^{+} + \left | a_{j+{\frac{1}{2}}}^{-} \right | }
\left[
a_{j+{\frac{1}{2}}}^{+}
 +  \frac{1}{2}
\frac
{\left(a_{j+{\frac{3}{2}}}^{+} \right)^{2} a_{j+{\frac{1}{2}}}^{+}
\max
\left(
\frac {\phi \left(\theta _{j+{\frac{3}{2}}}^{+}\right)}{\theta _{j+{\frac{3}{2}}}^{+}}
-\phi \left(\theta _{j+{\frac{1}{2}}}^{+} \right)
\right)       }
{\left(1+a_{j+{\frac{3}{2}} }^{+} a_{j+{\frac{3}{2}} }^{+} \right)}
\hspace{1in} \right.  \nonumber \\
& +& \left.  \frac{1}{2}
\frac
{
a_{j+{\frac{3}{2}}}^{+}
\max
\left(\frac{\phi \left(\theta _{j+{\frac{3}{2}}}^{+} \right)}{\theta _{j+{\frac{3}{2}}}^{+} } \right)
 -a_{j+{\frac{1}{2}}}^{+}
\min
\phi \left(\theta _{j+{\frac{1}{2}} }^{+} \right)}
{\left(1+a_{j+{\frac{3}{2}} }^{+} a_{j+{\frac{3}{2}} }^{+} \right)}
\right]
\nonumber \\
&+& \frac{\Delta t}{\Delta x}
\frac
{ \left(a_{j+{\frac{1}{2}} }^{+} -\frac{F_{j+1}^{n} -F_{j}^{n}}{C_{j+1}^{n} -C_{j}^{n}} \right)}
{a_{j+{\frac{1}{2}} }^{+} + \left | a_{j+{\frac{1}{2}}}^{-} \right |}
\left[
\left |  a_{j+{\frac{1}{2}}}^{-} \right |
+\frac{1}{2}
\frac
{\left |a_{j+{\frac{1}{2}}}^{-} \right | \left(a_{j-{\frac{1}{2}}}^{-} \right)^{2}
\max
\left(
\frac{\phi \left(\theta _{j-{\frac{1}{2}}}^{-} \right)}{\theta _{j-{\frac{1}{2}}}^{-}}
-\phi \left(\theta _{j+{\frac{1}{2}}}^{-} \right)
\right)}
{\left(1+ a_{j-{\frac{1}{2}} }^{-} a_{j-{\frac{1}{2}} }^{-} \right)}
\hspace{1in} \right. \nonumber \\
&+& \left.\frac{1}{2}
\frac
{
\left |a_{j-{\frac{1}{2}}}^{-} \right |
\max
\left(\frac{\phi \left(\theta _{j-{\frac{1}{2}}}^{-} \right)}{\theta _{j-{\frac{1}{2}}}^{-}} \right)
- \left |a_{j+{\frac{1}{2}}}^{-} \right |
\min
\phi \left(\theta _{j+{\frac{1}{2}}}^{-}
\right) }
{\left(1+ a_{j-{\frac{1}{2}} }^{-} a_{j-{\frac{1}{2}} }^{-} \right)}
\right]
\nonumber
\end{eqnarray}
For the van Leer limiter,
$\; \max \left[\frac{\phi \left(\theta _{j} \right)}{\theta _{j}}
-\phi \left(\theta _{j\pm 1} \right)\right] = 2 \; $,
$\; \max \left[\frac{\phi \left(\theta _{j} \right)}{\theta _{j}} \right] = 2\;$,
and $\; \min \phi \left(\theta _{j} \right) = 0 \; $. Substituting these in the above expression,
\begin{eqnarray}
\kappa {1}_{j}&+&\kappa {2}_{j} \hspace{6in} \nonumber \\
&\le&
\frac{\Delta t}{\Delta x} \frac
{\left | a_{j+{\frac{1}{2}} }^{-} \right |  + \frac{F_{j+1}^{n} -F_{j}^{n} }{C_{j+1}^{n} -C_{j}^{n}}}
{a_{j+{\frac{1}{2}}}^{+} + \left | a_{j+{\frac{1}{2}} }^{-} \right |}
\left[
a_{j+{\frac{1}{2}} }^{+} + \frac{1}{2} \frac{2\left(a_{j+{\frac{3}{2}} }^{+} \right)^{2} a_{j+{\frac{1}{2}} }^{+}
+2a_{j+{\frac{3}{2}} }^{+} }{1+\left(a_{j+{\frac{3}{2}} }^{+} \right)^{2} }
 \right]
\nonumber \\
&+&\frac{\Delta t}{\Delta x}
\frac{{\rm \; }a_{j+{\frac{1}{2}} }^{+} -\frac{F_{j+1}^{n} -F_{j}^{n} }{C_{j+1}^{n} -C_{j}^{n} } }
{a_{j+{\frac{1}{2}} }^{+} + \left |  a_{j+{\frac{1}{2}} }^{-} \right | }
\left[
\left | a_{j+{\frac{1}{2}} }^{-}  \right | +\frac{1}{2}
\frac
{2 \left | a_{j+{\frac{1}{2}}}^{-} \right | \left(a_{j-{\frac{1}{2}}}^{-} \right)^{2}
+ 2\left | a_{j-{\frac{1}{2}}}^{-} \right | }
{1+\left(a_{j-{\frac{1}{2}} }^{-} \right)^{2} }
 \right]
\nonumber
\\
&\le& \frac{\Delta t}{\Delta x} \frac
{\left | a_{j+{\frac{1}{2}}}^{-} \right | +\frac{F_{j+1}^{n} -F_{j}^{n} }{C_{j+1}^{n} -C_{j}^{n} }}
{a_{j+{\frac{1}{2}} }^{+} + \left | a_{j+{\frac{1}{2}} }^{-} \right |}
\left[
a_{j+{\frac{1}{2}} }^{+}
+\frac{1}{2}
\frac{2\left(a_{j+{\frac{3}{2}}}^{+} \right)^{2} a_{\max } + 2a_{\max }}
{1+\left(a_{j+{\frac{3}{2}} }^{+} \right)^{2} }
\right]
\nonumber \\
&+&\frac{\Delta t}{\Delta x}
\frac
{a_{j+{\frac{1}{2}}}^{+} -\frac{F_{j+1}^{n} -F_{j}^{n} }{C_{j+1}^{n} -C_{j}^{n}}}
{a_{j+{\frac{1}{2}}}^{+} + \left | a_{j+{\frac{1}{2}}}^{-}  \right |}
\left[
\left | a_{j+{\frac{1}{2}}}^{-} \right |
+\frac{1}{2}
\frac
{2a_{\max } \left(a_{j-{\frac{1}{2}} }^{-} \right)^{2} +2a_{\max } }
{1+\left(a_{j-{\frac{1}{2}} }^{-} \right)^{2} }
 \right]
\nonumber
\end {eqnarray}
\begin{eqnarray}
&\le& \frac{\Delta t}{\Delta x} \frac
{\left | a_{j+{\frac{1}{2}}}^{-} \right | +\frac{F_{j+1}^{n} -F_{j}^{n} }{C_{j+1}^{n} -C_{j}^{n} }}
{a_{j+{\frac{1}{2}} }^{+} + \left | a_{j+{\frac{1}{2}} }^{-} \right |}
a_{\max } \left[
\frac {a_{j+{\frac{1}{2}} }^{+} } {a_{\max } }
+ 1
\right]
\nonumber \\
&+&\frac{\Delta t}{\Delta x}
\frac
{a_{j+{\frac{1}{2}}}^{+} -\frac{F_{j+1}^{n} -F_{j}^{n} }{C_{j+1}^{n} -C_{j}^{n}}}
{a_{j+{\frac{1}{2}}}^{+} + \left | a_{j+{\frac{1}{2}}}^{-}  \right |}
 a_{\max }  \left[
\frac {\left | a_{j+{\frac{1}{2}}}^{-} \right | } {a_{\max }}
+ 1
 \right]
\nonumber
\\
&\le& \frac{\Delta t}{\Delta x} \frac
{\left | a_{j+{\frac{1}{2}}}^{-} \right | +\frac{F_{j+1}^{n} -F_{j}^{n} }{C_{j+1}^{n} -C_{j}^{n} }}
{a_{j+{\frac{1}{2}} }^{+} + \left | a_{j+{\frac{1}{2}} }^{-} \right |}
a_{\max } \left[
1 + 1
\right]
\nonumber \\
&+&\frac{\Delta t}{\Delta x}
\frac
{a_{j+{\frac{1}{2}}}^{+} -\frac{F_{j+1}^{n} -F_{j}^{n} }{C_{j+1}^{n} -C_{j}^{n}}}
{a_{j+{\frac{1}{2}}}^{+} + \left | a_{j+{\frac{1}{2}}}^{-}  \right |}
 a_{\max }  \left[
1 + 1
 \right]
\nonumber \\
&\le & 2 \frac{\Delta t}{\Delta x} a_{\max }
\nonumber
\end{eqnarray}
\qquad \end{proof}


\end{document}